\documentclass[10pt]{article}
\usepackage[T1]{fontenc}
\usepackage[utf8]{inputenc}
\usepackage{amsmath,amsfonts,amsthm,mathrsfs,amssymb}
\usepackage{cite}
\usepackage{multirow}
\usepackage{graphicx,float}
\usepackage{subfig}
\usepackage{placeins}
\usepackage{xcolor}
\usepackage{indentfirst}
\usepackage{enumitem}
\usepackage{caption}
\usepackage{float}
\usepackage{bm}
\usepackage{algorithm}
\usepackage{algorithmic}
\usepackage{tikz}
\usetikzlibrary{arrows.meta,patterns,calc,decorations.pathreplacing}

\usepackage[colorlinks,
            linkcolor=blue,
            anchorcolor=blue,
            citecolor=blue]{hyperref} 
\numberwithin{equation}{section}
\newcommand{\R}{{\mathbb R}}

\newcommand{\be}{\begin{eqnarray}}
\newcommand{\en}{\end{eqnarray}}
\newcommand{\ben}{\begin{eqnarray*}}
\newcommand{\enn}{\end{eqnarray*}}

\newcommand{\e}{\mathrm{e}}
\newcommand{\inc}{\mathrm{inc}}

\newtheorem{theorem}{Theorem}[section]
\newtheorem{lemma}[theorem]{Lemma}

\newtheorem{remark}[theorem]{Remark}

\definecolor{rot}{rgb}{1.000,0.000,0.000}
\definecolor{rot1}{rgb}{0.000,0.000,0.000}
\begin{document}
\renewcommand{\theequation}{\arabic{section}.\arabic{equation}}
\renewcommand{\figurename}{Figure}
\captionsetup[figure]{labelformat=simple}
\renewcommand{\tablename}{Table}
\captionsetup[table]{labelformat=simple}
\begin{titlepage}
\title{A new frequency-time hybrid RPML-BIE solver with few-mode corrections for the wave-equation problems in nonuniform closed waveguides}

\author{Yan Tan\thanks{School of Mathematical Sciences, University of Electronic Science and Technology of China, Chengdu, Sichuan 611731, China. Email: {\tt yantan@std.uestc.edu.cn
}}\;,
Liwei Xu\thanks{School of Mathematical Sciences, University of Electronic Science and Technology of China, Chengdu, Sichuan 611731, China. Email: {\tt xul@std.uestc.edu.cn
}}\;,
Tao Yin\thanks{State Key Laboratory of Mathematical Sciences and Institute of Computational Mathematics and Scientific/Engineering Computing, Academy of Mathematics and Systems Science, Chinese Academy of Sciences, Beijing 100190, China. Email:{\tt yintao@lsec.cc.ac.cn}}}

\date{}
\end{titlepage}
\maketitle

\begin{abstract} 
This paper is devoted to proposing a novel frequency-time hybrid (FTH) integral equation solver for time-dependent wave-equation problems in two-dimensional nonuniform closed waveguides with multiple channels. The FTH methodology, which combines the Fourier transform with frequency-domain boundary integral equation (BIE) formulations, has demonstrated its effectiveness in solving time-dependent wave scattering problems. However, when applied to waveguide problems, the modal expansions arising from multiple scattering introduce cut-off frequencies. Such frequencies pose significant challenges, requiring both a frequency-robust solver for frequency-domain waveguide problems and an efficient treatment of the solution's non-smooth behavior near the cut-off frequencies for the evaluation of the inverse Fourier transform. The first issue is addressed by developing a new frequency-domain BIE method based on the free-space fundamental solution and regionalized perfectly matched layer truncation, together with a few-mode Dirichlet-to-Robin modification. The method is rigorously proved to achieve uniform exponential convergence, even at frequencies arbitrarily close to the cut-off frequencies. To address the second issue---the recovery of the time-domain solution via the inverse Fourier transform---a sine-transform based stationary-point Filon--Clenshaw--Curtis quadrature is introduced to accurately handle the non-smoothness and square-root singularities of frequency-domain solutions. Additionally, theoretical result shows that the complex resonances closing to the real axis may significantly effect the convergence of the Fourier transform algorithm and then a novel strategy of complex resonance extraction based on AAA rational approximation and singularity subtraction is utilized to avoid this influence. Numerical experiments for both frequency-domain and time-domain problems illustrate the efficiency and high accuracy of the proposed method.

\end{abstract}

\section{Introduction}
\label{sec:1}

Developing accurate and efficient numerical solvers for the simulation of time-dependent wave propagation in complex media or structures---for example the nonuniform waveguides which find many important applications in acoustic devices, microwave engineering, integrated circuits, crystal fibers and so on---still remains challenging. This paper is devoted to propose a new high-accuracy solver for the wave-equation problems in nonuniform closed waveguides. One of the most popular ideas to treat the wave-equation problems uses volumetric discretization of the spatial domain such as the finite difference method~\cite{TH00}, finite element method (FEM)~\cite{LLC97}, in conjunction with time-stepping discretization such as the upwind/downwind scheme, Newmark method~\cite{GSS06} and Runge--Kutta method~\cite{SS19}. This idea has been extensively studied in open literature and applications, however, generally suffers from spatial and temporal numerical dispersion errors~\cite{BS97}. The recently developed space-time finite element method~\cite{BMPS21,LSZ22} provides an efficient, unconditionally stable scheme for solving the wave-equation problems, however, whose application to problems with complex structures still remains unresolved.

The boundary integral equation (BIE) method introduces an alternative way to transform the original problems into appropriate BIE formulations on the boundaries, for example the time-domain boundary integral equation (TDBIE) method~\cite{S16,HD03} and convolution quadrature (CQ) method~\cite{L88a,L88b,L94,B10}. The TDBIE method uses the retarded-potential fundamental solution or Green's function of the wave equation and hence, requires to treat singular integration over the space-time intersection of the light cone with the boundary of domain, making stable high-order discretization challenging. The CQ method, on the other hand, does not handle the wave-equation problem directly, but utilizes a combination of time-stepping discretization and Laplace transformation to reduce the wave equation to Helmholtz equation with complex wavenumbers, which can be solved by means of the well-constructed frequency-domain BIE methods~\cite{RS07} or other feasible solvers~\cite{LH21}. This method coupled with the FEM and boundary element discretization has also been extended to the wave-equation problem with an inhomogeneous medium~\cite{B22,BLS15}. However, it is known~\cite{BSS17} that the stability of the CQ method depends on the geometry as well as the choice of the BIEs and it inherits the dispersive character of the time-stepping discretization. Moreover, as far as we know, the efficiency and convergence of the CQ method to the complicated nonuniform waveguide problems are still unknown.

Besides the Laplace transformation, the Fourier transformation can alternatively be applied to conditionally reduce the time-dependent wave-equation problems to frequency-domain Helmholtz-equation problems with real frequencies and then the time-domain solutions can be evaluated by the inverse Fourier transform. For the exterior scattering problem, a frequency-time hybrid (FTH) BIE method~\cite{ABL20}, in conjunction with a novel windowing-and-recentering strategy, has been developed to provide essentially dispersion-free numerical solutions with spectral accuracy in time---supporting parallelization in both space and time and efficiently accommodating incident fields that persist indefinitely. However, strongly trapping configurations may produce complex resonances close to the real frequency axis, leading to nearly singular integrands and slow convergence of the numerical inverse Fourier transform. To address this difficulty, a frequency-domain singularity subtraction technique was developed in \cite{BS26} to regularize the integrand and evaluate the resonance contributions separately. The multi-patch multiple-scattering strategy proposed in \cite{BY24,PBYB25} provides an extension of the FTH method applicable to interior problems which proceeds by exploiting a multiple scattering technique to re-express the full time-evolution as problems of multiple scattering between multiple overlapping patches that covers the domain boundary. However, applying these ideas to nonuniform waveguide problems is highly nontrivial since both frequency-robust BIE solvers and efficient Fourier transform algorithm that account for both cut-off frequency singularities and complex resonances are necessary.

For frequency-domain waveguide problems, BIE formulations can be constructed using either the background waveguide Green's function or the free-space fundamental solution. The former incorporates the waveguide boundary conditions and the outgoing modal expansion, but the analytical form of the Green's function is unclear for general nonuniform waveguides and its kernel will blow up near the cut-off frequencies. The free-space fundamental solution does not suffer from such kind of singularities, but the resulting integral equations are posed on the whole unbounded surfaces of the waveguides and hence, require an appropriate truncation. The windowed Green function (WGF) method provides one truncation technique to achieve super-algebraically convergence and have been successfully applied to nonuniform open waveguides \cite{BGP17,GSB23}. Another possible approach uses perfectly matched layer (PML) truncation \cite{B94}, whose exponential convergence has been established for many classes of unbounded scattering problems \cite{LS01,CW03,CZ10} while considering the domain discretization schemes. Combining PML truncation with BIE formulations, leading to the so-called PML-BIE method, has been used to solve the problems of acoustic layered-medium \cite{LLQ18,LXYZ23}, electromagnetic multi-layered medium \cite{BLYZ24}, locally defected periodic surfaces \cite{YHLR22}, and water wave \cite{BFP24}. However, as indicated in Section~\ref{sec:2.3}, applying both of the WGF and PML truncation directly will suffer from rapidly deteriorated accuracy near the cut-off frequencies. 

To extend the PML-BIE method to the frequency-domain Helmholtz equation problems in nonuniform closed waveguides and noting that a single Cartesian stretching direction cannot align with all tilted ports, a regionalized perfectly matched layer (RPML) is designed by applying the stretching along the local longitudinal direction of each port of the nonuniform waveguide. Since the exponential convergence rate of standard RPML truncation will be destroyed near the cut-off frequencies and generally increasingly thick layers are needed to maintain a prescribed accuracy \cite{K19}, a few-mode Dirichlet-to-Robin (DtR) boundary condition, which is in fact inspired by \cite{ZW18}, is introduced to truncate the RPML at the outer boundaries instead of directly imposing a vanishing boundary condition. This ensures a uniform
exponential convergence for all frequencies, even at those arbitrarily close to the cut-off frequencies, and then allows us to establish a high-accuracy PML-BIE formulation for the frequency domain problems.

Further difficulties arise when the frequency-domain solutions are used to evaluate the time-domain solutions via the inverse Fourier transform along the real frequency axis. It can be observed from the explicit form of the background Green's function of the planar waveguide problem that the frequency-domain solutions suffer from a combined square-root singularity near the cut-off frequencies. This will cause slow convergence of the numerical evaluation of the inverse Fourier transform. To resolve this issue, we partition the integration interval at these cut-off frequencies and utilize a sine-transform to regularize the square-root singularities. The sine-transform introduces a nonlinear oscillatory phase with endpoint stationary points for which a stationary-point Filon--Clenshaw--Curtis (SP-FCC) quadrature \cite{DGS11,DGK13} is developed. However, the theoretical result shows that the convergence of this inverse Fourier transform algorithm depends on the size of a Bernstein ellipse in the complex plane and the complex resonances should be excluded. In addition, even when away from the cut-off frequencies, these complex resonance close to the real frequency axis can produce sharp, spike-like features in the Helmholtz frequency-domain solution, causing the integrand of the inverse Fourier integral to exhibit nearly singular behavior and resulting in extremely slow convergence of the numerical evaluation of the inverse Fourier transform. Resolving these features requires dense frequency sampling, which substantially increases the number of frequency-domain computations and the cost of evaluating the inverse Fourier transform \cite{BS26,LSW21}. Following~\cite{BS26}, we propose a complex resonance extraction strategy based on AAA rational approximation \cite{NST18} and singularity subtraction to separate the nearly singular contributions from the remainder which can be integrated by the SP-FCC quadrature, and the resonance contributions are evaluated separately and added back to obtain the time-domain solution with high accuracy.

This paper is organized as follows. Section~\ref{sec:2} introduces the considered time-domain wave-equation problem in nonuniform waveguides as well as the related frequency-domain model and the challenges in applying the FTH method. Section~\ref{sec:3} develops a frequency-robust RPML-BIE method based on a few-mode DtR correction technique and derives rigorous exponential convergence analysis. The Fourier-transform algorithm is presented in 
Section~\ref{sec:4} to resolve the issues arising from the cut-off frequency singularity and complex resonances. Section~\ref{sec:5} displays numerical experiments to demonstrate the efficiency and accuracy of the proposed method. Finally, Section~\ref{sec:6} gives a brief conclusion with further extensions.

\section{Preliminaries}
\label{sec:2}
\subsection{Time-domain Wave-equation problem}
\label{sec:2.1}

\begin{figure}[htbp]
\centering
\includegraphics[width=0.85\textwidth,trim=2pt 2pt 2pt 2pt,clip]{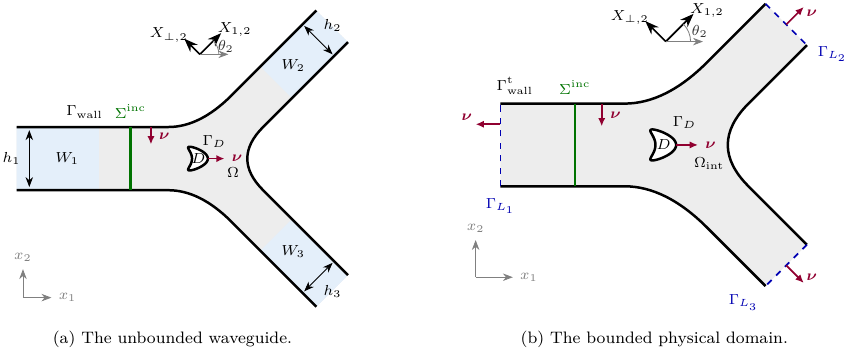}
\caption{Regions, boundaries and normals for (a) the unbounded waveguide and (b) the bounded physical domain.}
\label{fig:geometry}
\end{figure}

This paper considers time-dependent wave scattering in a two-dimensional closed, nonuniform multi-port waveguide; see Figure~\ref{fig:geometry}. The unbounded domain $\Omega\subset\mathbb R^2$ consists of a bounded region $\Omega_{\rm int}$, possibly containing impenetrable obstacles $D$, and $N_p$ non-overlapping straight semi-infinite ports $W_j$. The physical boundary is $\Gamma=\Gamma_D\cup\Gamma_{\rm wall}$, where $\Gamma_D=\partial D$. We take $\bm\nu$ to point into $\Omega$ on $\Gamma_D\cup\Gamma_{\rm wall}$ and write $\partial_{\bm\nu}=\bm\nu\cdot\nabla$.

For a port of width $h_j$ and inclination $\theta_j$, introduce local coordinates $\bm X^{(j)}=(X_{1,j},X_{\perp,j})^\top$ by
\be\label{eq:local_coords}
\bm X^{(j)}(\bm x)=\mathbb Q(\theta_j)\bm x-\bm d_j,\quad
\mathbb Q(\theta_j)=\begin{pmatrix}\cos\theta_j&\sin\theta_j\\-\sin\theta_j&\cos\theta_j\end{pmatrix}.
\en
The coordinate $X_{1,j}$ increases outward along the port and $0<X_{\perp,j}<h_j$. Choose a cross-section $\Gamma_{L_j}=\{\bm x:X_{1,j}=L_j,\ 0<X_{\perp,j}<h_j\}$ in each straight port, beyond the junction and the point source. Here $L_j$ need not be the beginning of the straight part: $\Omega_{\rm int}$ includes all physical port portions up to these sections, while $W_j$ extends through $\Gamma_{L_j}$ to infinity. On $\Gamma_{L_j}$ and the outer RPML boundaries, $\bm\nu$ points out of the bounded computational region.

To describe wave scattering in $\Omega$, we prescribe the temporal waveform independently of the spatial type of incidence. We consider the single Gaussian pulse
\be\label{eq:td_spectrum}
a(t)=\frac{\sigma}{2\sqrt{\pi}}\e^{-\sigma^2(t-t_0)^2/4}\e^{-\mathbf{i}\omega_0(t-t_0)},
\en
whose Fourier transform is given by $A(\omega)=\e^{-(\omega-\omega_0)^2/\sigma^2}\e^{\mathbf{i}\omega t_0}$. Here $\omega_0>0$ is the central frequency, $\sigma>0$ controls the bandwidth, and $t_0$ is chosen so that the pre-initial Gaussian tail is below the prescribed temporal tolerance. The Fourier pair is understood up to this negligible tail. We choose $0<W_-<W_+$ so that the response outside $[W_-,W_+]$ can be neglected. Signals requiring a zero-frequency contribution are outside the present scope. We consider the following two spatial types of incidence.

\begin{itemize}  
\item {Case I}. For a \emph{mode incidence}, the incident field
is prescribed by the $n_0$-th incoming guided mode in a selected
port $W_{j_0}$ such that $u^{\rm inc}(\bm x,t)$ is the inverse Fourier transform of
\be\label{eq:inc_mode}
U^{\rm inc}(\bm x,\omega)=U_{\mathrm{mode}}^{\rm inc}(\bm x,\omega):=\psi_{n_0}^{(j_0)}(X_{\perp,j_0})\,\e^{-\mathbf{i}\beta_{n_0}^{(j_0)}(\omega)X_{1,j_0}}\,A(\omega),
\en
where the formulas of the eigenfunctions $\psi_n^{(j)}$ and $\beta_n^{(j)}(\omega)$ corresponding to $W_j$ will be  introduced in Section~\ref{sec:2.2}. This incoming mode is a solution of the wave equation in $\Omega^\inc$ such that $W_{j_0}\subset \Omega^\inc$ for some $j_0=1,\cdots,N_p$ and satisfies $\partial_{\bm\nu} U^{\rm inc}=0$ on $\partial \Omega^\inc\cap\Gamma_{\rm wall}$. Denote $\Sigma^{\rm inc}=\partial \Omega^\inc\backslash(\partial \Omega^\inc\cap\Gamma_{\rm wall})$, see Figure~\ref{fig:geometry}. Introduce a characterized function
\ben
\chi^\inc(\bm x)=
\begin{cases}
1,& \bm x\in\Omega^\inc,\\
0,& \bm x\notin\Omega^\inc,
\end{cases}
\enn
so that the total field can be defined by
\ben
u^{\rm tot}(\bm x,t)=\chi^\inc(\bm x)u^{\rm inc}(\bm x,t)+u^{\rm sca}(\bm x,t),
\enn
where $u^{\rm sca}$ denotes the scattered field.

\item {Case II}. For a \emph{point source incidence}, the incident field is the free-space Green function such that $u^{\rm inc}(\bm x,t)$ is the inverse Fourier transform of
\be\label{eq:inc_source}
U^{\rm inc}(\bm x,\omega)=U_{\mathrm{point}}^{\rm inc}(\bm x,\omega):=A_0\,\Phi_\omega(\bm x,\bm z)\,A(\omega),
\quad
\Phi_\omega(\bm x,\bm y)=\frac{\mathbf{i}}{4}H_0^{(1)}\!\left(\frac{\omega}{c}|\bm x-\bm y|\right),
\en
where $\bm z$ denotes the source location, $A_0>0$ is an amplitude, and $H_0^{(1)}$ denotes the zeroth-order Hankel function of the first kind. In this case, the total field is given by
\ben
u^{\rm tot}(\bm x,t)=u^{\rm inc}(\bm x,t)+u^{\rm sca}(\bm x,t).
\enn
\end{itemize}

Then the studied time-dependent wave scattering problem is to find the sought field $u^\star$ such that
\be\label{initial-boundary value problem for wave equation2}
\begin{cases}
\partial_t^2 u^\star(\bm x,t)-c^2\Delta u^\star(\bm x,t)=0,
& (\bm x,t)\in\Omega\times\mathbb{R}^+,\\
\partial_{\bm\nu} u^\star(\bm x,t)+\mathbf{i}\eta_D u^\star(\bm x,t)=f_D^\star(\bm x,t),
& (\bm x,t)\in\Gamma_D\times\mathbb{R}^+,\\
\partial_{\bm\nu} u^\star(\bm x,t)=f^\star(\bm x,t),
& (\bm x,t)\in\Gamma_{\mathrm{wall}}\times\mathbb{R}^+,\\
u^\star(\bm x,0)=\partial_tu^\star(\bm x,0)=0, & \bm x\in\Omega,
\end{cases}
\en
with $c>0$ being the wave speed which is a constant and $\eta_D>0$, where 
\be\label{eq:ustar_def}
u^\star=
\begin{cases}
u^{\rm tot},& \text{mode incidence},\\
u^{\rm sca},& \text{point-source incidence}.
\end{cases}
\en
For mode incidence, $f_D^\star=f^\star=0$; for point-source incidence, $f_D^\star=-\partial_{\bm\nu} u^{\rm inc}-\mathbf{i}\eta_Du^{\rm inc}$ and $f^\star=-\partial_{\bm\nu} u^{\rm inc}$. The bounded obstacle $D$ and the impedance boundary condition on its boundary $\Gamma_D$ are included for the convenience of analysis of uniqueness, however, can be neglected in numerical calculation.
The discussion for other types of boundary conditions on $\Gamma_{\mathrm{wall}}$ follows analogously.

The derivation below only uses point-source incidence for simplicity; the modifications for mode incidence are summarized in Remark~\ref{remark:mode-incidence}.

\subsection{Frequency-domain model}
\label{sec:2.2}

For a given time-domain function $g\in L^2(\mathbb{R})$ with $g(t)=0, t\le 0$, we use the Fourier transform
\ben
G(\omega)=\mathbb{F}[g](\omega):=\int_{0}^\infty g(t)e^{\mathbf{i}\omega t}\mathrm{d}t,\quad\omega\in\R,
\enn
and, for $G\in L^1(\mathbb R)$, the inverse transform
\ben
g(t)=\mathbb{F}^{-1}[G](t):=\frac{1}{2\pi}\int_{-\infty}^\infty G(\omega)e^{-\mathbf{i}\omega t}\mathrm{d}\omega,\quad t\in\R.
\enn

Let $U^{\rm tot}$ and $U^{\rm sca}$ denote the Fourier transforms of $u^{\rm tot}$ and $u^{\rm sca}$, respectively. Then it gives that the scattered field $U^{\rm sca}$ satisfies the following problem of Helmholtz equation with wavenumber
$\kappa=\kappa(\omega)=\omega/c$:
\be\label{eq:fd_bvp}
\Delta U^{\rm sca}+\kappa^2U^{\rm sca}=0\ \text{ in }\Omega,
\quad
\partial_{\bm\nu} U^{\rm sca}+\mathbf{i}\eta_D U^{\rm sca}=F_D^\star\ \text{ on }\Gamma_D,
\quad
\partial_{\bm\nu} U^{\rm sca}=F^\star\ \text{ on }\Gamma_{\mathrm{wall}},
\en
where $F_D^\star=-\partial_{\bm\nu}U^{\rm inc}-\mathbf{i}\eta_D U^{\rm inc}$ and $F^\star=-\partial_{\bm\nu}U^{\rm inc}$ are the transforms of the point-source boundary data. To ensure the well-posedness of the frequency-domain model of waveguide scattering problem, additional mode expansions in each port is required and it leads to study the 
eigensystem $\{\psi_n^{(j)},\alpha_n^{(j)}\}_{n=0}^{\infty}$ satisfying
\ben
\begin{cases}
\frac{\mathrm{d}^2\psi_n^{(j)}}{\mathrm{d}X_{\perp,j}^2}(X_{\perp,j})
+(\alpha_n^{(j)})^2\psi_n^{(j)}(X_{\perp,j})=0,
& X_{\perp,j}\in(0,h_j),\\[2mm]
\frac{\mathrm{d}\psi_n^{(j)}}{\mathrm{d}X_{\perp,j}}(X_{\perp,j})=0,
& X_{\perp,j}=0,h_j.
\end{cases}
\enn 
Clearly,
\ben
\alpha_n^{(j)}=\frac{n\pi}{h_j},
\quad
\psi_n^{(j)}(X_{\perp,j})=
\sqrt{\frac{2-\delta_{n0}}{h_j}}
\cos\left(\frac{n\pi}{h_j}X_{\perp,j}\right),
\quad n\in\mathbb{N},
\enn
where $\delta_{n0}$ is the Kronecker delta. Then for the Helmholtz equation problem in each port $W_j$, separation of variables gives the outgoing and incoming guided modes $\e^{\pm\mathbf{i}\beta_n^{(j)}(\omega)X_{1,j}}\psi_n^{(j)}(X_{\perp,j}),n\in\mathbb N$, where
\ben
\beta_n^{(j)}(\omega)=\begin{cases}
  \sqrt{\kappa^2-(\alpha_n^{(j)})^2},&\mathrm{if}~(\alpha_n^{(j)})^2\leq \kappa^2,\\
  \mathbf{i}\sqrt{(\alpha_n^{(j)})^2-\kappa^2},&\mathrm{if}~(\alpha_n^{(j)})^2> \kappa^2.
\end{cases}
\enn
Then we impose the total field due to the point-source incidence to admit the strictly outgoing
mode expansion
\be\label{mode expansion}
U^{\rm tot}(\bm x,\omega)
=
\sum_{n=0}^{\infty}
a_{n}^{(j)}
\e^{\mathbf{i}\beta_n^{(j)}(\omega)X_{1,j}}
\psi_n^{(j)}(X_{\perp,j}),
\quad
X_{1,j}\ge L_j,\quad j=1,\ldots,N_p.
\en
Here $a_{n}^{(j)}$ denote the outgoing modal coefficients.

To study the well-posedness of the frequency-domain model (\ref{eq:fd_bvp})-(\ref{mode expansion}), we now introduce an artificial boundary condition on $\Gamma_{L_j}$ to truncate the infinite region $\Omega$ based on the so-called DtR map. Due to the continuity of $U^{\rm tot}$ as well as its normal derivative across $\Gamma_{L_j}$ and the mode expansion (\ref{mode expansion}), one can define a Robin-type boundary condition on $\Gamma_{L_j}$ for the total field $U^{\rm tot}$ as
\be\label{eq:exact-DtR-BC}
\partial_{\bm\nu} U^{\rm tot}+\mathbf{i}\eta U^{\rm tot}-\mathcal T_{j}U^{\rm tot}=0
\quad\text{on }\Gamma_{L_j},\quad j=1,\ldots,N_p.
\en
Here, $\eta\le 0$ is the Robin parameter shared by all ports. For any function $f(\bm x)$ on $\Gamma_{L_j}$ which admits the Fourier expansion $f(\bm x)=\sum_{n=0}^{\infty}f_n^{(j)}\psi_n^{(j)}(X_{\perp,j})$ with $f_n^{(j)}=\int_0^{h_j}f(\bm x)\overline{\psi_n^{(j)}(X_{\perp,j})}\,\mathrm dX_{\perp,j}$, $\mathcal T_{j}$ denotes a DtR map given by
\be\label{eq:exact-DtR-op}
\mathcal T_j f(\bm x)
=\sum_{n=0}^{\infty}\mathbf{i}\left(\beta_n^{(j)}(\omega)+\eta\right)f_n^{(j)}\psi_n^{(j)}(X_{\perp,j}),\quad \bm x\in\Gamma_{L_j}, j=1,\ldots,N_p,
\en
When $\eta=0$, the DtR map reduces to the more popular DtN map. Now we obtain an equivalent boundary value problem in $\Omega_{\mathrm{int}}$:
\be\label{eq:fd_model}
\begin{cases}
\Delta U^{\rm sca}+\kappa^2U^{\rm sca}=0,
& \mbox{in}\quad\Omega_{\mathrm{int}},\\
\partial_{\bm\nu} U^{\rm sca}+\mathbf{i}\eta_D U^{\rm sca}
=F_D^\star,
& \mbox{on}\quad\Gamma_D,\\
\partial_{\bm\nu} U^{\rm sca}=F^\star,
& \mbox{on}\quad\Gamma_{\mathrm{wall}}^{\mathrm{t}},\\
\partial_{\bm\nu} U^{\rm sca}+\mathbf{i}\eta U^{\rm sca}-\mathcal T_{j}U^{\rm sca}=G_j,
& \mbox{on}\quad\Gamma_{L_j},\quad j=1,\ldots,N_p.
\end{cases}
\en
Here, $\Gamma_{\mathrm{wall}}^{\mathrm{t}}=\Gamma_{\mathrm{wall}}\cap\partial\Omega_{\mathrm{int}}$ and $G_j=-\partial_{\bm\nu} U^{\rm inc}-\mathbf{i}\eta U^{\rm inc}+\mathcal T_{j}U^{\rm inc}$ on $\Gamma_{L_j}, j=1,\ldots,N_p$.

The variational formulation of the problem (\ref{eq:fd_model}) reads: given $F_D^\star\in H^{-1/2}(\Gamma_D)$, $F^\star\in H^{-1/2}(\Gamma_{\mathrm{wall}}^{\mathrm{t}})$ and $G_j\in H^{-1/2}(\Gamma_{L_j}), j=1,\ldots,N_p$, find $U^{\rm sca}\in H^1(\Omega_{\mathrm{int}})$ such that
\be\label{eq:weak_original}
b(U^{\rm sca},\psi)=\ell(\psi)
\quad\forall\,\psi\in H^1(\Omega_{\mathrm{int}}),
\en
where 
\be\label{eq:form_b}
\begin{aligned}
b(\varphi,\psi)
=&\int_{\Omega_{\rm int}}\left(\nabla\varphi\cdot\overline{\nabla\psi}-\kappa^2\varphi\bar\psi\right)\mathrm d\bm x
-\mathbf{i}\eta_D\int_{\Gamma_D}\varphi\bar\psi\,\mathrm ds+\sum_{j=1}^{N_p}\left[\mathbf{i}\eta\int_{\Gamma_{L_j}}\varphi\bar\psi\,\mathrm ds-\langle\mathcal T_j\varphi,\psi\rangle_{\Gamma_{L_j}}\right],\\
\ell(\psi)
=&-\langle F_D^\star,\psi\rangle_{\Gamma_D}
-\langle F^\star,\psi\rangle_{\Gamma_{\rm wall}^{\rm t}}
+\sum_{j=1}^{N_p}\langle G_j,\psi\rangle_{\Gamma_{L_j}},
\end{aligned}
\en
with $\langle\cdot,\cdot\rangle_\Gamma$ denoting the duality pairing between $H^{-1/2}(\Gamma)$ and $H^{1/2}(\Gamma)$ on some boundary $\Gamma$.

\begin{theorem}\label{thm:wellposed_original}
For any $\omega>0$ such that all $\beta_n^{(j)}(\omega)$ do not vanish, the variational problem \eqref{eq:weak_original} has a unique solution $U^{\rm sca}\in H^1(\Omega_{\mathrm{int}})$. Furthermore, there exists a constant $\gamma=\gamma(\omega)>0$ such that
\be\label{eq:infsup}
\sup_{0\neq\psi\in H^1(\Omega_{\mathrm{int}})} \frac{\left|b(\varphi,\psi)\right|}{\|\psi\|_{H^1(\Omega_{\mathrm{int}})}}\ge\gamma\,\|\varphi\|_{H^1(\Omega_{\mathrm{int}})}\quad\forall\,\varphi\in H^1(\Omega_{\mathrm{int}}).
\en
\end{theorem}

\begin{proof}
For $\varphi|_{\Gamma_{L_j}}=\sum_n\varphi_n^{(j)}\psi_n^{(j)}$, the DtR definition \eqref{eq:exact-DtR-op} and orthonormality give
\be\label{eq:exact-DtR-energy}
\begin{aligned}
\mathbf{i}\eta\|\varphi\|_{L^2(\Gamma_{L_j})}^2-\left\langle\mathcal T_j\varphi,\varphi\right\rangle_{\Gamma_{L_j}}
&=\sum_{n=0}^{\infty}\left[\mathbf{i}\eta-\mathbf{i}\left(\beta_n^{(j)}+\eta\right)\right]|\varphi_n^{(j)}|^2\\
&=\sum_{\alpha_n^{(j)}>\kappa}|\beta_n^{(j)}|\,|\varphi_n^{(j)}|^2-\mathbf{i}\sum_{\alpha_n^{(j)}<\kappa}\beta_n^{(j)}|\varphi_n^{(j)}|^2.
\end{aligned}
\en
Let $b(u,\psi)=0$ for all $\psi$ and take $\psi=u$. By \eqref{eq:exact-DtR-energy},
\ben
0=\operatorname{Im}b(u,u)
=-\eta_D\int_{\Gamma_D}|u|^2\,\mathrm ds
-\sum_{j=1}^{N_p}\sum_{n:\,\alpha_n^{(j)}<\kappa}\beta_n^{(j)}\left|u_n^{(j)}\right|^2.
\enn
Thus $u=\partial_{\bm\nu} u=0$ on $\Gamma_D$ which further gives $u\equiv0$ in $\Omega_{\rm int}$ due to the Holmgren's Theorem. Boundedness of the DtR and trace maps implies boundedness of $b$, while \eqref{eq:exact-DtR-energy} gives the G\aa rding inequality
\ben
\operatorname{Re}b(\varphi,\varphi)+(1+\kappa^2)\|\varphi\|_{L^2(\Omega_{\rm int})}^2
\ge\|\varphi\|_{H^1(\Omega_{\rm int})}^2.
\enn
Then the well-posedness of the variational problem follows from the standard Fredholm argument and then the inf--sup condition \eqref{eq:infsup} results.
\end{proof}

\begin{remark}
It should be pointed out that in real application and practical calculation, the DtR map should be truncated as a finite $N$-terms sum and the well-posedness of the reduced variational formulation does not follows trivially. When $\eta=0$, i.e., when the DtR map becomes the classical DtN map, analogous analysis as \cite{ZW18} shows that the reduced variational formulation admits a unique solution for sufficiently large $N$. 
\end{remark}

\subsection{FTH strategy and challenges}
\label{sec:2.3}

Let $f_B$ denote collectively the prescribed boundary data, and let $F_B=\mathbb F[f_B]$. Utilizing appropriate BIE method to solve the problem (\ref{eq:fd_model}), a direct application of FTH method would perform
\be\label{eq:direct_fth_process}
f_B(\bm x,t)\xrightarrow{\mathbb F}F_B(\bm x,\omega)
\xrightarrow{\mathrm{appropriate~BIE~solver}}U^{\rm sca}(\bm x,\omega)
\xrightarrow{\mathbb F^{-1}}u^{\rm sca}(\bm x,t).
\en
For the point-source incidence \eqref{eq:inc_source}, let $U^{\rm unit}(\bm x,\omega)$ denote the scattered-field solution of \eqref{eq:fd_model} with $A(\omega)$ replaced by one. Linearity gives $U^{\rm sca}(\bm x,\omega)=A(\omega) U^{\rm unit}(\bm x,\omega)$. Neglecting the small temporal and frequency tails, the complex scattered field is
\be\label{eq:td_window_sum}
u^{\rm sca}(\bm x,t)\approx\frac{1}{2\pi}\int_{\omega_-}^{W}
A(\omega)U^{\rm unit}(\bm x,\omega)\e^{-\mathbf{i}\omega t}\,\mathrm d\omega.
\en
Accurate evaluation of \eqref{eq:td_window_sum} is highly non-trivial due to the following four issues.

\begin{figure}[htbp]
\centering
\begin{tabular}{cc}
\subfloat[]{\includegraphics[width=0.24\textwidth]{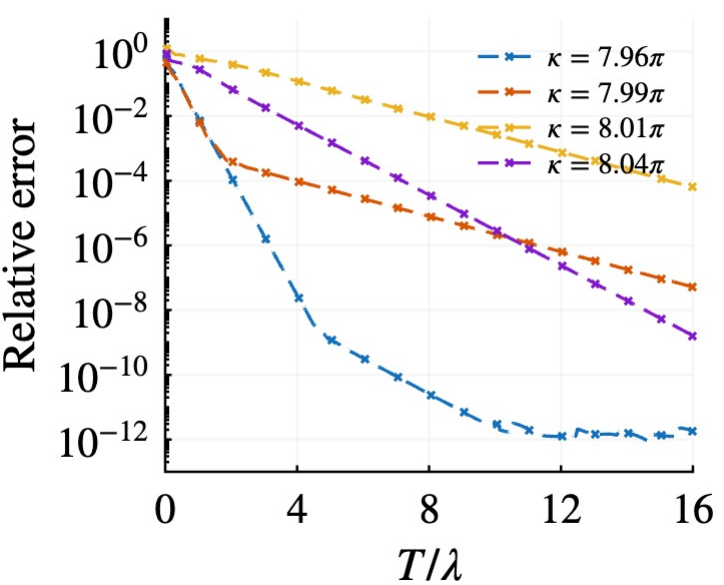}}
\subfloat[]{\includegraphics[width=0.24\textwidth]{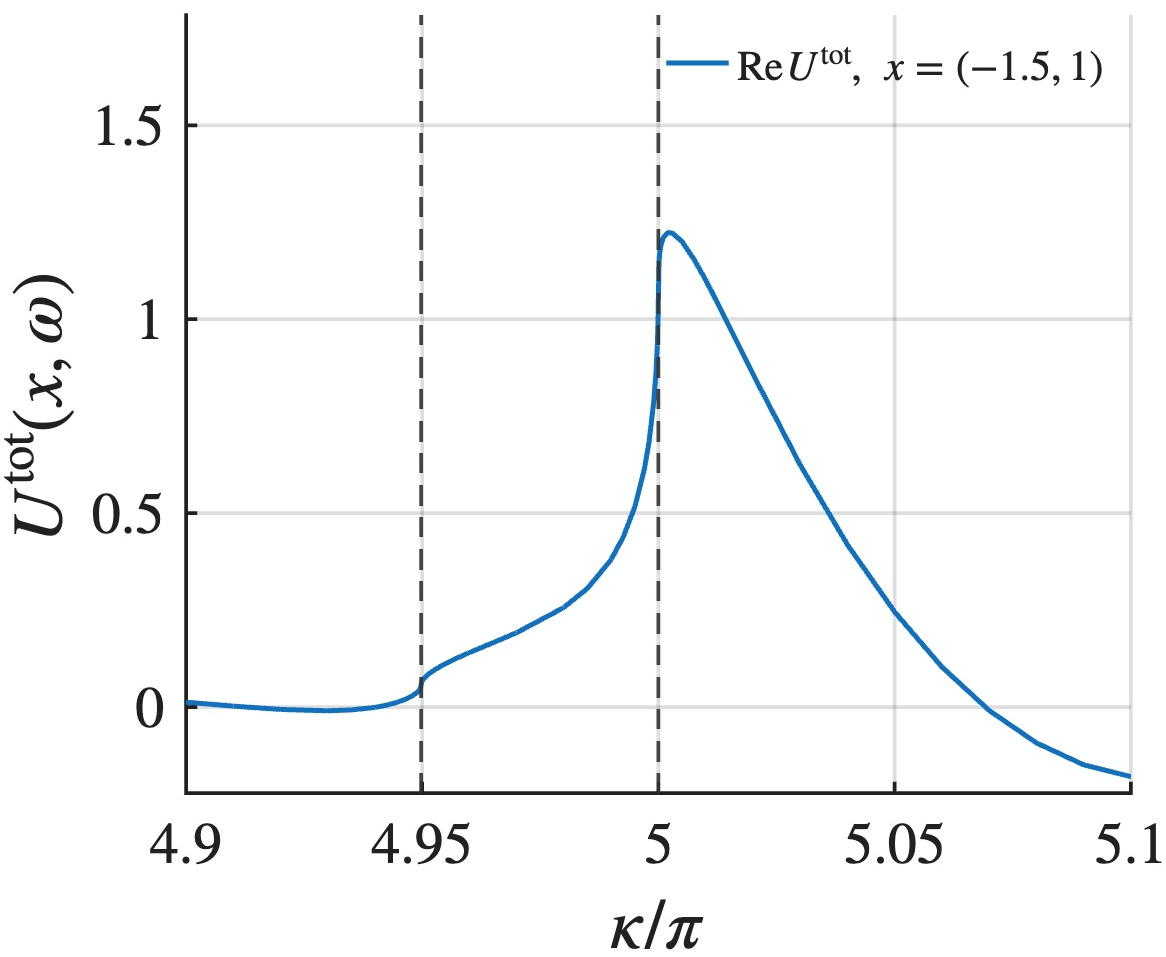}}
\subfloat[]{\includegraphics[width=0.24\textwidth]{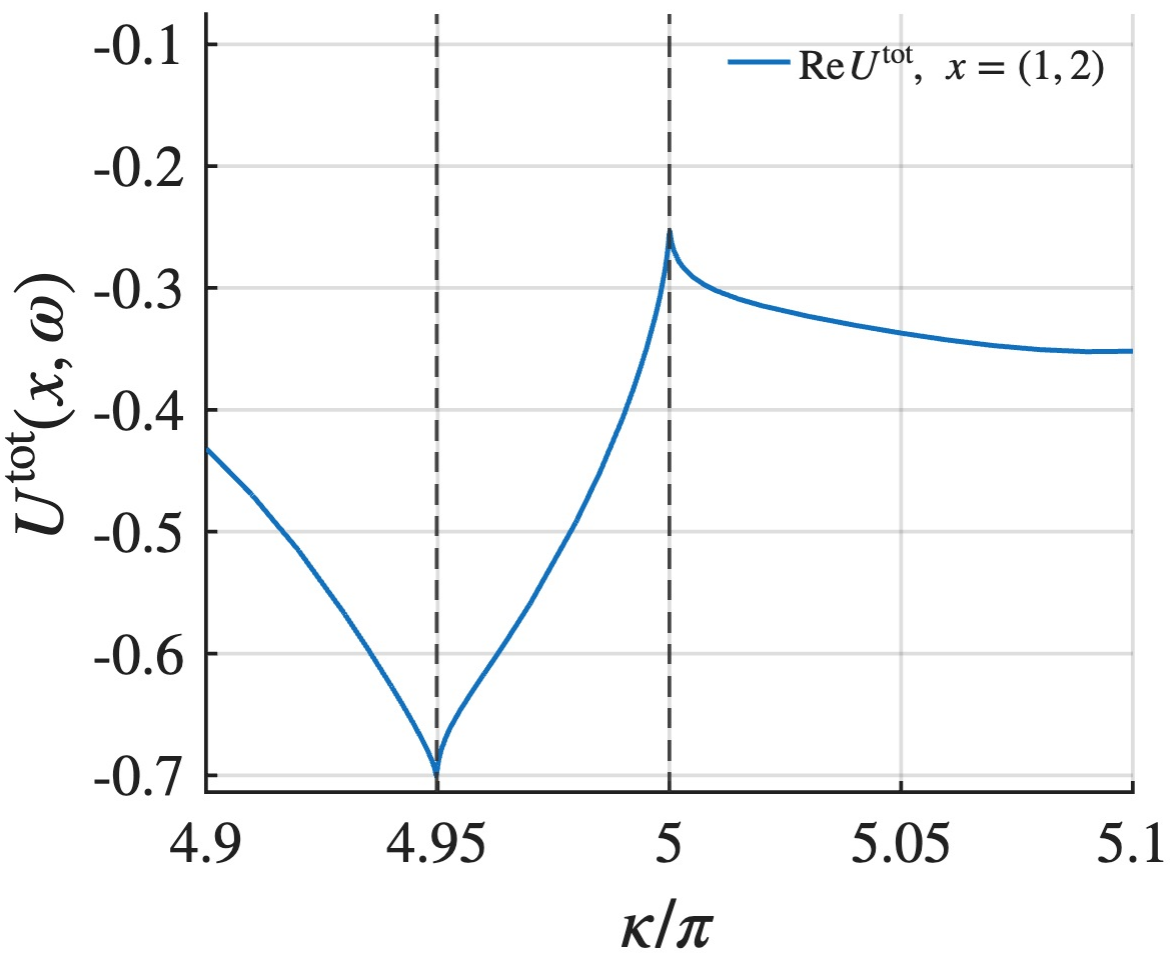}}
\subfloat[\label{fig:Ex3-3-AAA-poles}]{\includegraphics[width=0.24\textwidth]{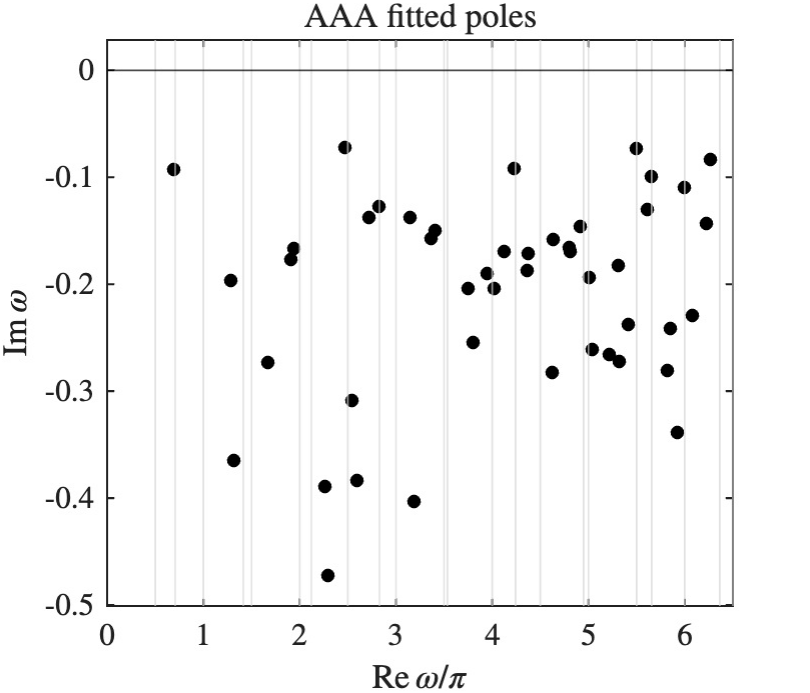}}
\end{tabular}
\caption{Challenges in developing FTH solver for waveguide problems: (a) destroyed convergence of the direct RPML truncation near the cut-off frequency $8\pi$ for the problem of a straight waveguide; (b)--(c) non-smoothness of the total field near cut-off frequencies at two points in the Y-junction structure shown in Figure~\ref{fig:numerical-geometries}(d) (dashed lines mark the cut-off frequencies); (d) locations of the complex resonances estimated by AAA for the Y-junction structure.}
\label{fig:cutoff}
\end{figure}

\begin{itemize}
\item[(I).] The BIE method using the background waveguide Green's function, whose mode expansion contains factors $1/\beta_n^{(j)}$, becomes singular at cut-off frequencies. Alternative approach based on the free-space fundamental solution can avoid this singularity at cut-off frequencies, but the resulting integral equations will be posed on the whole unbounded walls and hence require appropriate truncation strategies. The PML method provides 
one possible way to truncate the infinite domain, however, should be designed according to the direction along which the waveguide structure extends. A single Cartesian stretching direction cannot align with the local axes $X_{1,j},j=1,\cdots,N_p$ of all tilted ports.

\item[(II).] It is known that the PML method generally enjoys exponential convergence as the PML thickness increases. But this exponential convergence will be destroyed near the cut-off frequencies~\cite{K19}, i.e., $|\beta_n^{(j)}|\to0$. Maintaining a prescribed accuracy by the standard PML truncation requires an increasingly larger thickness near the cut-off frequencies; see Figure~\ref{fig:cutoff}(a) and Remark~\ref{rem:standard-RPML-cutoff}. Therefore, the BIE method based PML truncation should be corrected to ensure uniform accuracy within a wide range of frequencies. This issue also appears in the periodic structure scattering problems~\cite{SFFP23,BF20} near Rayleigh--Wood anomalies~\cite{TPY26}.

\item[(III).] The inverse Fourier transform requires integrating frequency-domain solutions; however, they may contain square-root singularities at cut-off frequencies{; see Figure~\ref{fig:cutoff}(b)--(c)}. These singularities must be treated to achieve high-order accuracy. The sine transform discussed in Section~\ref{sec:sine-SPFCC} can regularize the square-root singularities, but will transforms the linear Fourier moment into a nonlinear moment with stationary points.

\item[(IV).] The inverse Fourier transform is evaluated along the real frequency axis. The frequency-domain solution, however, may have complex resonances in the lower half-plane. This phenomenon is not specific to waveguides and also occurs for strongly trapping obstacles\cite{BS26}. In waveguide scattering, multiple scattering might keep waves localized near a junction for a long time. If energy escapes slowly through the ports, the associated complex resonances can lie extremely close to the real axis. Such resonances resonances manifest the frequency-domain solutions as sharp, spike-like features along the real frequency axis~\cite{BS26,LSW21}. As a result, achieving a prescribed accuracy then requires a dense sampling of frequency points in the affected sub-intervals---prohibitively increasing the computing cost of inverse Fourier transform. Figure~\ref{fig:cutoff}(d) shows the locations of complex resonances estimated by AAA method for the proposed RPML-BIE solver for the waveguide problem with a
Y-junction structure shown in Figure~\ref{fig:numerical-geometries}(d).
\end{itemize}

\section{Frequency-robust RPML-BIE method}
\label{sec:3}

We now first study the family of frequency-domain problems (\ref{eq:fd_bvp}) in nonuniform closed waveguides and propose a novel frequency-robust RPML-BIE formulation to resolve the challenges (I) and (II) raised in Section~\ref{sec:2.3}.

\subsection{RPML truncation and limitation}
\label{sec:3.1}

\begin{figure}[htbp]
\centering
\includegraphics[width=\textwidth,trim=2pt 2pt 2pt 2pt,clip]{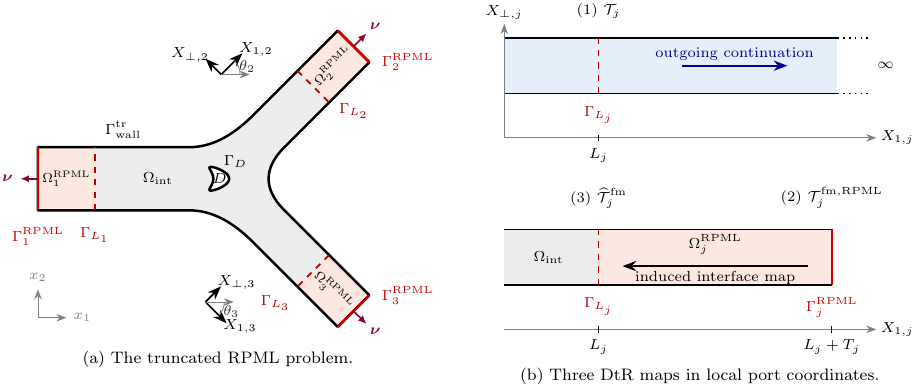}
\caption{RPML truncation and a port in rotated local coordinates, showing
(1) the exact DtR map, (2) the outer few-mode DtR truncation, and
(3) the induced interface DtR map.}
\label{fig:rpml-dtr-maps}
\end{figure}

To ensure the convergence of the RPML truncation in all channels with different inclinations, it is necessary to construct the complex stretching directly along the local longitudinal coordinate
$X_{1,j}$ of each port $W_j$. For the $j$-th port, we define the stretched local coordinates
$\widetilde{\bm X}^{(j)}
=(\widetilde X_{1,j},\widetilde X_{\perp,j})^\top$
as
\ben
\widetilde{X}_{1,j}(X_{1,j}) = X_{1,j}+\textbf{i}\int_{L_j}^{X_{1,j}}\sigma_j(t)\mathrm{d}t, \quad \text{and} \quad \widetilde{X}_{\perp,j} = X_{\perp,j},
\enn
where $L_j$ is the selected physical--RPML interface position (see for example Figure~\ref{fig:rpml-dtr-maps}(a)) and the stretching profile $\sigma_j(t)$ is
extended by zero throughout the physical region and defined by
\be\label{sigma}
\sigma_j(t)=\begin{cases}
   0, & t < L_j,\\
   \frac{2S_j\chi_1^{P_j}}{\chi_1^{P_j}+\chi_2^{P_j}}, & L_j\leq t\leq L_j+T_j,\\
   2S_j, & t>L_j+T_j.
\end{cases}
\en
Here $T_j>0$ is the RPML thickness, $P_j\geq 2$ is a positive integer, $S_j>0$ is the RPML scaling parameter, and 
\ben
\chi_1(s)=\left(\frac{1}{2}-\frac{1}{P_j}\right)s^3+\frac{s}{P_j}+\frac{1}{2},\quad s=\frac{2t-(2L_j+T_j)}{T_j},\quad\chi_2=1-\chi_1.
\enn

The locally stretched coordinates are mapped back to the global frame to define the piecewise global complex coordinates $\tilde{\bm x} = (\tilde{x}_1, \tilde{x}_2)^\top$:
\be\label{eq:global_stretch}
\tilde{\bm x}(\bm x) = 
\begin{cases}
\bm x, & \bm x \in \Omega_{\mathrm{int}},\\
\mathbb Q(-\theta_j) \left[\big(\widetilde{X}_{1,j}(\bm x), X_{\perp,j}(\bm x)\big)^\top + \bm d_j\right], & \bm x \in W_j,\ X_{1,j}>L_j,\quad j=1,\ldots,N_p.
\end{cases}
\en
Let $\alpha_j(X_{1,j}) = 1 + \textbf{i}\sigma_j(X_{1,j})$. By the chain rule, $\mathbb{B}_j(\bm x) = \frac{\partial \tilde{\bm x}}{\partial \bm x}$ is simply the identity matrix $\mathbb{I}$ in $\Omega_{\mathrm{int}}$, whereas for $X_{1,j}>L_j$, the fully coupled symmetric complex Jacobian matrix takes the form:
\ben
\mathbb{B}_j(\bm x) = \mathbb Q(-\theta_j) \begin{pmatrix} \alpha_j(X_{1,j}) & 0 \\ 0 & 1 \end{pmatrix} \mathbb Q(\theta_j)=\mathbb{B}_j^\top(\bm x).
\enn
Define the stretched incident, scattered and total fields by
\ben
\widetilde U^{\rm inc}(\bm x,\omega)=U^{\rm inc}(\widetilde{\bm x},\omega),\quad
\widetilde U^{\rm sca}(\bm x,\omega)=U^{\rm sca}(\widetilde{\bm x},\omega),\quad
\widetilde U^{\rm tot}=\widetilde U^{\rm inc}+\widetilde U^{\rm sca}.
\enn
The chain rule gives
\ben
\mathscr L_j\widetilde U^{\rm sca}
:=\nabla_{\bm x}\cdot(\mathbb A_j^{-1}\nabla_{\bm x}\widetilde U^{\rm sca})
+\kappa^2\det(\mathbb B_j)\widetilde U^{\rm sca}=0,\quad
\mathbb A_j^{-1}=\det(\mathbb B_j)\mathbb B_j^{-2},
\enn
with $\mathbb A_j^{-1}=\mathbb I$ in $\Omega_{\rm int}$. Writing $\widetilde\partial_{\bm\nu}v=\bm\nu\cdot(\mathbb A_j^{-1}\nabla v)$, we introduce the complex-stretching boundary conditions
\ben
\widetilde\partial_{\bm\nu}\widetilde U^{\rm sca}+\mathbf{i}\eta_D\widetilde U^{\rm sca}=\widetilde F_D^\star
\quad\text{on }\Gamma_D,\quad
\widetilde\partial_{\bm\nu}\widetilde U^{\rm sca}=\widetilde F^\star
\quad\text{on }\Gamma_{\rm wall},
\enn
where $\widetilde F_D^\star=-\widetilde\partial_{\bm\nu}\widetilde U^{\rm inc}-\mathbf{i}\eta_D\widetilde U^{\rm inc}$ and $\widetilde F^\star=-\widetilde\partial_{\bm\nu}\widetilde U^{\rm inc}$.
In addition, the outgoing total field $\widetilde U^{\rm tot}$ satisfies the stretched mode expansion
\ben
\widetilde U^{\rm tot}(\bm x,\omega)
=\sum_{n=0}^{\infty}a_n^{(j)}\e^{\mathbf{i}\beta_n^{(j)}(\omega)\widetilde X_{1,j}}\psi_n^{(j)}(X_{\perp,j}),
\quad X_{1,j}>L_j.
\enn

Given $T_j>0$, we truncate the RPML medium at $\Gamma_j^{\mathrm{RPML}}$ defined by 
\ben
\Gamma_j^{\mathrm{RPML}}=\{\bm x\in W_j:X_{1,j}=L_j+T_j,0<X_{\perp,j}<h_j\}.
\enn
The truncated computational domain $\Omega^{\mathrm{tr}}$ is then the union of $\Omega_{\rm int}$ and the closed layers
\be\label{eq:RPML_region}
\Omega_j^{\rm RPML}=\{\bm x\in W_j:L_j<X_{1,j}<L_j+T_j,\ 0<X_{\perp,j}<h_j\}.
\en
Denote $\Gamma_{\rm wall}^{\mathrm{tr}}=\Gamma_{\rm wall}\cap\partial\Omega^{\mathrm{tr}}$. The truncated geometry are illustrated in Figure~\ref{fig:rpml-dtr-maps}(a).

\begin{remark}[Direct RPML truncation near cut-off frequencies]\label{rem:standard-RPML-cutoff}
Analogous to~\cite{K19}, it can be shown that the convergence of the RPML truncation using a homogeneous Dirichlet boundary $\widetilde U^{\rm tot}=0$ on $\Gamma_j^{\rm RPML}, j=1,\cdots,N_p$ should depend on the smallest nonzero longitudinal wavenumber $\beta_n^{(j)}$. In particular, for the RPML stretching described in this section, a homogeneous Dirichlet boundary induces the truncation error
\ben
\mathcal{O}\left\{ \frac{2\beta_{\min}^{(j)}(\omega)}{\e^{2\beta_{\min}^{(j)}(\omega)\min(1,S_j)T_j}-1}\right\}\quad\mbox{with}\quad \beta_{\min}^{(j)}(\omega):=\min_{n\ge0}|\beta_n^{(j)}(\omega)|.
\enn
This gives exponential convergence at any fixed frequency with $\beta_{\min}^{(j)}>0$, but the rate degenerates as $\beta_{\min}^{(j)}\to0$ for fixed RPML parameters. In deed, when $\beta_{\min}^{(j)}\to0$, the convergence of the RPML truncation behaviors only algebraically as $\mathcal O((\min(1,S_j)T_j)^{-1})$, see also \cite[Remark~4.6]{K19}, which also happens when applying direct PML truncation to periodic structure scattering problems~\cite{TPY26}.
\end{remark}

\subsection{Few-mode DtR correction}
\label{sec:fm-DtR}

To remove this deterioration near the cut-off frequencies, instead of the direct truncation in terms of homogeneous boundary condition, we impose a DtR condition on a few selected modes and a homogeneous Robin condition on the others on $\Gamma_j^{\rm RPML}, j=1,\cdots,N_p$. 

Let $r$ be some prescribed positive integer. In each port $j,j=1,\cdots,N_p$, define the index minimizing the longitudinal-wavenumber magnitude by
\ben
n_{\min}^{(j)}(\omega)\in\operatorname*{arg\,min}_{n\ge0}|\beta_n^{(j)}(\omega)|.
\enn
Denote by $\mathcal C_{{\rm near}}^{(j)}(\omega)$ the set of index of modes those will be used to construct the few-mode DtR condition and in particular, the selected set $\mathcal C_{{\rm near}}^{(j)}(\omega)$ consists of these $n_{\min}^{(j)}$, the $r$ closest remaining propagating modes, and the $r$ closest remaining evanescent modes. The total number of indexes is denoted by $N_{\rm near}^{(j)}(\omega)=\operatorname{card}\big(\mathcal C_{{\rm near}}^{(j)}(\omega)\big)$. Define the remaining propagating and evanescent sets by
\ben
\mathcal P^{(j)}(\omega)
=\{n\ge0:\kappa>\alpha_n^{(j)}\}
\setminus\mathcal C_{{\rm near}}^{(j)}(\omega),
\quad
\mathcal Q^{(j)}(\omega)
=\{n\ge0:\kappa<\alpha_n^{(j)}\}
\setminus\mathcal C_{{\rm near}}^{(j)}(\omega).
\enn
The nearest-mode rule provides a frequency-independent separation for all unselected modes. Let
\ben
\beta_{\rm gap}^{(j)}:=\left[\frac12\inf_{n\ge 0}\left((\alpha_{n+1}^{(j)})^2-(\alpha_n^{(j)})^2\right)\right]^{1/2}=\frac{\pi}{\sqrt2h_j}.
\enn
The minimizing property gives, for every $n\ne n_{\min}^{(j)}$,
\ben
\left|\kappa^2-(\alpha_n^{(j)})^2\right|
\ge\frac12\left|(\alpha_n^{(j)})^2-(\alpha_{n_{\min}^{(j)}}^{(j)})^2\right|
\ge(\beta_{\rm gap}^{(j)})^2.
\enn
Therefore every unselected mode satisfies
\be\label{eq:beta_lower_bound}
|\beta_n^{(j)}(\omega)|
=\sqrt{|\kappa^2-(\alpha_n^{(j)})^2|}
\ge\beta_{\rm gap}^{(j)},
\quad n\notin\mathcal C_{{\rm near}}^{(j)}(\omega).
\en
This lower bound is independent of both $\omega$ and the distance to the cut-off frequencies.

For any function $f(\bm x)$ on $\Gamma_j^{\rm RPML}$ which admits the Fourier expansion $f(\bm x)=\sum_{n=0}^{\infty}f_n^{(j)}\psi_n^{(j)}(X_{\perp,j})$ with $f_n^{(j)}=\int_0^{h_j}f(\bm x)\overline{\psi_n^{(j)}(X_{\perp,j})}\,\mathrm dX_{\perp,j}$, $\mathcal T_{j}$, we introduce a \emph{few-mode DtR operator} on $\Gamma_j^{\rm RPML}$ defined by
\be\label{eq:few-mode-DtR}
\left(\mathcal T_{j}^{\rm fm,RPML}f\right)(\bm x,\omega)
=\sum_{n\in\mathcal C_{{\rm near}}^{(j)}(\omega)}\mathbf{i}\left(\beta_n^{(j)}(\omega)+\eta\right)f_n^{(j)}(\omega)\psi_n^{(j)}(X_{\perp,j}),
\en
which allows us to truncate the RPML region by imposing
\be\label{eq:few-mode-DtR-BC}
{\widetilde\partial_{\bm\nu}}\widetilde U^{\rm tot}+\mathbf{i}\eta\widetilde U^{\rm tot}-\mathcal T_{j}^{\rm fm,RPML}\widetilde U^{\rm tot}=0
\quad\text{on }\Gamma_j^{\rm RPML},\quad j=1,\ldots,N_p.
\en
This outer boundary condition implies that
\be\label{eq:DtR-effective-Neumann}
({\widetilde\partial_{\bm\nu}}\widetilde U^{\rm tot})_n^{(j)}+\mathbf{i}\eta(\widetilde U^{\rm tot})_n^{(j)}
=\begin{cases}
\mathbf{i}(\beta_n^{(j)}+\eta)(\widetilde U^{\rm tot})_n^{(j)},&n\in\mathcal C_{{\rm near}}^{(j)}(\omega),\\
0,&n\notin\mathcal C_{{\rm near}}^{(j)}(\omega).
\end{cases}
\en
which means that the Robin terms cancel on selected modes, giving the exact outgoing condition; unselected modes satisfy the homogeneous Robin condition.

The truncated RPML problem is then to find an approximate solution $\widetilde U^{\rm sca}_{\rm RPML}$ such that
\be\label{eq:truncated_RPML_DtR}
\begin{cases}
\mathscr L_j\widetilde U^{\rm sca}_{\rm RPML}=0,& \text{in }{\Omega^{\mathrm{tr}}},\\
{\widetilde\partial_{\bm\nu}}\widetilde U^{\rm sca}_{\rm RPML}+\mathbf{i}\eta_D\widetilde U^{\rm sca}_{\rm RPML}=\widetilde F_D^\star,& \text{on }\Gamma_D,\\
{\widetilde\partial_{\bm\nu}}\widetilde U^{\rm sca}_{\rm RPML}=\widetilde F^\star,& \text{on }{\Gamma_{\rm wall}^{\mathrm{tr}}},\\
{\widetilde\partial_{\bm\nu}}\widetilde U^{\rm sca}_{\rm RPML}+\mathbf{i}\eta\widetilde U^{\rm sca}_{\rm RPML}-\mathcal T_j^{\rm fm,RPML}\widetilde U^{\rm sca}_{\rm RPML}=\widetilde G_j,& \text{on }\Gamma_j^{\rm RPML},\ j=1,\ldots,N_p.
\end{cases}
\en
Here $\widetilde G_j=-{\widetilde\partial_{\bm\nu}}\widetilde U^{\rm inc}-\mathbf{i}\eta\widetilde U^{\rm inc}+\mathcal T_j^{\rm fm,RPML}\widetilde U^{\rm inc}$. For the variational formulation on the full computational domain, let $\varphi,\psi\in H^1({\Omega^{\mathrm{tr}}})$ and define
\ben
\begin{aligned}
b_{\rm RPML}(\varphi,\psi)
={}&\int_{\Omega_{\rm int}}
\left(\nabla\varphi\cdot\overline{\nabla\psi}
-\kappa^2\varphi\bar\psi\right)\,\mathrm d\bm x+\sum_{j=1}^{N_p}\int_{\Omega_j^{\rm RPML}}
\left[(\mathbb A_j^{-1}\nabla\varphi)\cdot\overline{\nabla\psi}
-\kappa^2\det(\mathbb B_j)\varphi\bar\psi\right]\,\mathrm d\bm x\\
&-\mathbf{i}\eta_D\int_{\Gamma_D}\varphi\bar\psi\,\mathrm ds+\sum_{j=1}^{N_p}\left[
\mathbf{i}\eta\int_{\Gamma_j^{\rm RPML}}\varphi\bar\psi\,\mathrm ds
-\left\langle\mathcal T_j^{\rm fm,RPML}\varphi,\psi\right\rangle_
{\Gamma_j^{\rm RPML}}\right],\\
\ell_{\rm RPML}(\psi)
={}&-\langle\widetilde F_D^\star,\psi\rangle_{\Gamma_D}
-\langle\widetilde F^\star,\psi\rangle_{{\Gamma_{\rm wall}^{\mathrm{tr}}}}
+\sum_{j=1}^{N_p}\langle\widetilde G_j,\psi\rangle_
{\Gamma_j^{\rm RPML}}.
\end{aligned}
\enn
The weak formulation of \eqref{eq:truncated_RPML_DtR} is to find
$\widetilde U^{\rm sca}_{\rm RPML}\in H^1({\Omega^{\mathrm{tr}}})$ such that
\be
\label{eq:weak_RPML0}
b_{\rm RPML}(\widetilde U^{\rm sca}_{\rm RPML},\psi)=\ell_{\rm RPML}(\psi)
\quad\forall\,\psi\in H^1({\Omega^{\mathrm{tr}}}).
\en

As illustrated in Figure~\ref{fig:rpml-dtr-maps}(b), we distinguish the exact DtR map $\mathcal T_j$ on $\Gamma_{L_j}$ defined in (\ref{eq:exact-DtR-op}), the few-mode DtR map $\mathcal T_{j}^{\rm fm,RPML}$ on $\Gamma_j^{\rm RPML}$ defined in (\ref{eq:few-mode-DtR}), and the interface map induced by the truncated RPML problem in $\Omega_j^{\rm RPML}$ which will be derived as follows. Since $\bm z\in\Omega_{\rm int}$, each RPML layer is source-free. For prescribed boundary values $f\in H^{1/2}(\Gamma_{L_j})$, let $\widetilde U^{\rm tot}_{\rm RPML}=\widetilde U^{\rm inc}+ \widetilde U^{\rm sca}_{\rm RPML}$ satisfy $\widetilde U^{\rm tot}_{\rm RPML}=f$ on $\Gamma_{L_j}$, the homogeneous Neumann condition on the waveguide walls within the layer $\Omega_j^{\rm RPML}$, and the outer condition \eqref{eq:few-mode-DtR-BC}. Define
\be\label{eq:fm-induced-definition}
\widehat{\mathcal T}_{j}^{\rm fm}f
:=\left.({\widetilde\partial_{\bm\nu}}\widetilde U^{\rm tot}_{\rm RPML}+\mathbf{i}\eta\widetilde U^{\rm tot}_{\rm RPML})\right|_{\Gamma_{L_j}}.
\en
For fixed layer geometry and parameters, $f$ and the outer condition will fix all modal coefficients of $\widetilde U^{\rm tot}_{\rm RPML}$ in $\Omega_j^{\rm RPML}$. 

Denote
\be\label{eq:Ttilde}
\widetilde T_j=T_j+ \mathbf{i} T_j^I,\quad T_j^I:=\int_{L_j}^{L_j+T_j}\sigma_j(t)\,\mathrm dt,
\en
\begin{lemma}\label{Lemma1}
It holds that $T_j^I = S_j T_j$.
\end{lemma}
\begin{proof}
Let $\Phi(s)=\frac{s^{P_j}}{s^{P_j}+(1-s)^{P_j}}$ for $s\in[0,1]$. Performing the change of variables $s=\frac{2t-(2L_j+T_j)}{T_j}$ maps $t\in[L_j,L_j+T_j]$ onto $s\in[-1,1]$ which yields
\ben
T_j^I= S_j T_j \int_{-1}^{1} \Phi(\chi_1(s))\,\mathrm{d}s.
\enn
Using the symmetry properties $\chi_1(-s)=1-\chi_1(s)$ and $\Phi(1-\tau)=1-\Phi(\tau)$, we obtain $\Phi(\chi_1(-s)) = 1 - \Phi(\chi_1(s))$. The change of variable $s \mapsto -s$ implies $\int_{-1}^{1} \Phi(\chi_1(s))\,ds = 2 - \int_{-1}^{1} \Phi(\chi_1(s))\,\mathrm{d}s$ and as a result, $\int_{-1}^{1} \Phi(\chi_1(s))\,ds = 1$. This completes the proof.
\end{proof}

By Lemma~\ref{Lemma1}, we obtain $\widetilde X_{1,j}(L_j+T_j)=L_j+\widetilde T_j$. For a fixed port $j$, the solution $\widetilde U^{\rm tot}_{\rm RPML}$ in $\Omega_j^{\rm RPML}$ satisfies a modal expansion as
\ben
\widetilde U^{\rm tot}_{\rm RPML}(\bm x,\omega)
=\sum_{n=0}^\infty \left(a_n^{(j)}\e^{\mathbf{i}\beta_n^{(j)}\widetilde X_{1,j}}+b_n^{(j)}\e^{-\mathbf{i}\beta_n^{(j)}\widetilde X_{1,j}}\right)\psi_n^{(j)}(X_{\perp,j}),\quad \bm x\in \Omega_j^{\rm RPML},
\enn
with some coefficients $a_n^{(j)},b_n^{(j)}$ to be determined. The boundary condition of $\widetilde U^{\rm tot}_{\rm RPML}$ on $\Gamma_j^{\rm RPML}$ gives the coefficient relation
\be\label{eq:DtR-coefficient-ratio1}
b_n^{(j)}=\begin{cases}
0,&n\in\mathcal C_{{\rm near}}^{(j)}(\omega),\\[1mm]
\displaystyle\frac{\beta_n^{(j)}+\eta}{\beta_n^{(j)}-\eta}\e^{2\mathbf{i}\beta_n^{(j)}(L_j+\widetilde T_j)}a_n^{(j)},&n\notin\mathcal C_{{\rm near}}^{(j)}(\omega).
\end{cases}
\en
In addition, from the boundary condition $\widetilde U^{\rm tot}_{\rm RPML}=f$ on $\Gamma_{L_j}$ it can be shown that
\be
\label{eq:DtR-coefficient-ratio2}
a_n^{(j)}\e^{\mathbf{i}\beta_n^{(j)}L_j}+b_n^{(j)}\e^{-\mathbf{i}\beta_n^{(j)}L_j}=f_n^{(j)}.
\en
A combination of (\ref{eq:DtR-coefficient-ratio1})-(\ref{eq:DtR-coefficient-ratio2}) yields
\be\label{eq:fm-layer-expansion}
\widetilde U^{\rm tot}_{\rm RPML}(\bm x,\omega)
=\sum_{n\in\mathcal C_{{\rm near}}^{(j)}(\omega)}f_n^{(j)}
\frac{\e^{\mathbf{i}\beta_n^{(j)}\widetilde X_{1,j}}}{\e^{\mathbf{i}\beta_n^{(j)}L_j}}\psi_n^{(j)}(X_{\perp,j})
+\sum_{n\notin\mathcal C_{{\rm near}}^{(j)}(\omega)}
\frac{\zeta_n^{(j)}(\widetilde X_{1,j})}{\zeta_n^{(j)}(L_j)}f_n^{(j)}\psi_n^{(j)}(X_{\perp,j}),\quad \bm x\in \Omega_j^{\rm RPML},
\en
where
\ben
\zeta_n^{(j)}(\widetilde X_{1,j})
=(\beta_n^{(j)}-\eta)\e^{-\mathbf{i}\beta_n^{(j)}(L_j+\widetilde T_j-\widetilde X_{1,j})}
+(\beta_n^{(j)}+\eta)\e^{\mathbf{i}\beta_n^{(j)}(L_j+\widetilde T_j-\widetilde X_{1,j})}.
\enn
Evaluating \eqref{eq:fm-induced-definition} using \eqref{eq:fm-layer-expansion} gives the series expansion of the induced interface DtR map as follows:
\be\label{eq:effective-RPML-DtR}
\widehat{\mathcal T}_{j}^{\rm fm}f
=\sum_{n\in\mathcal C_{{\rm near}}^{(j)}(\omega)}\mathbf{i}(\beta_n^{(j)}+\eta)f_n^{(j)}\psi_n^{(j)}
+\sum_{n\notin\mathcal C_{{\rm near}}^{(j)}(\omega)}
\frac{\mathbf{i}\left[(\beta_n^{(j)})^2-\eta^2\right]}
{\beta_n^{(j)}\coth(-\mathbf{i}\beta_n^{(j)}\widetilde T_j)-\eta}f_n^{(j)}\psi_n^{(j)}.
\en

Although the outer map $\mathcal T_{j}^{\rm fm,RPML}$ is finite-rank, the induced interface map $\widehat{\mathcal T}_{j}^{\rm fm}$ generally is not: unselected modes are attenuated by the RPML, satisfy the homogeneous Robin condition at its outer boundary, and also contribute to the interface map. Using the induced interface map $\widehat{\mathcal T}_{j}^{\rm fm}$, the RPML truncated problem (\ref{eq:truncated_RPML_DtR}) in $\Omega_j^{\rm RPML}$ can be reduced to an equivalent problem in $\Omega_{\mathrm{int}}$ as follows:
\be\label{eq:truncated_RPML_DtR1}
\begin{cases}
\Delta \widetilde U^{\rm sca}_{\rm RPML}+\kappa^2\widetilde U^{\rm sca}_{\rm RPML}=0,
& \text{in }\Omega_{\mathrm{int}},\\
{\partial_{\bm\nu}}\widetilde U^{\rm sca}_{\rm RPML}+\mathbf{i}\eta_D\widetilde U^{\rm sca}_{\rm RPML}=F_D^\star,& \text{on }\Gamma_D,\\
{\partial_{\bm\nu}}\widetilde U^{\rm sca}_{\rm RPML}=F^\star,& \text{on }{\Gamma_{\rm wall}^{\mathrm{t}}},\\
\partial_{\bm\nu}\widetilde U^{\rm sca}_{\rm RPML}+\mathbf{i}\eta \widetilde U^{\rm sca}_{\rm RPML}
-\widehat{\mathcal T}_{j}^{\rm fm}\widetilde U^{\rm sca}_{\rm RPML}= \widehat G_j, &\text{on }\Gamma_{L_j},\ j=1,\ldots,N_p,
\end{cases}
\en
with $\widehat G_j=-{\partial_{\bm\nu}}U^{\rm inc}-\mathbf{i}\eta U^{\rm inc}+\widehat{\mathcal T}_{j}^{\rm fm} U^{\rm inc}$. 

The next goal of this section lies in investigating the well-posedness of the reduced problem (\ref{eq:truncated_RPML_DtR1}) as well as the convergence of the approximate solution $\widetilde U^{\rm sca}_{\rm RPML}$ to the exact solution $U^{\rm sca}$ in the physical region $\Omega_{\mathrm{int}}$ with respect to the RPML parameters.

\subsection{Well-posedness and convergence analysis}
\label{sec:3.2}

We first derive an error estimate of the difference $\mathcal T_j-\widehat{\mathcal T}_{j}^{\rm fm}$. For any $\varphi\in H^{1/2}(\Gamma_{L_j})$, it holds that
\be\label{eq:fm-operator-error}
\mathcal E_j^{\rm fm}\varphi
:=\left(\mathcal T_j-\widehat{\mathcal T}_{j}^{\rm fm}\right)\varphi
=\sum_{n\notin\mathcal C_{{\rm near}}^{(j)}(\omega)}
\frac{2\mathbf{i}\beta_n^{(j)}(\beta_n^{(j)}+\eta)}
{(\beta_n^{(j)}-\eta)\e^{-2\mathbf{i}\beta_n^{(j)}\widetilde T_j}+\beta_n^{(j)}+\eta}\varphi_n^{(j)}\psi_n^{(j)}.
\en
To estimate this difference, we set
\ben
T_j^{\min}:=\min\left(T_j,T_j^I\right)=\min(1,S_j)T_j,
\enn
where the last equality follows from Lemma~\ref{Lemma1}, and we define
\ben
\beta_-^{(j)}(\omega):=\min_{n\in\mathcal P^{(j)}(\omega)}\beta_n^{(j)}(\omega),
\quad
\beta_+^{(j)}(\omega):=\min_{n\in\mathcal Q^{(j)}(\omega)}|\beta_n^{(j)}(\omega)|.
\enn

\begin{theorem}\label{thm:RPML-convergence}
For all $j=1,\cdots,N_p$ and $\varphi,\psi\in H^{1/2}(\Gamma_{L_j})$,
\be\label{eq:DtR_operator_bound}
\left|\left\langle\left(\mathcal T_j-\widehat{\mathcal T}_{j}^{\rm fm}\right)\varphi,\psi\right\rangle_{\Gamma_{L_j}}\right|
\leq M_j\,\|\varphi\|_{L^2(\Gamma_{L_j})}\|\psi\|_{L^2(\Gamma_{L_j})},
\en
where
\be\label{eq:Mj_definition}
M_j=\max\left(
\frac{2\beta_-^{(j)}(\omega)}
{\e^{2T_j^I\beta_-^{(j)}(\omega)}-1},
\frac{2\beta_+^{(j)}(\omega)}
{\e^{2 T_j\beta_+^{(j)}(\omega)}-1}
\right).
\en
Furthermore, the monotonicity of $\tau\mapsto2\tau/(\e^{2\sigma\tau}-1)$ and the lower bound estimate (\ref{eq:beta_lower_bound}) give
\be\label{eq:Mj_uniform}
M_j\le\frac{2\beta_{\rm gap}^{(j)}}{\e^{2\beta_{\rm gap}^{(j)}T_j^{\min}}-1}
=\mathcal O\!\left(\e^{-2\beta_{\rm gap}^{(j)}T_j^{\min}}\right)\quad\mbox{as}\quad T_j\to+\infty .
\en
in which the exponential convergence is independent of $\omega$ and its distance to the cut-off frequencies.
\end{theorem}
\begin{proof}
For $n\notin\mathcal C_{{\rm near}}^{(j)}(\omega)$, $\eta\le0$ implies $|\beta_n^{(j)}+\eta|\le|\beta_n^{(j)}-\eta|$, and
\ben
\left|\e^{-2\mathbf{i}\beta_n^{(j)}\widetilde T_j}\right|
=\begin{cases}
\e^{2\beta_n^{(j)}T_j^I},&n\in\mathcal P^{(j)}(\omega),\\
\e^{2|\beta_n^{(j)}|T_j},&n\in\mathcal Q^{(j)}(\omega).
\end{cases}
\enn
The reverse triangle inequality gives
\be\label{eq:fm-denominator-bound}
\begin{aligned}
\left|\e^{-\mathbf{i}\beta_n^{(j)}\widetilde T_j}\zeta_n^{(j)}(L_j)\right|
=\left|(\beta_n^{(j)}-\eta)\e^{-2\mathbf{i}\beta_n^{(j)}\widetilde T_j}+\beta_n^{(j)}+\eta\right|\ge |\beta_n^{(j)}-\eta|\left(\e^{2|\beta_n^{(j)}|T_j^{\min}}-1\right)>0,
\end{aligned}
\en
which also verifies that the denominators in \eqref{eq:fm-layer-expansion} are nonzero. Then the factor in (\ref{eq:fm-operator-error}) satifies
\ben
\left|\frac{2\mathbf{i}\beta_n^{(j)}(\beta_n^{(j)}+\eta)}
{(\beta_n^{(j)}-\eta)\e^{-2\mathbf{i}\beta_n^{(j)}\widetilde T_j}+\beta_n^{(j)}+\eta}\right|
\le\begin{cases}
\displaystyle\frac{2\beta_n^{(j)}}{\e^{2\beta_n^{(j)}T_j^I}-1},&n\in\mathcal P^{(j)}(\omega),\\[2mm]
\displaystyle\frac{2|\beta_n^{(j)}|}{\e^{2|\beta_n^{(j)}|T_j}-1},&n\in\mathcal Q^{(j)}(\omega).
\end{cases}
\enn
Since the function $\tau\mapsto2\tau/(\e^{2a\tau}-1)$ for $a,\tau>0$ is monotonically decreasing, we obtain
\ben
\left|\frac{2\mathbf{i}\beta_n^{(j)}(\beta_n^{(j)}+\eta)}
{(\beta_n^{(j)}-\eta)\e^{-2\mathbf{i}\beta_n^{(j)}\widetilde T_j}+\beta_n^{(j)}+\eta}\right|
\le M_j \quad\mbox{for all}\quad n\notin\mathcal C_{{\rm near}}^{(j)}(\omega).
\enn
Hence,
\ben
\left|\left\langle\mathcal E_j^{\rm fm}\varphi,\psi\right\rangle_{\Gamma_{L_j}}\right|
\le M_j\left(\sum_{n\notin\mathcal C_{N_{\rm near}}^{(j)}(\omega)}|\varphi_n^{(j)}|^2\right)^{1/2}\|\psi\|_{L^2(\Gamma_{L_j})}
\le M_j\|\varphi\|_{L^2(\Gamma_{L_j})}\|\psi\|_{L^2(\Gamma_{L_j})}.
\enn
The proof is complete.
\end{proof}

Let $b_{\rm fm}(\cdot,\cdot): H^1(\Omega_{\rm int})\times H^1(\Omega_{\rm int})\rightarrow \mathbb C$ be a sesquilinear form defined by
\be\label{eq:form_bRPML}
\begin{aligned}
b_{\rm fm}(\varphi,\psi)
=\int_{\Omega_{\rm int}}\left(\nabla\varphi\cdot\overline{\nabla\psi}-\kappa^2\varphi\bar\psi\right)\mathrm d\bm x-\mathbf{i}\eta_D\int_{\Gamma_D}\varphi\bar\psi\,\mathrm ds+\sum_{j=1}^{N_p}\left[\mathbf{i}\eta\int_{\Gamma_{L_j}}\varphi\bar\psi\,\mathrm ds-\left\langle\widehat{\mathcal T}_{j}^{\rm fm}\varphi,\psi\right\rangle_{\Gamma_{L_j}}\right].
\end{aligned}
\en
The variational formulation of the RPML truncated problem (\ref{eq:truncated_RPML_DtR1}) reads: find $\widetilde U^{\rm sca}_{\rm RPML}\in H^1(\Omega_{\rm int})$ such that
\be\label{eq:weak_RPML}
b_{\rm fm}(\widetilde U^{\rm sca}_{\rm RPML},\psi)=\ell_{\rm fm}(\psi)
\quad\forall\,\psi\in H^1(\Omega_{\rm int}),
\en
where 
\be\label{eq:fm-incident-functional}
\ell_{\rm fm}(\psi)=-\langle F_D^\star,\psi\rangle_{\Gamma_D}
-\langle F^\star,\psi\rangle_{\Gamma_{\rm wall}^{\rm t}}
+\sum_{j=1}^{N_p}\langle\widehat G_j,\psi\rangle_{\Gamma_{L_j}}
=\ell(\psi)-\sum_{j=1}^{N_p}\langle\mathcal E_j^{\rm fm}U^{\rm inc},\psi\rangle_{\Gamma_{L_j}}.
\en

The well-posedness of the variational problem (\ref{eq:weak_RPML}) and the convergence of $U^{\rm sca}-\widetilde U^{\rm sca}_{\rm RPML}$ is concluded in the following theorem. Note that, analogous to \cite[Lemma 2.3]{CW03}, there exists a constant $C>0$, depending only on $\Omega_{\rm int}$, such that
\be\label{eq:port-trace-bound}
\|\psi\|_{L^2(\Gamma_{L_j})}\le C\,\|\psi\|_{H^1(\Omega_{\rm int})},
\quad \psi\in H^1(\Omega_{\rm int}),\quad j=1,\ldots,N_p.
\en

\begin{theorem}\label{thm:RPML-error}
For any $\omega>0$ such that all $\beta_n^{(j)}(\omega)$ do not vanish, let $\gamma=\gamma(\omega)>0$ be the inf-sup constant in \eqref{eq:infsup} and $C$ be the constant in \eqref{eq:port-trace-bound}. If
\be\label{eq:smallness}
C^2\sum_{j=1}^{N_p}M_j<\gamma,
\en
then the variational problem (\ref{eq:weak_RPML}) has a unique solution $\widetilde U^{\rm sca}_{\rm RPML}\in H^1(\Omega_{\rm int})$ satisfying
\be\label{eq:RPML_error}
\|U^{\rm sca}-\widetilde U^{\rm sca}_{\rm RPML}\|_{H^1(\Omega_{\rm int})}
\le\frac{C}{\gamma}\sum_{j=1}^{N_p}M_j\|\widetilde U^{\rm tot}_{\rm RPML}\|_{L^2(\Gamma_{L_j})}.
\en
\end{theorem}
\begin{proof}
Let $u\in H^1(\Omega_{\rm int})$ satisfy $b_{\rm fm}(u,\psi)=0$ for all $\psi\in H^1(\Omega_{\rm int})$. By \eqref{eq:infsup}, Theorem~\ref{thm:RPML-convergence} and \eqref{eq:port-trace-bound},
\ben
\gamma\|u\|_{H^1(\Omega_{\rm int})}
\le \sup_{0\ne\psi\in H^1(\Omega_{\rm int})}
\frac{|b(u,\psi)|}{\|\psi\|_{H^1(\Omega_{\rm int})}}=
\sup_{0\ne\psi\in H^1(\Omega_{\rm int})}
\frac{\left|\sum_{j=1}^{N_p}
\langle\mathcal E_j^{\rm fm}u,\psi\rangle_{\Gamma_{L_j}}\right|}
{\|\psi\|_{H^1(\Omega_{\rm int})}}
\le C^2\sum_{j=1}^{N_p}M_j\|u\|_{H^1(\Omega_{\rm int})}.
\enn
Condition \eqref{eq:smallness} therefore gives $u=0$ which shows the uniqueness. Moreover, \eqref{eq:exact-DtR-energy} and Theorem~\ref{thm:RPML-convergence} imply
\ben
\operatorname{Re}b_{\rm fm}(\varphi,\varphi)
\ge
\|\varphi\|_{H^1(\Omega_{\rm int})}^2
-(1+\kappa^2)\|\varphi\|_{L^2(\Omega_{\rm int})}^2
-\sum_{j=1}^{N_p}M_j
\|\varphi\|_{L^2(\Gamma_{L_j})}^2.
\enn
which means that $b_{\rm fm}$ satisfies a G{\aa}rding's inequality. The well-posedness of the variational problem (\ref{eq:weak_RPML}) then follows directly from the standard Fredholm alternative.

Subtracting \eqref{eq:weak_RPML} from \eqref{eq:weak_original} and using \eqref{eq:fm-incident-functional} give
\be\label{eq:fm-error-identity}
b(U^{\rm sca}-\widetilde U^{\rm sca}_{\rm RPML},\psi)=\ell(\psi)-\ell_{\rm fm}(\psi)
+\sum_{j=1}^{N_p}
\left\langle\mathcal E_j^{\rm fm}
\widetilde U^{\rm sca}_{\rm RPML},\psi\right\rangle_{\Gamma_{L_j}}=\sum_{j=1}^{N_p}
\left\langle\mathcal E_j^{\rm fm}
\widetilde U^{\rm tot}_{\rm RPML},\psi\right\rangle_{\Gamma_{L_j}}.
\en
Consequently,
\ben
|b(U^{\rm sca}-\widetilde U^{\rm sca}_{\rm RPML},\psi)|
\le
C\sum_{j=1}^{N_p}M_j
\|\widetilde U^{\rm tot}_{\rm RPML}\|_{L^2(\Gamma_{L_j})}
\|\psi\|_{H^1(\Omega_{\rm int})}.
\enn
Applying \eqref{eq:infsup} proves \eqref{eq:RPML_error}.
\end{proof}


\begin{theorem}
\label{prop:uniqueness}
For any $\omega>0$ such that all $\beta_n^{(j)}(\omega)$ do not vanish, if the condition \eqref{eq:smallness} holds, the truncated RPML problem \eqref{eq:truncated_RPML_DtR} has a unique weak solution $\widetilde U^{\rm sca}_{\rm RPML}\in H^1({\Omega^{\mathrm{tr}}})$. Its restriction to $\Omega_{\rm int}$ is the solution of \eqref{eq:weak_RPML}.
\end{theorem}
\begin{proof}
Let $u\in H^1(\Omega^{\mathrm{tr}})$ solve the homogeneous problem of \eqref{eq:truncated_RPML_DtR}. Its restriction satisfies $b_{\rm fm}(u|_{\Omega_{\rm int}},\psi)=0$ for all
$\psi\in H^1(\Omega_{\rm int})$. Theorem~\ref{thm:RPML-error} gives $u=0$ in
$\Omega_{\rm int}$ and hence on each $\Gamma_{L_j}$. By \eqref{eq:fm-layer-expansion} and \eqref{eq:fm-denominator-bound}, $u=0$ in every RPML layer which shows the uniqueness. Set
\ben
c_0=\left(
1+\max_j\|\sigma_j\|_{L^\infty(L_j,L_j+T_j)}^2
\right)^{-1}>0.
\enn
The sesquilinear form $b_{\rm RPML}$ is bounded and, by \eqref{eq:few-mode-DtR}, satisfies the G{\aa}rding's inequality
\ben
\operatorname{Re}b_{\rm RPML}(\varphi,\varphi)
+(\kappa^2+c_0)
\|\varphi\|_{L^2(\Omega^{\mathrm{tr}})}^2
\ge
c_0\|\varphi\|_{H^1(\Omega^{\mathrm{tr}})}^2.
\enn
Existence then follows from uniqueness and the standard Fredholm alternative. This completes the proof.
\end{proof}

\subsection{RPML-BIE formulations}
\label{sec:3.3}

We now derive the RPML-BIE formulations for the truncated RPML problem (\ref{eq:truncated_RPML_DtR}). By \cite[Theorem~2.8]{LS01}, the RPML-transformed free-space fundamental solution under the complex-stretching map \eqref{eq:global_stretch} is given by
\ben
\widetilde{\Phi}_\omega(\bm x,\bm y)=\Phi_\omega(\widetilde{\bm x},\widetilde{\bm y})=\frac{\mathbf{i}}{4}H_0^{(1)}\!\left(\kappa\rho(\widetilde{\bm x},\widetilde{\bm y})\right),
\enn
where
\ben
\rho(\widetilde{\bm x},\widetilde{\bm y})=\left[(\widetilde x_1-\widetilde y_1)^2+(\widetilde x_2-\widetilde y_2)^2\right]^{1/2}.
\enn
Here the square root is taken with non-negative real part on $\mathbb C\setminus(-\infty,0]$. Although the RPML stretching \eqref{eq:global_stretch} is defined regionwise, it is continuous and agrees with the identity on each $\Gamma_{L_j},j=1,\cdots,N_p$. Hence $\widetilde{\Phi}_\omega$ satisfies the Helmholtz equation in $\Omega_{\rm int}$ and the
transformed Helmholtz equation in each $\Omega_j^{\rm RPML}$.

For notational convenience, denote the boundary components by
\ben
\Sigma_q=\begin{cases}\Gamma_D,&q=1,\\ {\Gamma_{\rm wall}^{\mathrm{tr}}},& q=2,\\ \Gamma_{q-2}^{\rm RPML},&3\le q\le N_p+2.\end{cases}.
\enn
For any given density function $\varphi:\Sigma_q\to\mathbb{C}, q=1,\cdots,N_p+2$ and $\bm x\in\Omega^{\rm tr}$, we define the RPML single- and double-layer potentials $\widetilde{\mathcal{S}}$ and $\widetilde{\mathcal{D}}$ by
\ben
\widetilde{\mathcal{S}}_{\Sigma_q}[\varphi](\bm x,\omega)=\int_{\Gamma_i} \widetilde{\Phi}_\omega(\bm x,\bm y)\varphi(\bm y,\omega)\mathrm{d}s_{\bm y},\quad\widetilde{\mathcal{D}}_{\Sigma_q}[\varphi](\bm x,\omega)=\int_{\Gamma_i} {\widetilde\partial_{\bm\nu_{\bm y}}}\widetilde{\Phi}_\omega(\bm x,\bm y)\varphi(\bm y,\omega)\mathrm{d}s_{\bm y}.
\enn
For $\bm x\in\Sigma_{q'}, q'=1,\cdots,N_p+2$, the associated RPML single-layer, double-layer, adjoint double-layer and hyper-singular BIOs are denoted by $\widetilde S_{\Sigma_{q'},\Sigma_q}$, $\widetilde D_{\Sigma_{q'},\Sigma_q}$, $\widetilde K_{\Sigma_{q'},\Sigma_q}$ and $\widetilde N_{\Sigma_{q'},\Sigma_q}$, respectively, and are defined as
\ben
\begin{aligned}
\widetilde{S}_{\Sigma_{q'},\Sigma_q}[\varphi](\bm x,\omega)&=\int_{\Sigma_q} \widetilde{\Phi}_\omega(\bm x,\bm y)\varphi(\bm y,\omega)\mathrm{d}s_{\bm y},\quad &
\widetilde{D}_{\Sigma_{q'},\Sigma_q}[\varphi](\bm x,\omega)&=\int_{\Sigma_q} {\widetilde\partial_{\bm\nu_{\bm y}}}\widetilde{\Phi}_\omega(\bm x,\bm y)\varphi(\bm y,\omega)\mathrm{d}s_{\bm y},\\
\widetilde{K}_{\Sigma_{q'},\Sigma_q}[\varphi](\bm x,\omega)&=\int_{\Sigma_q} {\widetilde\partial_{\bm\nu_{\bm x}}}\widetilde{\Phi}_\omega(\bm x,\bm y)\varphi(\bm y,\omega)\mathrm{d}s_{\bm y},\quad &
\widetilde{N}_{\Sigma_{q'},\Sigma_q}[\varphi](\bm x,\omega)&=\int_{\Sigma_q} {\widetilde\partial_{\bm\nu_{\bm x}}}{\widetilde\partial_{\bm\nu_{\bm y}}}\widetilde{\Phi}_\omega(\bm x,\bm y)\varphi(\bm y,\omega)\mathrm{d}s_{\bm y}.
\end{aligned}
\enn

Denoting $\mathcal M_j=\mathcal T_j^{\rm fm,RPML}-\mathbf{i}\eta I$, then ${\widetilde\partial_{\bm\nu}}\widetilde U^{\rm tot}_{\rm RPML}=\mathcal M_j\widetilde U^{\rm tot}_{\rm RPML}$ on $\Gamma_j^{\rm RPML}$ for $j=1,\cdots,N_p$. Then the Green's representation formulas for the problem (\ref{eq:truncated_RPML_DtR}) gives
\be\label{eq:point-source-rpml-representation}
\begin{aligned}
\widetilde U^{\rm sca}_{\rm RPML}(\bm x,\omega)=&\left({\widetilde{\mathcal D}_{\Gamma_D}}+\mathbf{i}\eta_D{\widetilde{\mathcal S}_{\Gamma_D}}\right)[\widetilde U^{\rm tot}_{\rm RPML}](\bm x,\omega)+{\widetilde{\mathcal D}_{\Gamma_{\rm wall}^{\mathrm{tr}}}}[\widetilde U^{\rm tot}_{\rm RPML}](\bm x,\omega)\\
&-\sum_{j=1}^{N_p}\left({\widetilde{\mathcal D}_{\Gamma_j^{\rm RPML}}}-{\widetilde{\mathcal S}_{\Gamma_j^{\rm RPML}}}\mathcal M_j\right)[\widetilde U^{\rm tot}_{\rm RPML}](\bm x,\omega),\quad \bm x\in{\Omega^{\mathrm{tr}}}.
\end{aligned}
\en
Taking the boundary traces and using the boundary conditions in (\ref{eq:truncated_RPML_DtR}) yield the following RPML-BIE system
\be\label{eq:rpml-bie-system}
-\frac{1}{2}\phi_{q'}
+\sum_{q=1}^{N_p+2}\mathbb T_{q',q}\phi_q
=R_{q'},
\quad q'=1,\ldots,N_p+2.
\en
where
\ben
\mathbb T_{q',q}=\begin{cases}
{\widetilde D_{\Sigma_{q'},\Sigma_1}}+\mathbf{i}\eta_D{\widetilde S_{\Sigma_{q'},\Sigma_1}},&q=1,\\
{\widetilde D_{\Sigma_{q'},\Sigma_q}},&q=2,\\
{\widetilde S_{\Sigma_{q'},\Sigma_q}}\mathcal M_{q-2}-{\widetilde D_{\Sigma_{q'},\Sigma_q}},&3\le q\le N_p+2,
\end{cases}\quad q'=1,\cdots,N_p+2,
\enn
and
\ben
\phi_q=\widetilde U^{\rm tot}_{\rm RPML}|_{\Sigma_q},\quad
R_{q'}=-\widetilde U^{\rm inc}_{\rm RPML}|_{\Sigma_{q'}},\quad q',q=1,\ldots,N_p+2.
\enn

In this work, we utilized the Chebyshev-based rectangular-polar solver developed in~\cite{LXYZ23} to achieve the numerical discretization of the RPML-stretched BIOs. Note that the few-mode DtR operators $\mathcal T_j^{\rm fm,RPML}$ on $\Gamma_j^{\rm RPML},j=1,\cdots,N_p$ are defined through mode expansions
\ben
\left(\mathcal T_{j}^{\rm fm,RPML}f\right)(\bm x,\omega)
=\sum_{n\in\mathcal C_{{\rm near}}^{(j)}(\omega)}\mathbf{i}\left(\beta_n^{(j)}(\omega)+\eta\right)f_n^{(j)}(\omega)\psi_n^{(j)}(X_{\perp,j}),
\enn
for any function $f(\bm x)=\sum_{n=0}^{\infty}f_n^{(j)}\psi_n^{(j)}(X_{\perp,j})$ on $\Gamma_j^{\rm RPML}$ with coefficients $f_n^{(j)}=\int_0^{h_j}f(\bm x)\overline{\psi_n^{(j)}(X_{\perp,j})}\,\mathrm dX_{\perp,j}$. Let $\{\bm X_k^{(j)}\}_{k=1}^{N_d^{(j)}}$, $\bm X_k^{(j)}=(L_j+T_j,X_{\perp,j}^{(k)})$ be the set of collocation nodes on $\Gamma_j^{\rm RPML}$ that used in Chebyshev-based rectangular-polar solver and denote the vector $\bm f^{(j)}=[\widehat f(\bm X_k^{(j)})]_{k=1}^{N_d^{(j)}}$ where we use the notation $\widehat f(\bm X^{(j)})= f(\bm x)$ for $\bm x=\in \Gamma_j^{\rm RPML}$. Then the numerical evaluation of the operator $\mathcal M_j$ can be proceeded as follows. For $n\in\mathcal C_{\rm near}^{(j)}(\omega)$, we define the discretization matrixes $\mathbb{P}_D^{(j)},\mathbb{P}_R^{(j)}$ with elements
\ben
(\mathbb{P}_D^{(j)})_{k,n}
=\psi_n^{(j)}(X_{\perp,j}^{(k)}),
\quad
(\mathbb{P}_R^{(j)})_{k,n}
=\mathbf{i}\big(\beta_n^{(j)}+\eta\big)
\psi_n^{(j)}(X_{\perp,j}^{(k)}),\quad k=1,\cdots,N_d^{(j)},\; n\in\mathcal C_{\rm near}^{(j)}(\omega)
\enn
Then the coefficients ${\bm F}^{(j)}=[f_n^{(j)}]_{n\in\mathcal C_{\rm near}^{(j)}(\omega)}$ can be obtained through the least-squares approach, i.e.,
\ben
\bm{F}^{(j)}=\left(\overline{\mathbb{P}_D^{(j)}}^\top \mathbb{P}_D^{(j)}\right)^{-1} \overline{\mathbb{P}_D^{(j)}}^\top \bm{f}^{(j)}.
\enn
Hence, the discretized version of the operator $\mathcal M_j$ takes the form
\ben
\mathbb{P}_R^{(j)}\left(\overline{\mathbb{P}_D^{(j)}}^\top \mathbb{P}_D^{(j)}\right)^{-1} \overline{\mathbb{P}_D^{(j)}}^\top -\mathbf{i}\eta \mathbb I,
\enn
where $\mathbb I$ is an identity matrix. 

\begin{remark}\label{remark:mode-incidence}
For mode incidence \eqref{eq:inc_mode}, the outgoing mode expansions apply to the scattered field in the incident port $W_{j_0}$ and to the total field in the other ports which means that the few-mode DtR conditions should be imposed on the RPML stretching of the scattered field in the incident port $W_{j_0}$ and on the RPML stretching of the total field in the other ports. Define
\ben
\mathcal J_\Sigma^{\rm inc}(\bm x,\omega)
=\widetilde{\mathcal D}_{\Sigma^{\rm inc}}[U^{\rm inc}](\bm x,\omega)
-\widetilde{\mathcal S}_{\Sigma^{\rm inc}}
[\partial_{\bm\nu}U^{\rm inc}](\bm x,\omega),\quad \bm x\notin\Sigma^{\rm inc}.
\enn
Let $\phi_q,q=1,\cdots,N_p+2$ be the boundary traces such that $\phi_1= \widetilde U^{\rm tot}_{\rm RPML}|_{\Gamma_D}$, $\phi_2|_{\partial \Omega^\inc\cap\Gamma_{\rm wall}^{\rm tr}}= \widetilde U^{\rm sca}_{\rm RPML}|_{\partial \Omega^\inc\cap\Gamma_{\rm wall}^{\rm tr}}$, $\phi_2|_{\Gamma_{\rm wall}^{\rm tr}\backslash\partial \Omega^\inc}= \widetilde U^{\rm tot}_{\rm RPML}|_{\Gamma_{\rm wall}^{\rm tr}\backslash\partial \Omega^\inc}$, $\phi_{j_0+2}=\widetilde U^{\rm sca}_{\rm RPML}|_{\Gamma_{j_0}^{\rm RPML}}$ and $\phi_{q'+2}=\widetilde U^{\rm tot}_{\rm RPML}|_{\Gamma_{q'}^{\rm RPML}}$ for $q'=1,\cdots,N_p, q'\ne j_0$. Then the Green's representation gives
\ben
\begin{aligned}
\widetilde U^{\rm sca}_{\rm RPML}(\bm x,\omega)
={}&
\left(
\widetilde{\mathcal D}_{\Gamma_D}
+\mathbf{i}\eta_D\widetilde{\mathcal S}_{\Gamma_D}
\right)[\phi_1](\bm x,\omega)
+\widetilde{\mathcal D}_{\Gamma_{\rm wall}^{\rm tr}}
[\phi_2](\bm x,\omega)\\
&-\sum_{j=1}^{N_p}
\left(
\widetilde{\mathcal D}_{\Gamma_j^{\rm RPML}}
-\widetilde{\mathcal S}_{\Gamma_j^{\rm RPML}}\mathcal M_j
\right)[\phi_{j+2}](\bm x,\omega)
-\mathcal J_\Sigma^{\rm inc}(\bm x,\omega)\quad\mbox{in}\quad\Omega^{\rm inc}\cap\Omega^{\rm tr}
\end{aligned}
\enn
and
\ben
\begin{aligned}
\widetilde U^{\rm tot}_{\rm RPML}(\bm x,\omega)
={}&
\left(
\widetilde{\mathcal D}_{\Gamma_D}
+\mathbf{i}\eta_D\widetilde{\mathcal S}_{\Gamma_D}
\right)[\phi_1](\bm x,\omega)
+\widetilde{\mathcal D}_{\Gamma_{\rm wall}^{\rm tr}}
[\phi_2](\bm x,\omega)\\
&-\sum_{j=1}^{N_p}
\left(
\widetilde{\mathcal D}_{\Gamma_j^{\rm RPML}}
-\widetilde{\mathcal S}_{\Gamma_j^{\rm RPML}}\mathcal M_j
\right)[\phi_{j+2}](\bm x,\omega)
-\mathcal J_\Sigma^{\rm inc}(\bm x,\omega)\quad \mbox{in}\quad \Omega^{\rm tr}\setminus\overline{\Omega^{\rm inc}}.
\end{aligned}
\enn
Taking boundary traces and applying the boundary conditions yield another BIE system with the same left-hand side of \eqref{eq:rpml-bie-system} and the right-hand side is given by
\ben
R_{q'}
=\left.\mathcal J_\Sigma^{\rm inc}\right|_{\Sigma_{q'}}
=\left.
\left(
\widetilde{\mathcal D}_{\Sigma^{\rm inc}}[U^{\rm inc}]
-\widetilde{\mathcal S}_{\Sigma^{\rm inc}}
[\partial_{\bm\nu}U^{\rm inc}]
\right)\right|_{\Sigma_{q'}},
\quad q'=1,\ldots,N_p+2.
\enn
\end{remark}

\begin{remark}\label{remark:dirichlet}
If the sound-hard condition on $\Gamma_{\rm wall}$ is replaced by the sound-soft condition $U^{\rm tot}=0$, the transverse eigensystem is replaced by
\ben
\psi_n^{(j)}(X_{\perp,j})
=\sqrt{\frac{2}{h_j}}
\sin\!\left(\frac{n\pi X_{\perp,j}}{h_j}\right),
\quad
\alpha_n^{(j)}=\frac{n\pi}{h_j},
\quad n\ge1.
\enn
The modal expansions now start at $n=1$, and the lower bound in \eqref{eq:beta_lower_bound} will be $\beta_{\rm gap}^{(j)}=\sqrt{3/2}\,\pi/h_j$. The few-mode DtR operators can be constructed in an analogouly manner using this new eigensystem.
\end{remark}

\section{Fourier transform algorithm}
\label{sec:4}

Once the the solutions of the frequency-domain problem (\ref{eq:truncated_RPML_DtR}) have been obtained through the frequency-robust RPML-BIE solver, the time domain solution $u^{\rm sca}(\bm x,t)\approx \widetilde u^{\rm tot}_{\rm RPML}(\bm x,t)$ in physical domain $\Omega_{\rm int}$ can be evaluated via inverse Fourier transform
\be
\label{eq:inverse-Fourier-sol}
\widetilde u^{\rm tot}_{\rm RPML}(\bm x,t)=\frac{1}{2\pi}\int_{W_-}^{W_+} \widetilde U^{\rm tot}_{\rm RPML}(\bm x,\omega)e^{-\mathbf{i}\omega t}\mathrm{d}\omega,\quad t>0.
\en
This section is devoted to proposing highly accurate Fourier transform algorithm to resolve the challenges (III) and (IV) raised in Section~\ref{sec:2.3}.

\subsection{Cut-off frequency singularity}
\label{sec:sine-SPFCC}

Due to the singularity of the waveguide Green's function near the cut-off frequencies, it can be expected that the frequency-domain solutions also suffer from the same type of singularity near the cut-off frequencies. In particular, let the unbounded waveguide $\Omega$ be the simplest planar plate with width $h>0$, i.e., $\Omega=\{\bm x=(x_1,x_2)\in\mathbb R^2: 0<x_2<h\}$. The background Green's function of the Helmholtz equation of this sound-hard planar waveguide is given by
\be\label{eq:exact-green}
G(\bm x,\bm z;\kappa)=\sum_{n\ge 0}\frac{\mathbf{i}}{2\beta_n}
\psi_n(x_2)\psi_n(z_2)\e^{\mathbf{i}\beta_n|x_1-z_1|},
\quad \beta_n=\begin{cases}
  \sqrt{\kappa^2-\alpha_n^2},&\mathrm{if}\alpha_n^2\leq \kappa^2,\\
  \mathbf{i}\sqrt{\alpha_n^2-\kappa^2},&\mathrm{if}~\alpha_n^2> \kappa^2.
\end{cases},\quad \bm x\ne\bm z,
\en
which clearly indicates a combined square-root singularity near the cut-off frequencies $n\pi/h,n\ge 0$. Then the solutions $\widetilde U^{\rm tot}_{\rm RPML}(\bm x,\omega)$ also admit such kind of square-root singularity which means that for sufficiently small $\xi$ with $\xi=\sqrt{\omega-\varpi}$ for $\omega>\varpi$ and $\xi=\mathbf{i}\sqrt{\varpi-\omega}$ for $\omega<\varpi$, $\widetilde U^{\rm tot}_{\rm RPML}(\bm x,\omega)$ can be asymptotically expanded as
\ben
\widetilde U^{\rm tot}_{\rm RPML}(\bm x,\omega)=\frac{a_{-1}(\bm x)}{\xi}
+\sum_{p=0}^{\infty}a_p(\bm x)\xi^p.
\enn
where $\varpi$ denotes a positive cut-off frequency. To effectively treat this singularity for the evaluation of the inverse Fourier transform (\ref{eq:inverse-Fourier-sol}), we introduce a sine-transform described as follows. 

We first partition the integration interval $[W_-,W_+]$ at the distinct port cut-off frequencies:
\ben
W_-=\varpi_0<\varpi_1<\cdots<\varpi_Q=W_+,
\quad q=1,\ldots,Q.
\enn
This gives
\ben
\widetilde u^{\rm tot}_{\rm RPML}(\bm x,t)=\sum_{q=1}^Q I_q(\bm x,t),\quad I_q(\bm x,t):=\frac{1}{2\pi}\int_{\varpi_{q-1}}^{\varpi_q} \widetilde U^{\rm tot}_{\rm RPML}(\bm x,\omega)e^{-\mathbf{i}\omega t}\mathrm{d}\omega,\quad t>0.
\enn
To regularize the square-root singularity at the sub-interval endpoints, we set
$\Delta_q=\varpi_q-\varpi_{q-1}$ and
$\bar\omega_q=(\varpi_q+\varpi_{q-1})/2$, and define the sine transform
\be\label{eq:sine_nodes}
\omega_q(u)=\bar\omega_q+\frac{\Delta_q}{2}
\sin\!\left(\frac{\pi u}{2}\right),
\quad -1\le u\le1.
\en
Its Jacobian is $\omega_q'(u)=(\Delta_q/2)p(u)$, where $p(u)=(\pi/2)\cos(\pi u/2)$. Then it can be verified that $G_q(\bm x,u)=p(u)\widetilde U^{\rm tot}_{\rm RPML}(\bm x,\omega_q(u))$ is smooth in $[\varpi_{q-1},\varpi_q]$ which further implies that $G_q(\bm x,u)$ can be expressed by means of the Chebyshev polynomial
\ben
G_q(\bm x,u)\approx \mathcal I_{N_\omega}G_q(\bm x,u)
=\sum_{m=0}^{N_\omega-1}c_{q,m}(\bm x)T_m(u)
\enn
where
\be\label{eq:fejer_nodes}
u_\ell=\cos\!\left(\frac{(2\ell-1)\pi}{2N_\omega}\right),
\quad \omega_{q,\ell}=\omega_q(u_\ell),
\quad \ell=1,\ldots,N_\omega.
\en
denotes the $N_\omega$ sampling Chebyshev nodes. Discrete orthogonality at the nodes $u_\ell$ gives
\ben
c_{q,m}(\bm x)=\frac{2-\delta_{m0}}{N_\omega}
\sum_{\ell=1}^{N_\omega}G_q(\bm x,u_\ell)T_m(u_\ell),
\quad m=0,\ldots,N_\omega-1,
\enn
where $\delta_{m0}$ is the Kronecker delta. To estimate the convergence of this expansion, we extend $G_q(\bm x,u)$ onto the complex plane $u\in\mathbb C$ and denote $E_\varrho$, $\varrho>1$ a Bernstein ellipse in $\mathbb C$ with foci $\pm1$ and semiaxes$(\varrho+\varrho^{-1})/2,(\varrho-\varrho^{-1})/2$. Set 
\ben
r_{q,q'}:=\begin{cases}
\frac{\varpi_{q-1}-\varpi_{q'}}{\Delta_q}, & \mbox{if}\; q'<q-1,\cr
\frac{\varpi_{q'}-\varpi_q}{\Delta_q}, & \mbox{if}\; q'>q.
\end{cases}
\enn
Then by solving the equation $\omega_q(u)=\varpi_{q'}$ it can be shown that the cut-off frequencies $\varpi_{q'}$ in the complex plane $u\in\mathbb C$ corresponding to the points
\ben
\label{eq:rho_sine}
u_{q,q'}=\begin{cases}
-1\pm\frac{2\mathbf{i}}{\pi}\operatorname{arcosh}(1+2r_{q,q'}),& \mbox{if}\; q'<q-1,\cr
1\pm\frac{2\mathbf{i}}{\pi}\operatorname{arcosh}(1+2r_{q,q'}), & \mbox{if}\; q'>q.
\end{cases}
\enn
Define $\varrho^*=\min_{q'}\max\{|u_{q,q'}+\sqrt{u_{q,q'}^2-1}|,\; |u_{q,q'}-\sqrt{u_{q,q'}^2-1}|\}$ satisfying $\varrho^*>1$. The convergence of the Chebyshev expansion is concluded in the following theorem.

\begin{theorem}\cite[Theorem~8.2]{T13}
\label{thm:cheb_geometric1}
Let $1<\varrho<\varrho^*$ and assume that $G_q(\bm x,u)$ is analytic on $\overline{E_\varrho}$. Then 
\be\label{eq:fejer_geometric_error}
\left\|G_q-\mathcal I_{N_\omega}G_q\right\|_{L^\infty(-1,1)}
\le\frac{4M_q(\varrho)}{\varrho-1}\varrho^{-(N_\omega-1)},
\en
where $M_q(\varrho)=\max_{z\in E_\varrho}|G_q(\bm x,z)|$.
\end{theorem}

Using the sine-transform and the Chebyshev expansion of $G_q$, the integration $I_q$ is approximated by
\be\label{eq:subband_integral}
I_q(\bm x,t)
=\frac{\Delta_q}{4\pi}\e^{-\mathbf{i}\bar\omega_qt}
\int_{-1}^{1}G_q(\bm x,u)
\e^{-\mathbf{i}\lambda_q\sin(\pi u/2)}\,\mathrm du\approx I_q^{N_\omega}(\bm x,t)
\en
where
\be\label{eq:fcc_rule}
I_q^{N_\omega}(\bm x,t)=\frac{\Delta_q}{4\pi}
\e^{-\mathbf{i}\bar\omega_qt}
\sum_{m=0}^{N_\omega-1}c_{q,m}(\bm x)\mu_m(\lambda_q),\quad \lambda_q=\frac{\Delta_qt}{2},
\en
with moments
\be\label{eq:fcc_moments}
\mu_m(\lambda_q)=\int_{-1}^{1}T_m(u)
\e^{-\mathbf{i}\lambda_q\sin(\pi u/2)}\,\mathrm du{\color{red}.}
\en
The nonlinear phase in \eqref{eq:fcc_moments} has stationary points at $u=\pm1$. This integration can be numerically evaluated via the SP-FCC algorithm described below. 

Using the change of variable $v=\sin(\pi u/2)$,  the nonlinear phase $\e^{-\mathbf{i}\lambda_q\sin(\pi u/2)}$ are transformed into the linear phase $\e^{-\mathbf{i}\lambda_q v}$ and we obtain the linear oscillation moment
\be\label{eq:fcc_stationary_transform}
\mu_m(\lambda_q)=\frac{2}{\pi}\int_{-1}^{1}
\frac{T_m\!\left(\frac{2}{\pi}\arcsin v\right)}{\sqrt{1-v^2}}
\e^{-\mathbf{i}\lambda_q v}\,\mathrm dv,
\en
which therefore introduces a square-root singularity at $v=\pm1$. Given the symmetry of $T_m$, it suffices to compute the integral over the interval $0 \le v \le 1$. Setting $s=1-v$ gives
\ben
\begin{aligned}
\mu_m(\lambda_q)=\mu_m^+(\lambda_q)+(-1)^m\overline{\mu_m^+(\lambda_q)},\quad\mu_m^+(\lambda_q)=\frac{2}{\pi}\int_0^1\frac{T_m\!\left(1-\frac{4}{\pi}\arcsin\sqrt{s/2}\right)}
{\sqrt{s(2-s)}}\e^{-\mathbf{i}\lambda_q(1-s)}\,\mathrm ds,\quad \lambda_q\in\mathbb R.
\end{aligned}
\enn
To treat the square-root singularity at $s=0$, a graded mesh $s_j=(j/M_F)^{q_F}$ is used with $q_F\ge1$, $j=0,\ldots,M_F$, which is refined near zero when $q_F$ increases. The first sub-interval $[0,s_1]$ contains the singular endpoint and is omitted as in \cite[Eq.~(3.5)]{DGK13}. Since the integrand is $\mathcal O(s^{-1/2})$, this introduces an error $\mathcal O(\sqrt{s_1})=\mathcal O(M_F^{-q_F/2})$ uniformly in real $\lambda_q$. On each remaining sub-interval $[s_{j-1},s_j]$, FCC quadrature approximates the smooth amplitude at the Clenshaw--Curtis nodes
\ben
s_{j,\ell}=s_{j-1}+\frac{s_j-s_{j-1}}{2}\bigl(1+\cos(\ell\pi/N_F)\bigr),\quad j=2,\ldots,M_F,\quad \ell=0,\ldots,N_F.
\enn
The local FCC weights are obtained by integrating the normalized Chebyshev polynomials with a linear oscillatory factor. When $0<|\lambda_q|(s_j-s_{j-1})<2$, the local oscillation is weak, and Clenshaw--Curtis quadrature is applied to the local integrand using the weights
\ben
w_\ell\!\left(\frac{\lambda_q(s_j-s_{j-1})}{2}\right)
=w_\ell(0)
\e^{\mathbf{i}\lambda_q(s_{j,\ell}-(s_j+s_{j-1})/2)}.
\enn
Analogous to the approach in \cite[Section~3]{DGS11}, the weights for the case $|\lambda_q|(s_j-s_{j-1})\ge2$ can be computed by
\ben
\begin{aligned}
w_\ell\!\left(\frac{\lambda_q(s_j-s_{j-1})}{2}\right)
={}&\frac{4}{N_F(s_j-s_{j-1})}
\sum_{n=0}^{N_F}{}'\cos(n\ell\pi/N_F)\\
&\quad\times
\int_{s_{j-1}}^{s_j}
T_n\!\left(\frac{2s-s_j-s_{j-1}}{s_j-s_{j-1}}\right)
\e^{\mathbf{i}\lambda_q(s-(s_j+s_{j-1})/2)}
\,\mathrm ds.
\end{aligned}
\enn
A prime on a sum means that its first and last terms are divided by 2. Combining the above local quadrature rules over
$j=2,\ldots,M_F$ gives
\ben
\begin{aligned}
\mu_m^+(\lambda_q)
\approx\frac{1}{\pi}\sum_{j=2}^{M_F}(s_j-s_{j-1})
\e^{-\mathbf{i}\lambda_q(1-(s_j+s_{j-1})/2)}
\sum_{\ell=0}^{N_F}{}'
w_\ell\!\left(\frac{\lambda_q(s_j-s_{j-1})}{2}\right)
\frac{T_m\!\left(1-\frac{4}{\pi}\arcsin\sqrt{s_{j,\ell}/2}\right)}
{\sqrt{s_{j,\ell}(2-s_{j,\ell})}}.
\end{aligned}
\enn
Substituting their numerical values into
\eqref{eq:fcc_rule} together with the symmetry relation of $\mu_m$ and summing over the frequency sub-intervals
yields the inverse Fourier transform
\ben
\widetilde u^{\rm tot}_{\rm RPML}(\bm x,t)
\approx \sum_{q=1}^Q\frac{\Delta_q}{4\pi}
\e^{-\mathbf{i}\bar\omega_q t}
\sum_{m=0}^{N_\omega-1}c_{q,m}(\bm x)\mu_m(\lambda_q).
\enn

\subsection{Complex resonance extraction and singularity subtraction}
\label{sec:rational-SPFCC}

Note that Theorem~\ref{thm:cheb_geometric1} requires that the transformed solutions $G_q(\bm x,u)$ to be analytic on $\overline{E_\varrho}$, however, those complex resonances close to the real axis might either destroy this assumption or restrict to choose smaller $\rho>1$ and as a result, a dense sampling of frequency points is necessary to achieve high accuracy. In light of the singularity subtraction technique proposed in \cite{BS26}, here we propose a new methodology to resolve this issue for the current waveguide problems by utilizing the AAA method~\cite{NST18} to estimate the locations and residues of selected complex resonances within $\overline{E_\varrho}$.

\begin{figure}[htbp]
\centering
\includegraphics[width=0.75\textwidth]{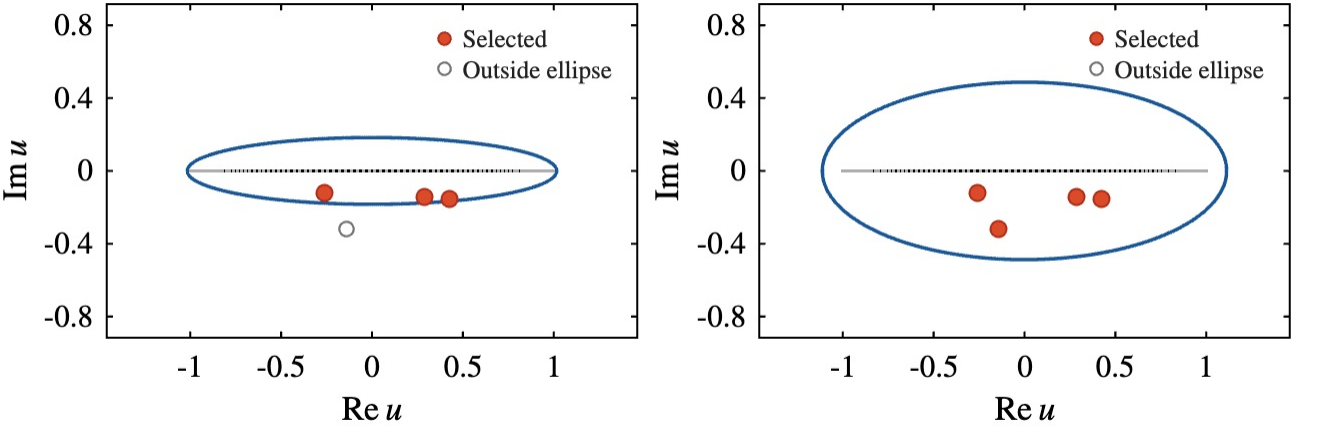}
\caption{Selection of complex resonances estimated by the AAA method: (a) $\varrho=1.2$; (b) $\varrho=1.6$. Filled and open circles denote selected and unselected complex resonances.}
\label{fig:bernstein-pole-selection}
\end{figure}

For a fixed sub-interval $[\varpi_{q-1},\varpi_q]$, let $V_q^{\rm r}(u)$ be the AAA rational approximation to the total-field solution $\widetilde U^{\rm tot}_{\rm RPML}(\bm x,\omega_q(u))$, obtained from its values at the $N_\omega$ Chebyshev nodes in \eqref{eq:fejer_nodes}:
\be\label{eq:aaa-weight}
V_q^{\rm r}(u)=
\Bigl[\sum_{j=1}^{N_{\rm s}}
\vartheta_j\widetilde U^{\rm tot}_{\rm RPML}(\bm x,\hat\omega_{q,j})/(u-\hat u_{j})\Bigr]
\Big/\Bigl[\sum_{j=1}^{N_{\rm s}}\vartheta_j/(u-\hat u_{j})\Bigr].
\en
Here $1\le N_{\rm s}\le N_\omega$ is the number of support nodes $\hat u_j$ selected from these Chebyshev nodes, $\hat\omega_{q,j}=\omega_q(\hat u_j)$, and $\vartheta_j$ are the barycentric weights. Among the remaining Chebyshev nodes, support nodes are added sequentially in descending order of approximation error. For a given relative fitting tolerance $\epsilon_{\rm A}>0$, the iteration stops when the following condition is satisfied
\ben
\max_{1\le\ell\le N_\omega}
\left|
\widetilde U^{\rm tot}_{\rm RPML}
(\bm x,\omega_{q,\ell})
-V_q^{\rm r}(u_\ell)
\right|
\le
\epsilon_{\rm A}
\max_{1\le\ell\le N_\omega}
\left|
\widetilde U^{\rm tot}_{\rm RPML}
(\bm x,\omega_{q,\ell})
\right|,
\enn
or when the fixed upper bound of $N_s$ is reached. To estimate the locations $\widehat u_b$ in the complex $u$-plane corresponding to the complex resonances $\widehat\omega_b=\omega_q(\widehat u_b)$, we set
\ben
\sum_{j=1}^{N_{\rm s}}
\frac{\vartheta_j}{\widehat u_b-\hat u_j}=0,\quad \mathrm{Im}\,\omega_q(\widehat u_b)<0.
\enn
In particular, a necessary condition to ensure $\mathrm{Im}\,\omega_q(\widehat u_b)<0$ is $\mathrm{Im}\,u<0$. Since a sine-transform is utilized to treat the square-root singularity near the cut-off frequencies, we additionally impose that the denominator has a nonzero derivative at $\widehat u_b$ and that the numerator is nonzero at those $\widehat u_b$, i.e.,
\ben
-\sum_{j=1}^{N_{\rm s}}
\frac{\vartheta_j}{(\widehat u_b-\hat u_j)^2}\ne0,
\quad
\sum_{j=1}^{N_{\rm s}}
\frac{\vartheta_j\widetilde U^{\rm tot}_{\rm RPML}
(\bm x,\hat\omega_{q,j})}{\widehat u_b-\hat u_j}\ne0, 
\enn
Then the residue corresponding to $\widehat\omega_b$ is
\be\label{eq:nonzero-pole-extraction}
\widehat R_b(\bm x)=-\omega_q'(\widehat u_b)
\Bigl[\sum_{j=1}^{N_{\rm s}}
\vartheta_j\widetilde U^{\rm tot}_{\rm RPML}
(\bm x,\hat\omega_{q,j})/(\widehat u_b-\hat u_j)\Bigr]
\Big/\Bigl[\sum_{j=1}^{N_{\rm s}}
\vartheta_j/(\widehat u_b-\hat u_j)^2\Bigr].
\en
According to Theorem~\ref{thm:cheb_geometric1}, we select those $\widehat u_b$ lying inside the chosen Bernstein ellipse $\overline{E_\varrho}$ for a given $\varrho>1$. Figure~\ref{fig:bernstein-pole-selection} presents an illustration for the selection of $\widehat u_b$ in the complex $u$-plane. 

Subtracting the singular terms, it is natural to obtain the residual regular term
\be\label{eq:nonzero-pole-subtraction}
\widetilde U_{\rm reg}(\bm x,\omega)
=\widetilde U^{\rm tot}_{\rm RPML}(\bm x,\omega)
-\sum_b\frac{\widehat R_b(\bm x)}{\omega-\widehat\omega_b}.
\en
Then the inverse Fourier integral $I_q$ on $[\varpi_{q-1},\varpi_q]$ can be split as $I_q=I_{1,q}+I_{2,q}$, where
\be\label{eq:regular-sub-interval-integral}
I_{1,q}(\bm x,t)=\frac{1}{2\pi}
\int_{\varpi_{q-1}}^{\varpi_q}
\widetilde U_{\rm reg}(\bm x,\omega)
\e^{-\mathbf{i}\omega t}\,\mathrm d\omega
\en
and
\be\label{eq:nonzero-pole-time}
I_{2,q}(\bm x,t)=\frac{1}{2\pi}\sum_b\widehat R_b(\bm x)
\int_{\varpi_{q-1}}^{\varpi_q}
\frac{\e^{-\mathbf{i}\omega t}}{\omega-\widehat\omega_b}\,\mathrm d\omega.
\en
Replacing $\widetilde U^{\rm tot}_{\rm RPML}(\bm x,\omega_q(u))$ by $\widetilde U_{\rm reg}(\bm x,\omega_q(u))$ in the Chebyshev interpolation and SP-FCC quadrature rule in Section~\ref{sec:sine-SPFCC} gives the following approximation to $I_{1,q}(\bm x,t)$:
\be\label{eq:regular-SPFCC}
I_{1,q}(\bm x,t)
\approx I_{1,q}^{N_\omega}(\bm x,t)
=\frac{\Delta_q}{4\pi}\e^{-\mathbf{i}\bar\omega_q t}
\sum_{m=0}^{N_\omega-1}c_{q,m}(\bm x)\mu_m(\lambda_q),
\en
For $I_{2,q}$, the near-singularity can be removed using the identity \cite[Section~4.2.1]{BS26}
\be\label{eq:finite-band-resonance-integral}
\int_{\varpi_{q-1}}^{\varpi_q}
\frac{\e^{-\mathbf{i}\omega t}}{\omega-\widehat\omega_b}\,\mathrm d\omega
=\int_{\varpi_{q-1}}^{\varpi_q}
\frac{\e^{-\mathbf{i}\omega t}-\e^{-\mathbf{i}\widehat\omega_b t}}
{\omega-\widehat\omega_b}\,\mathrm d\omega+\e^{-\mathbf{i}\widehat\omega_b t}
\left[\log(\varpi_q-\widehat\omega_b)
-\log(\varpi_{q-1}-\widehat\omega_b)\right].
\en
After subtraction, the integrand $(\e^{-\mathbf{i}\omega t}-\e^{-\mathbf{i}\widehat\omega_b t})/(\omega-\widehat\omega_b)$ is smooth, with the limit $-\mathbf{i}t\e^{-\mathbf{i}\widehat\omega_b t}$ at $\omega=\widehat\omega_b$. Hence $I_{2,q}$ can be numerical evaluated by using the classical Fejer's quadrature rule based on the same $N_\omega$ Chebyshev nodes and we denote $I_{2,q}^{N_\omega}$ the approximation. Combining the two contributions gives
\be\label{eq:nonzero-pole-reconstruction}
\widetilde u^{\rm tot}_{\rm RPML}(\bm x,t)
\approx\sum_{q=1}^Q\left[I_{1,q}^{N_\omega}(\bm x,t)+I_{2,q}^{N_\omega}(\bm x,t)\right].
\en

In conclusion, the numerical evaluation of the time-domain solution proceeds through the following sequence:
\ben
\begin{aligned}
f_B(\bm x,t)\xrightarrow{\mathbb F}F_B(\bm x,\omega)
&\xrightarrow[\text{Section~\ref{sec:3}}]
{\text{RPML-BIE with few-mode DtR}}
\widetilde U^{\rm tot}_{\rm RPML}
(\bm x,\omega_{q,\ell})\\
& \xrightarrow[\text{Section~\ref{sec:4}}]
{\text{AAA-based singularity subtraction + SP-FCC}}
\widetilde u^{\rm tot}_{\rm RPML}(\bm x,t).
\end{aligned}
\enn


\section{Numerical experiments}
\label{sec:5}

In this section, we present several numerical examples to demonstrate the efficiency and accuracy of the proposed frequency-time hybrid RPML-BIE solver for solving time-dependent wave-equation problems in two-dimensional nonuniform closed waveguides with
multiple channels. The resulting linear system is solved using the fully complex version of the iterative solver GMRES with a residual tolerance of $\epsilon_r=10^{-14}$. All numerical results presented in this paper were obtained using Fortran implementations of the algorithms parallelized with OpenMP. In all examples, we set $\eta_D=1$, $\eta=-1$, $r=2$, $S_j=1$ and $P_j=8$ for $j=1,\cdots,N_p$. The parameter $T_j,j=1,\cdots,N_p$ for different ports are set to be the same and we denote $T_j=T,j=1,\cdots,N_p$. For the incident Gauss pulse~\eqref{eq:td_spectrum}, we set $\omega_0=3.5\pi$, $\sigma=\sqrt{2}$, $t_0=4$ and $A_0=10$ which allows us to truncate the frequency interval as $[W_-,W_+]=[0.5\pi,6.5\pi]$.

\begin{figure}[!htbp]
\centering
\setlength{\tabcolsep}{1pt}
\begin{tabular}{cccc}
\subfloat[two-port]{\includegraphics[width=0.235\textwidth]{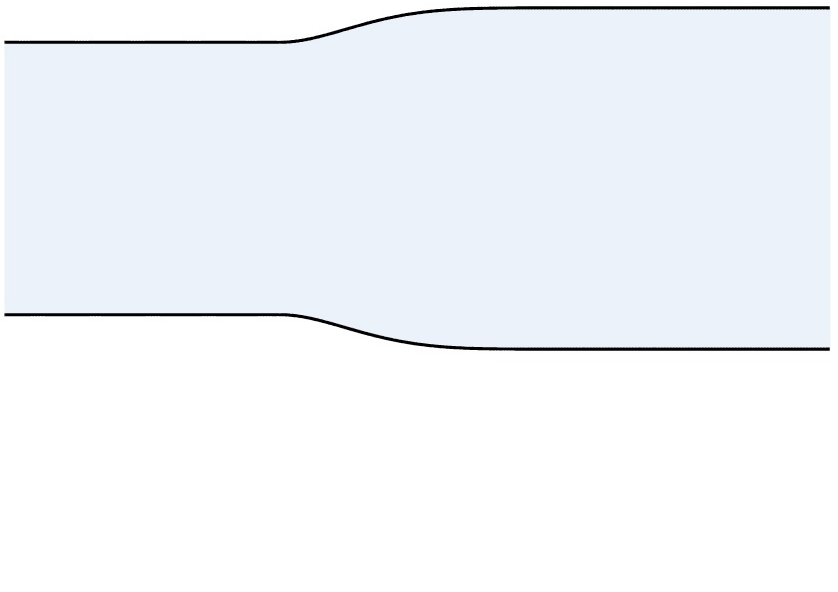}}&
\subfloat[parallel three-port]{\includegraphics[width=0.235\textwidth]{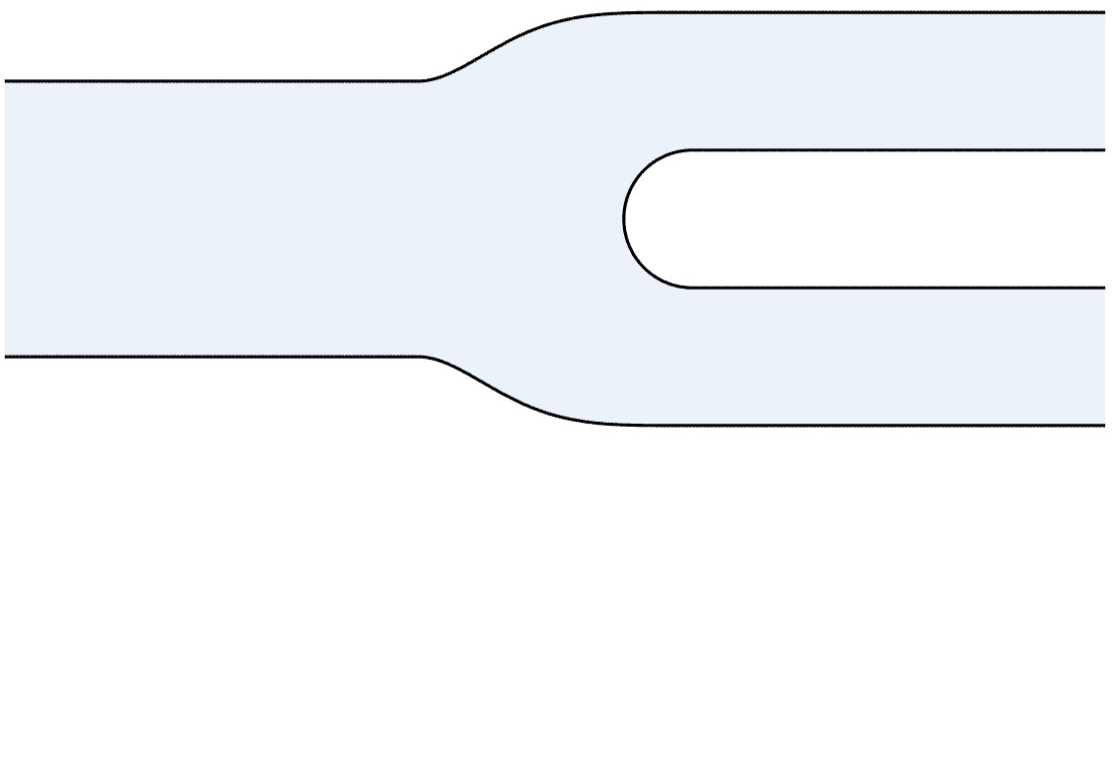}}&
\subfloat[T-junction]{\includegraphics[width=0.235\textwidth]{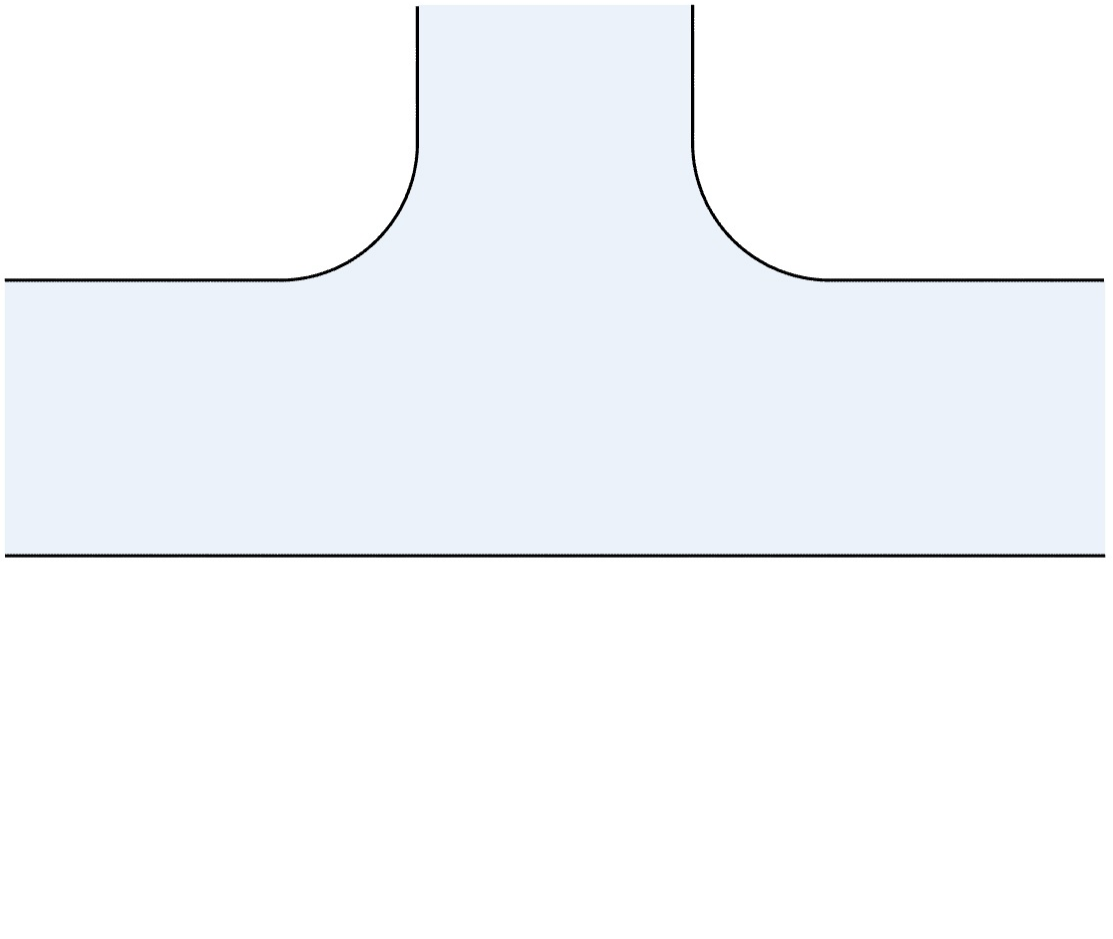}}&
\subfloat[Y-junction]{\includegraphics[width=0.235\textwidth]{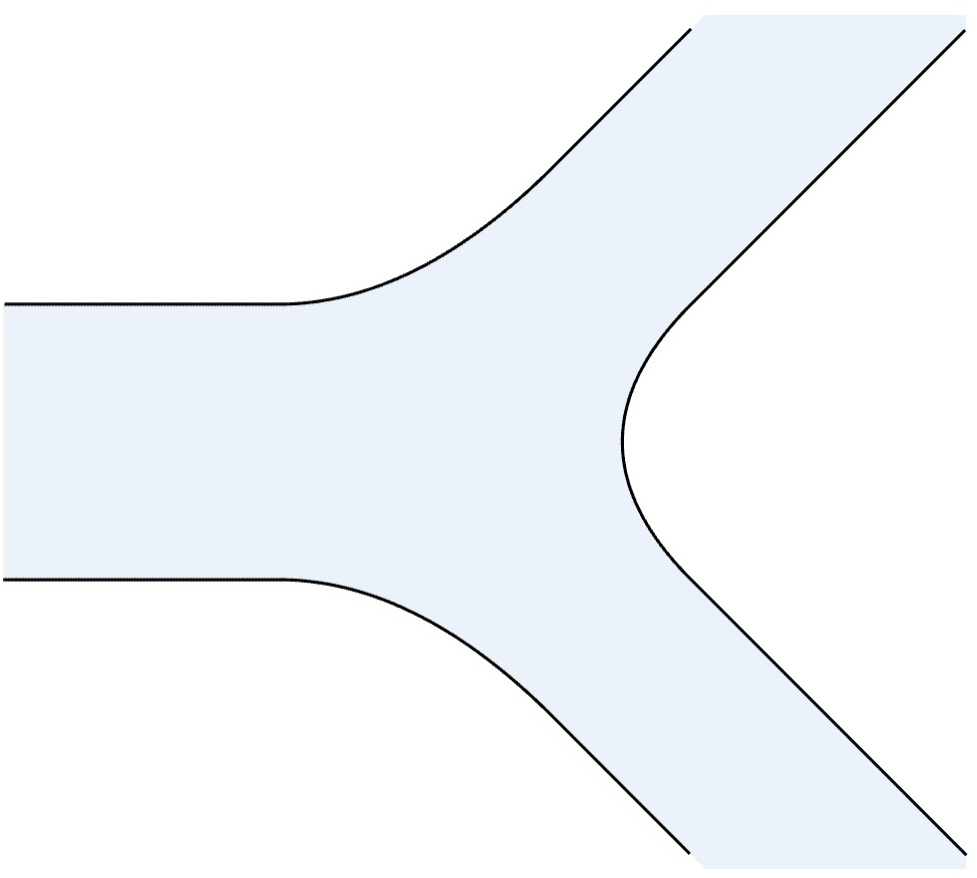}}
\end{tabular}
\caption{Waveguide geometries in numerical tests.}
\label{fig:numerical-geometries}
\end{figure}

To illustrate the performance of the method, we consider four kinds of nonuniform waveguides shown in Figure~\ref{fig:numerical-geometries} each of which may contains a bounded obstacle with an impedance condition. The relative error in frequency domain is defined as
\be\label{eq:FD-relative-error}
\varepsilon_{\rm FD}=\frac{\max_{\bm x}|U_{\rm num}(\bm x_{\rm ob})-U_{\rm ref}(\bm x_{\rm ob})|}
{\max_{\bm x_{\rm ob}}|U_{\rm ref}(\bm x_{\rm ob})|},
\en
where $\bm x_{\rm ob}$ denotes the observation points in the physical domain $\Omega_{\rm int}$, $U_{\rm num}$ denotes the numerical solution resulting from the RPML-BIE method and $U_{\rm ref}$ denotes the reference solution which is obtained from an exact representation when available or from a sufficiently refined numerical solution. For the time-domain tests, the numerical error at one observation point $\bm x_{\rm ob}$ is defined as 
\be\label{eq:errors}
\varepsilon_{\rm TD}(\bm x_{\rm ob},t)=|u_{\rm num}(\bm x_{\rm ob},t)-u_{\rm ref}(\bm x_{\rm ob},t)|,
\en
where $u_{\rm num}$ denotes the numerical solution resulting from the proposed frequency-time hybrid RPML-BIE solver and $u_{\rm ref}$ denotes the time-domain reference solution.

\FloatBarrier

\begin{figure}[!ht]
\centering
\setlength{\tabcolsep}{1pt}
\begin{tabular}{cc}
\subfloat[Sound-hard case.]{\includegraphics[width=0.3\textwidth]{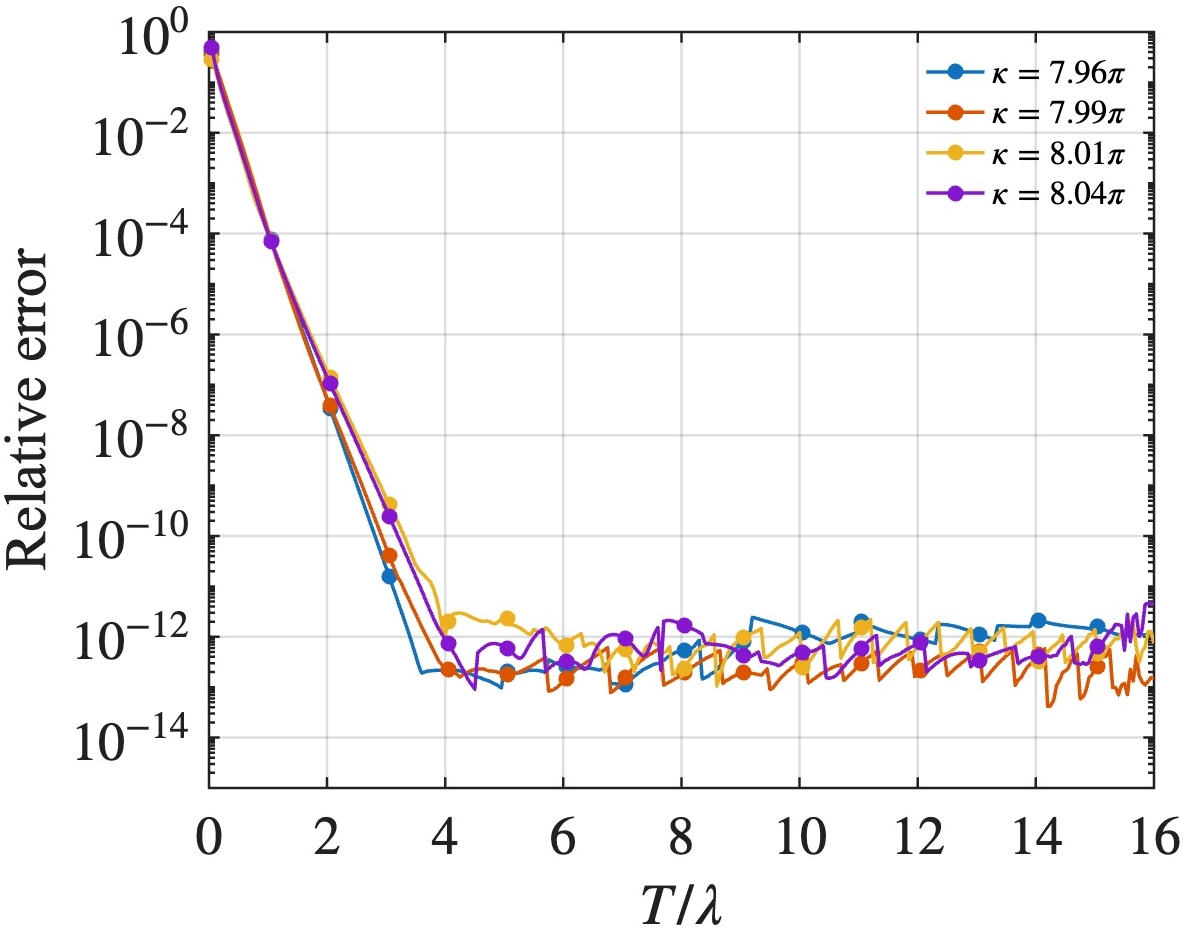}}&
\subfloat[Sound-soft case.]{\includegraphics[width=0.3\textwidth]{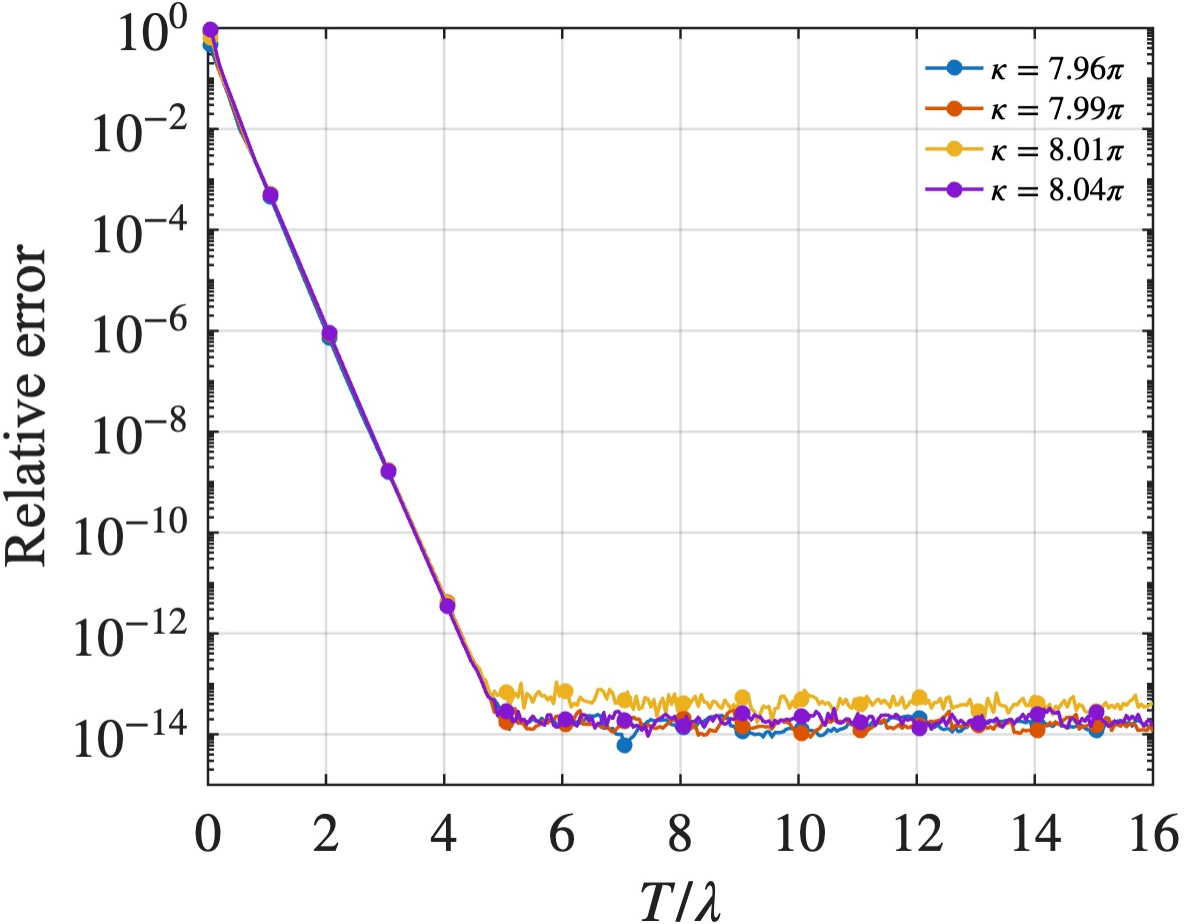}}
\end{tabular}
\caption{Example 1. Relative errors at different frequencies of the few-mode corrected RPML-BIE method for the planar waveguide problems.}
\label{fig:error-exact}
\end{figure}

{\bf Example~1} ({\it frequency-domain test for a planar waveguide}). We first consider a point-source incidence placed at $\bm z=(-2,1.3)$ in an obstacle-free planar waveguide with width $h=2$ and sound-hard or sound-soft walls. The exact total field for sound-hard case is given by the modal Green function \eqref{eq:exact-green} and the exact for sound-soft case can be expressed analogously. We choose the observation points $\bm x_{\rm ob}=(\pm 0.5,\pm 0.5)$. Figure~\ref{fig:cutoff}(a) shows the destroyed convergence of the uncorrected RPML-BIE method near the cut-off frequencies, whereas the results displayed in Figure~\ref{fig:error-exact} verify that the proposed few-mode corrected RPML-BIE method can uniformly achieve exponential convergence with respect to $T$ and high accuracy near the cut-off frequencies.

\begin{figure}[!ht]
\centering
\captionsetup[subfloat]{labelformat=empty}
\setlength{\tabcolsep}{1pt}
\begin{tabular}{cccc}
\subfloat[(a)]{\includegraphics[width=0.2\textwidth]{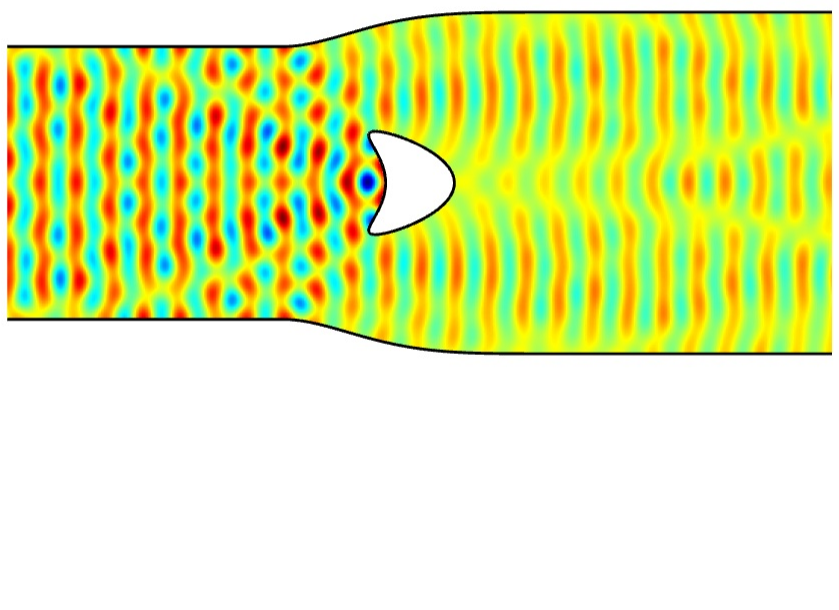}}&
\subfloat[(b)]{\includegraphics[width=0.2\textwidth]{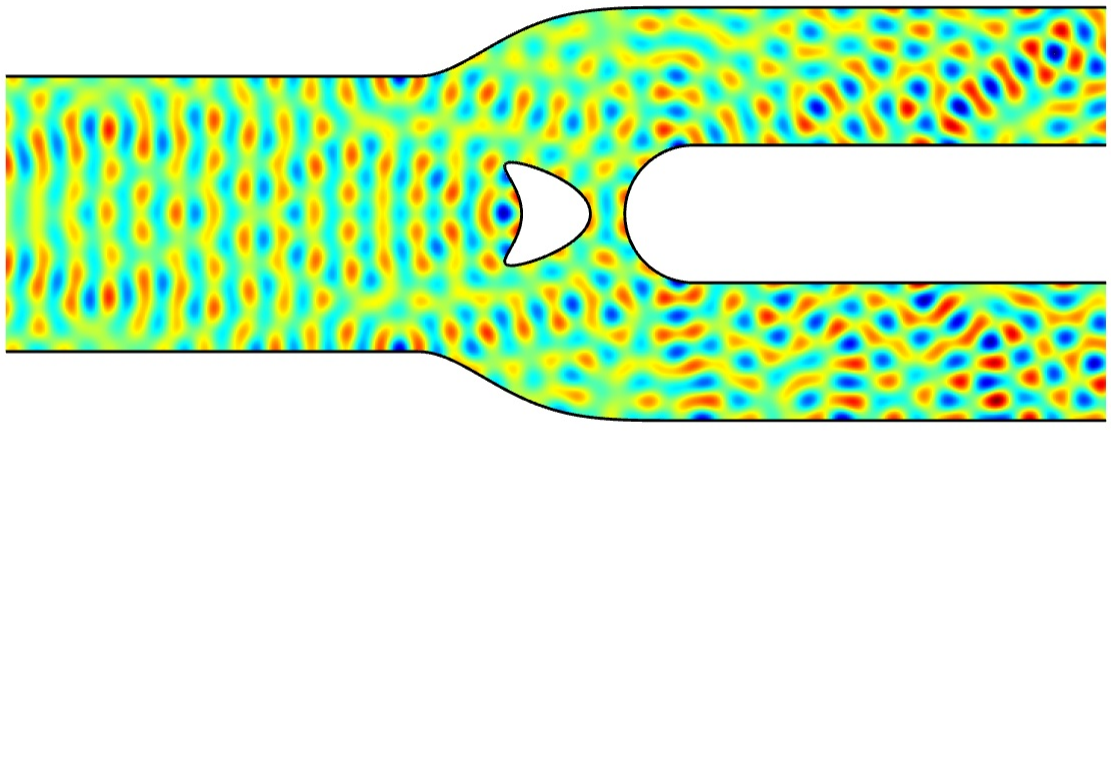}}&
\subfloat[(c)]{\includegraphics[width=0.2\textwidth]{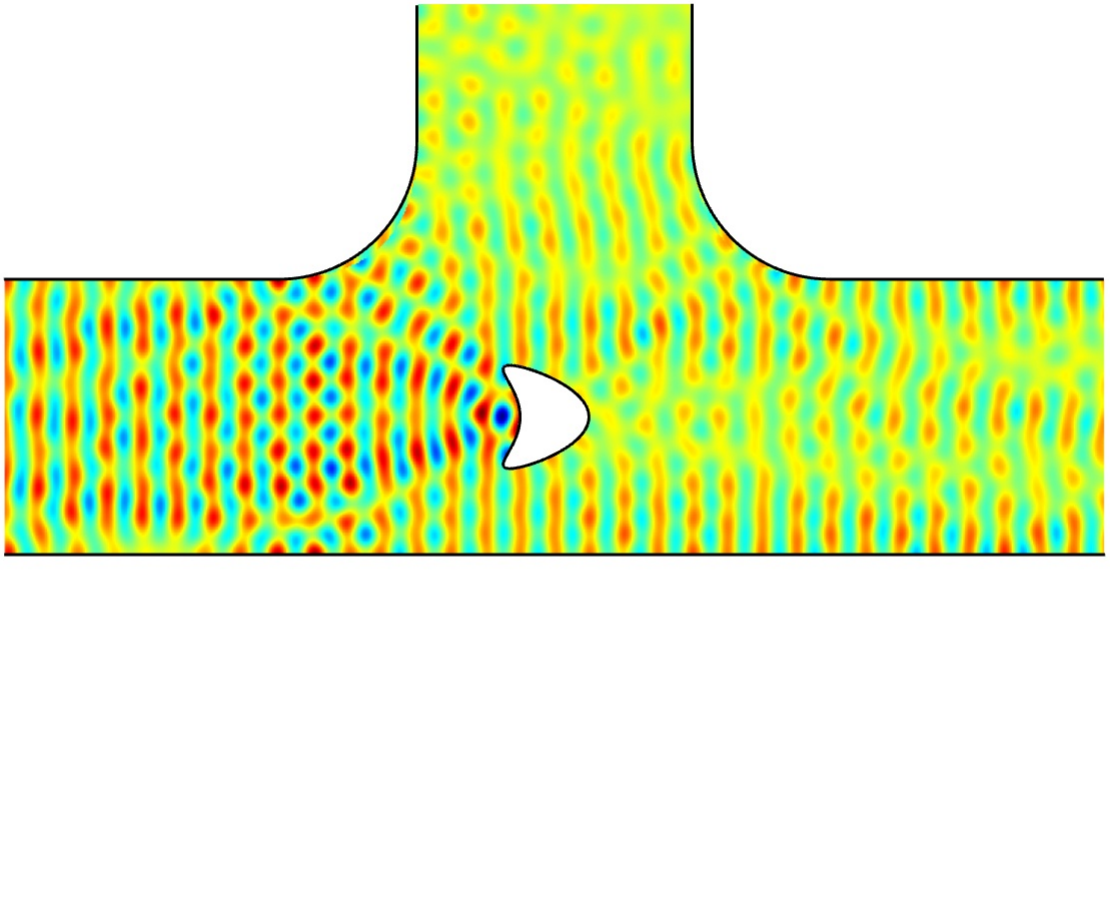}}&
\subfloat[(d)]{\includegraphics[width=0.2\textwidth]{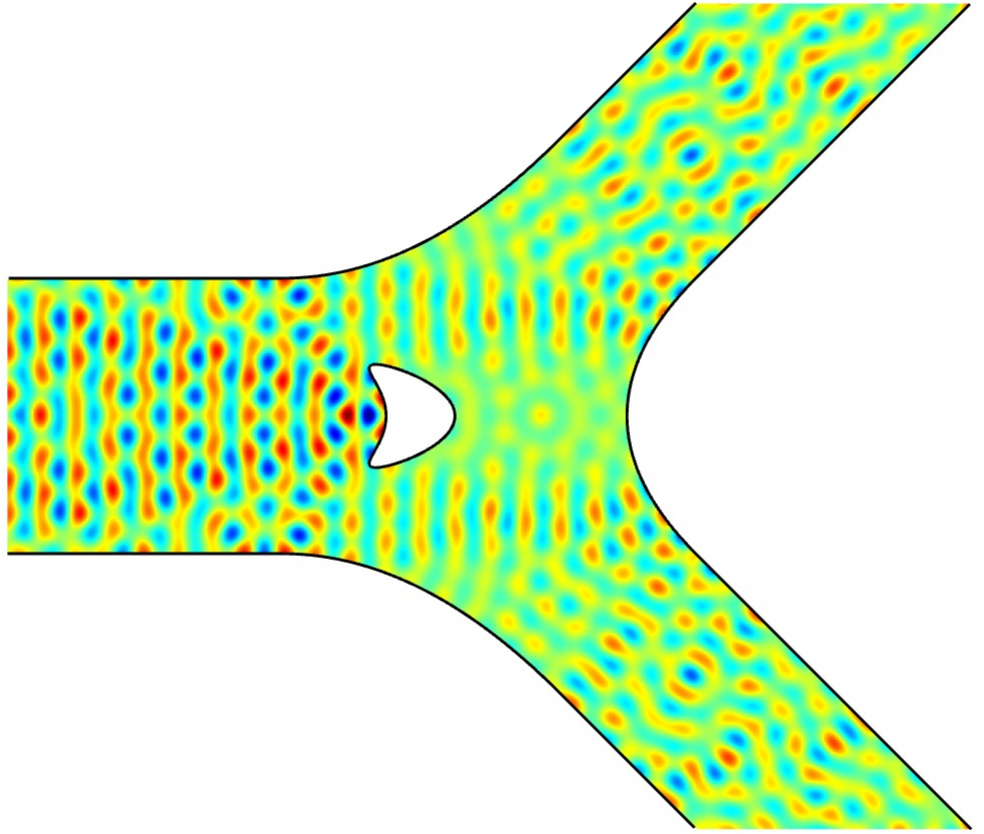}}\\
\subfloat[(e)]{\includegraphics[width=0.2\textwidth]{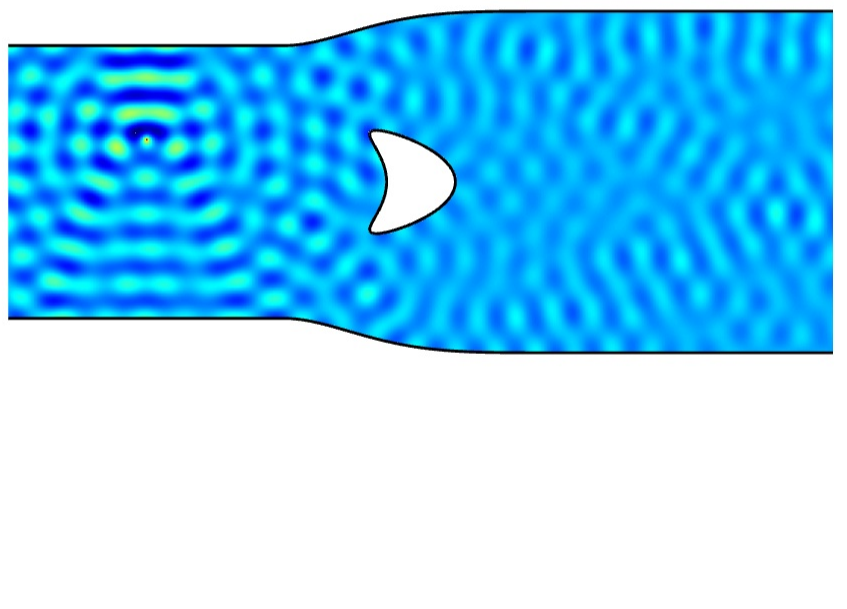}}&
\subfloat[(f)]{\includegraphics[width=0.2\textwidth]{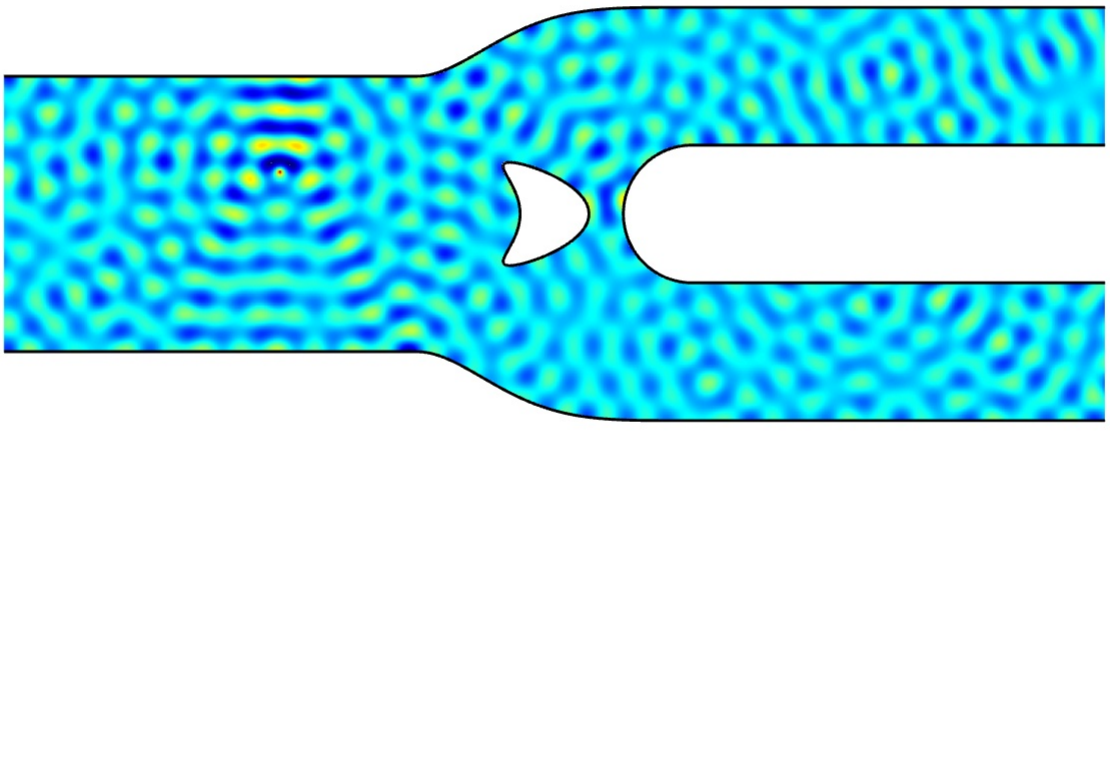}}&
\subfloat[(g)]{\includegraphics[width=0.2\textwidth]{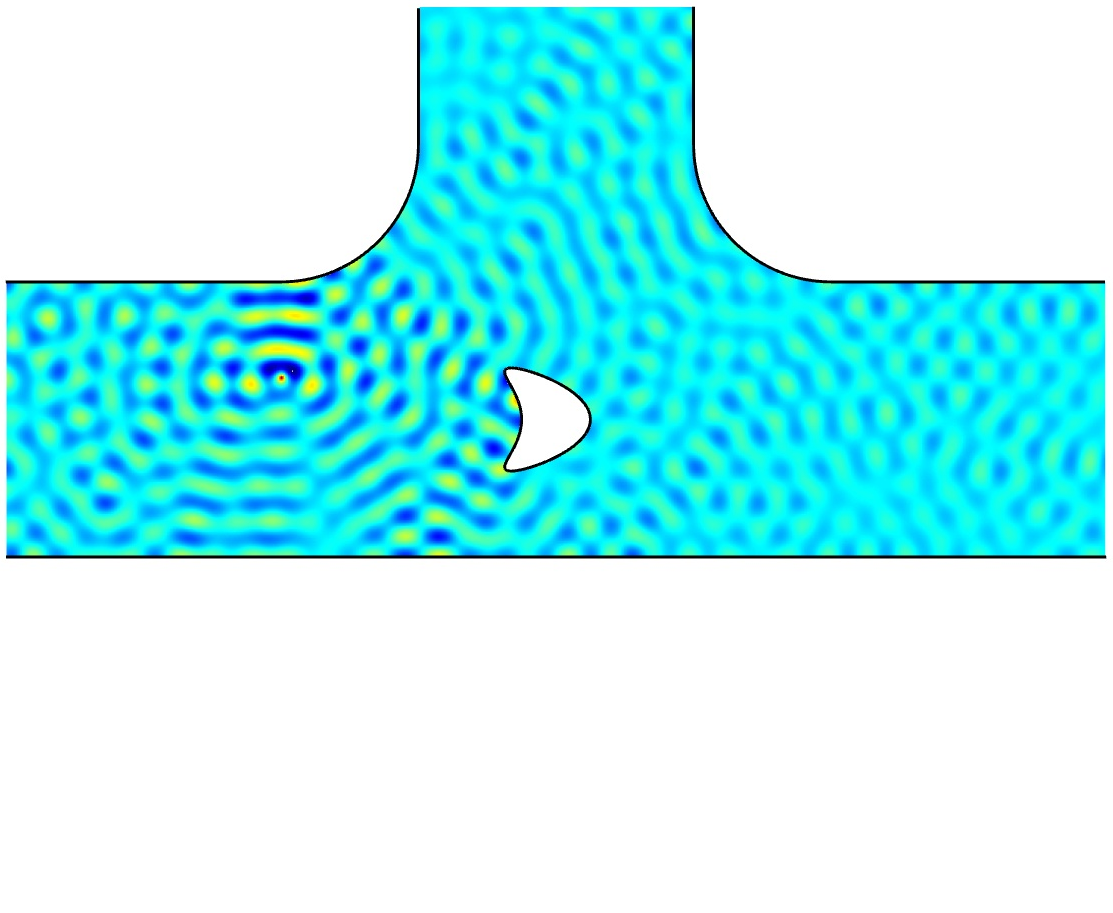}}&
\subfloat[(h)]{\includegraphics[width=0.2\textwidth]{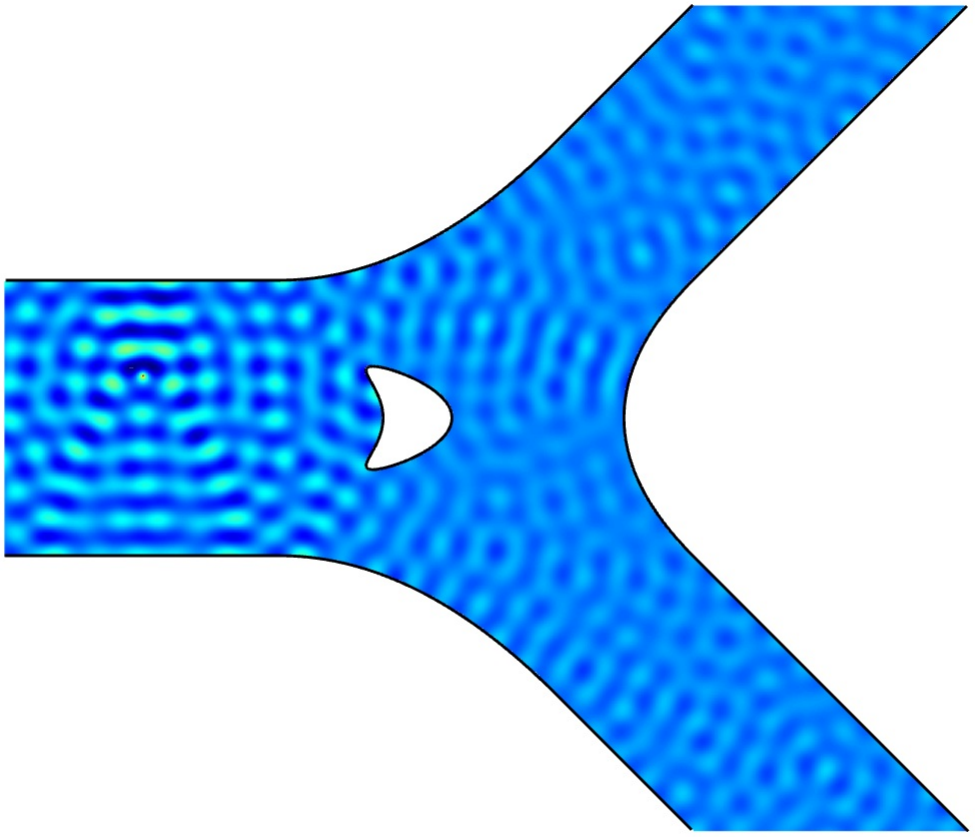}}
\end{tabular}
\caption{Example~2. Real parts of the total fields for mode (a-d) and point-source (e-h) incidence for the problems in nonuniform waveguides.}
\label{fig:FD-utot}
\end{figure}

\begin{figure}[!ht]
\centering
\setlength{\tabcolsep}{2pt}
\begin{tabular}{cccc}
\subfloat[two-port]{\includegraphics[width=0.235\textwidth]{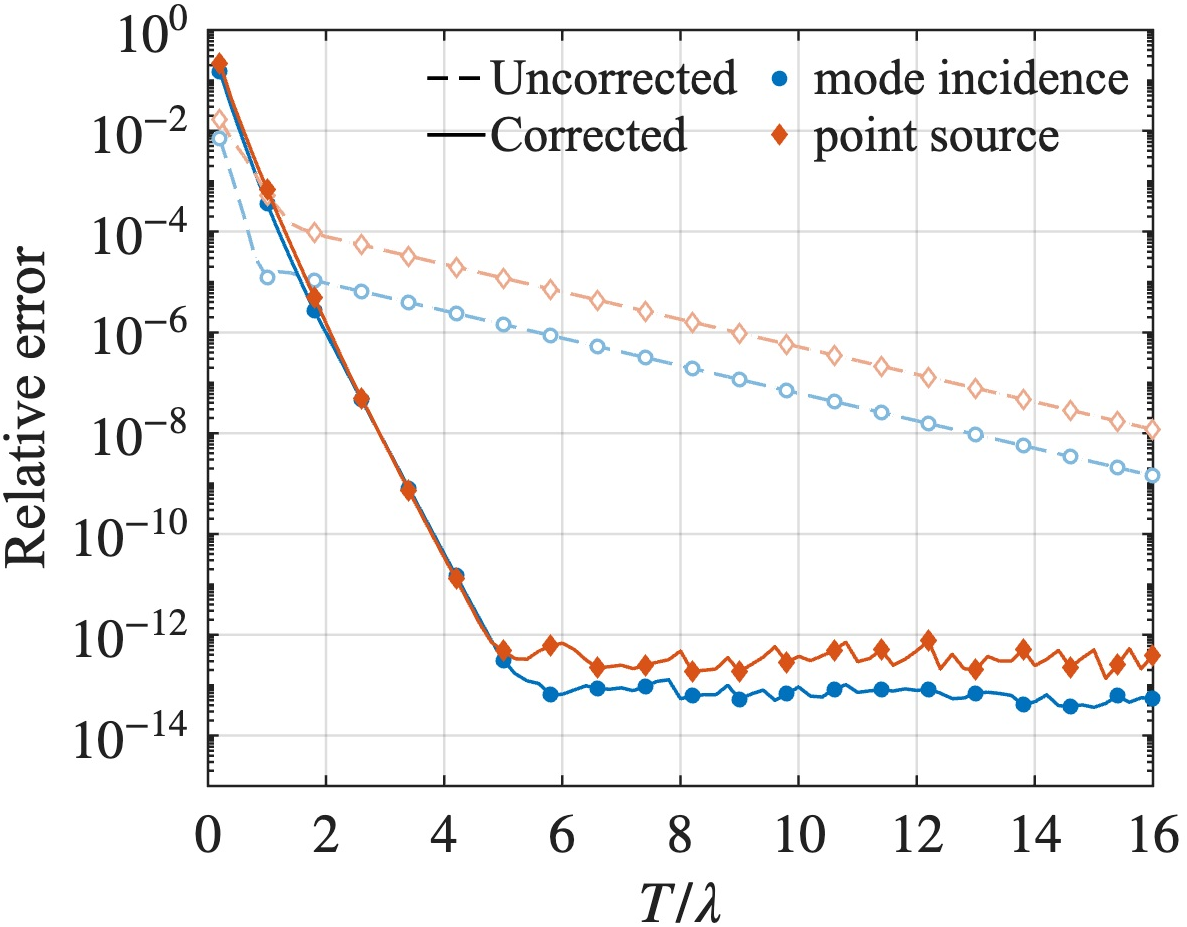}}&
\subfloat[parallel three-port]{\includegraphics[width=0.235\textwidth]{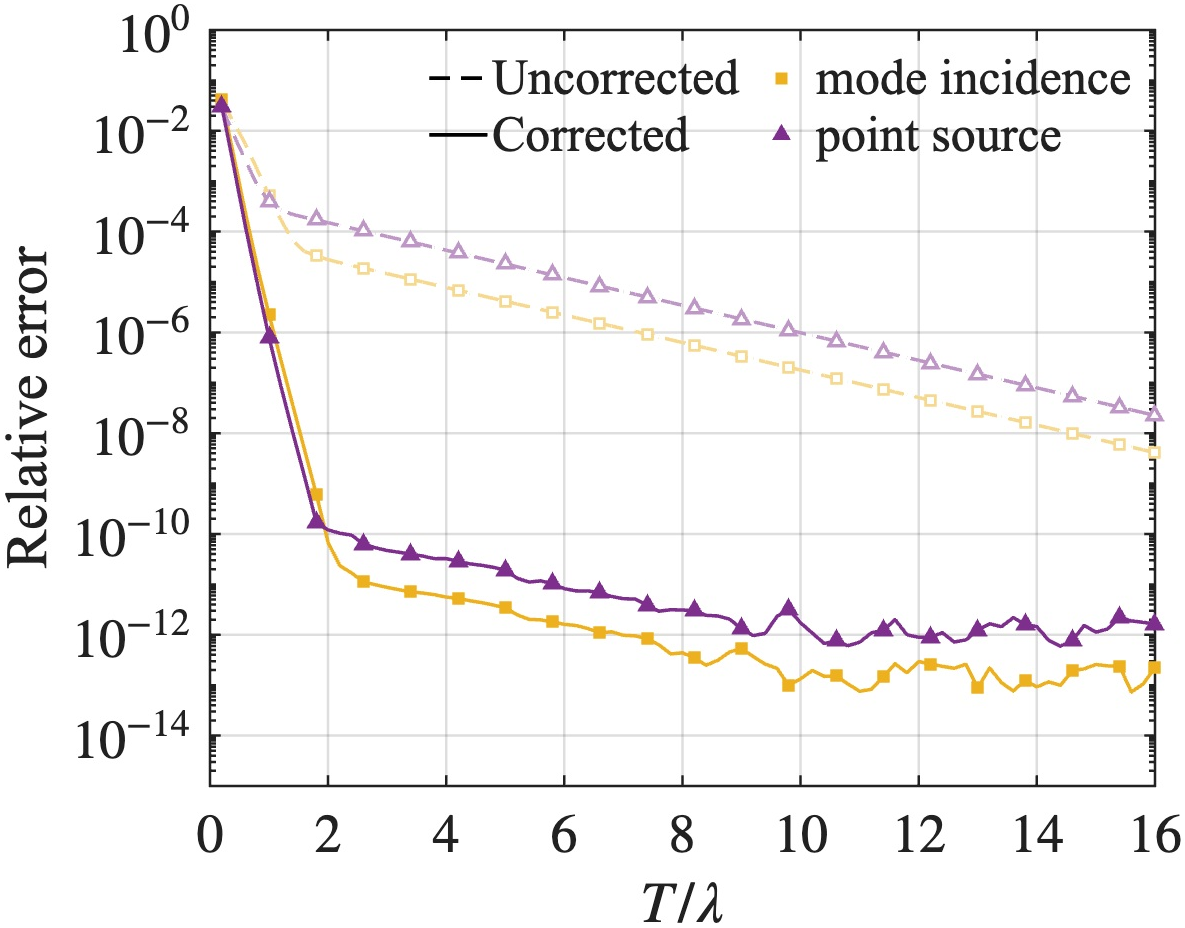}}&
\subfloat[T-junction]{\includegraphics[width=0.235\textwidth]{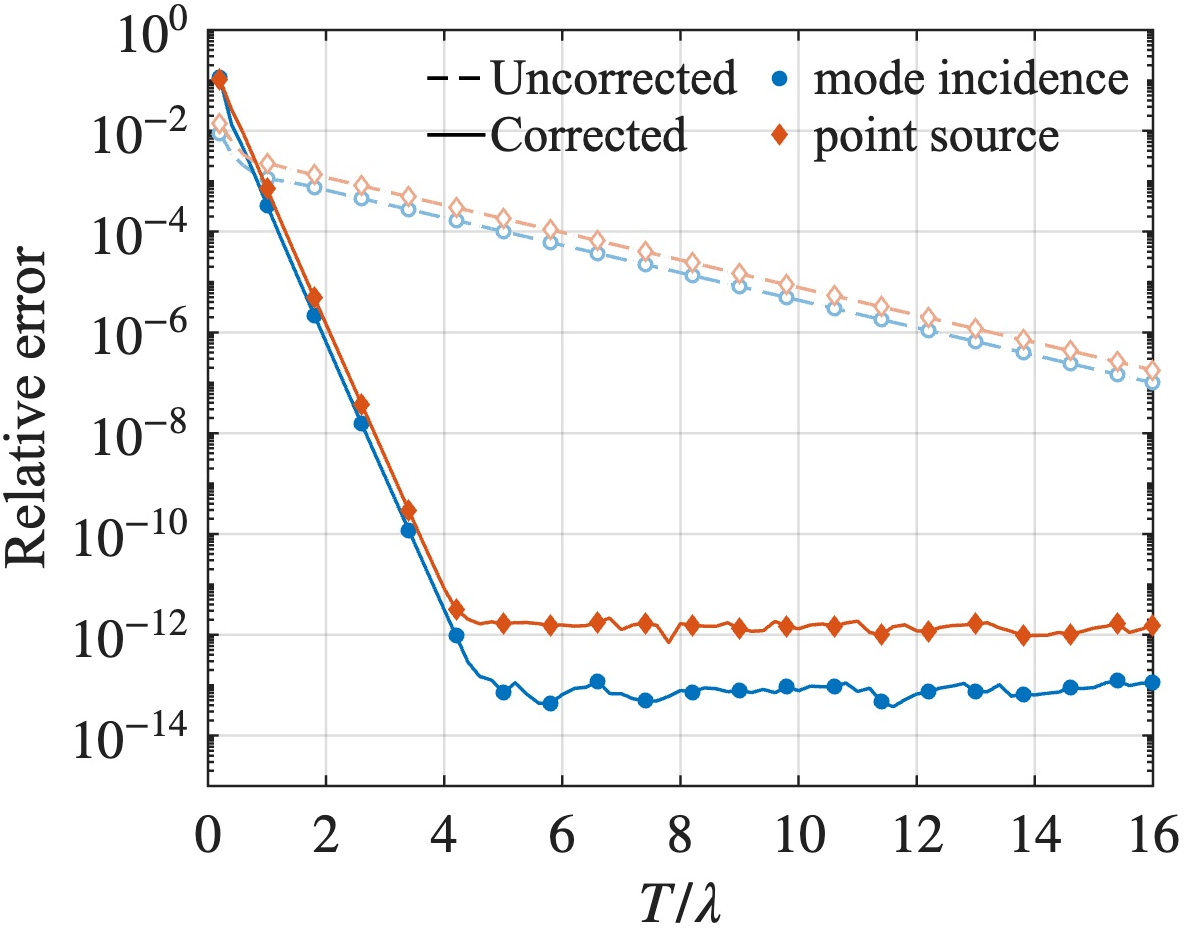}}&
\subfloat[Y-junction]{\includegraphics[width=0.235\textwidth]{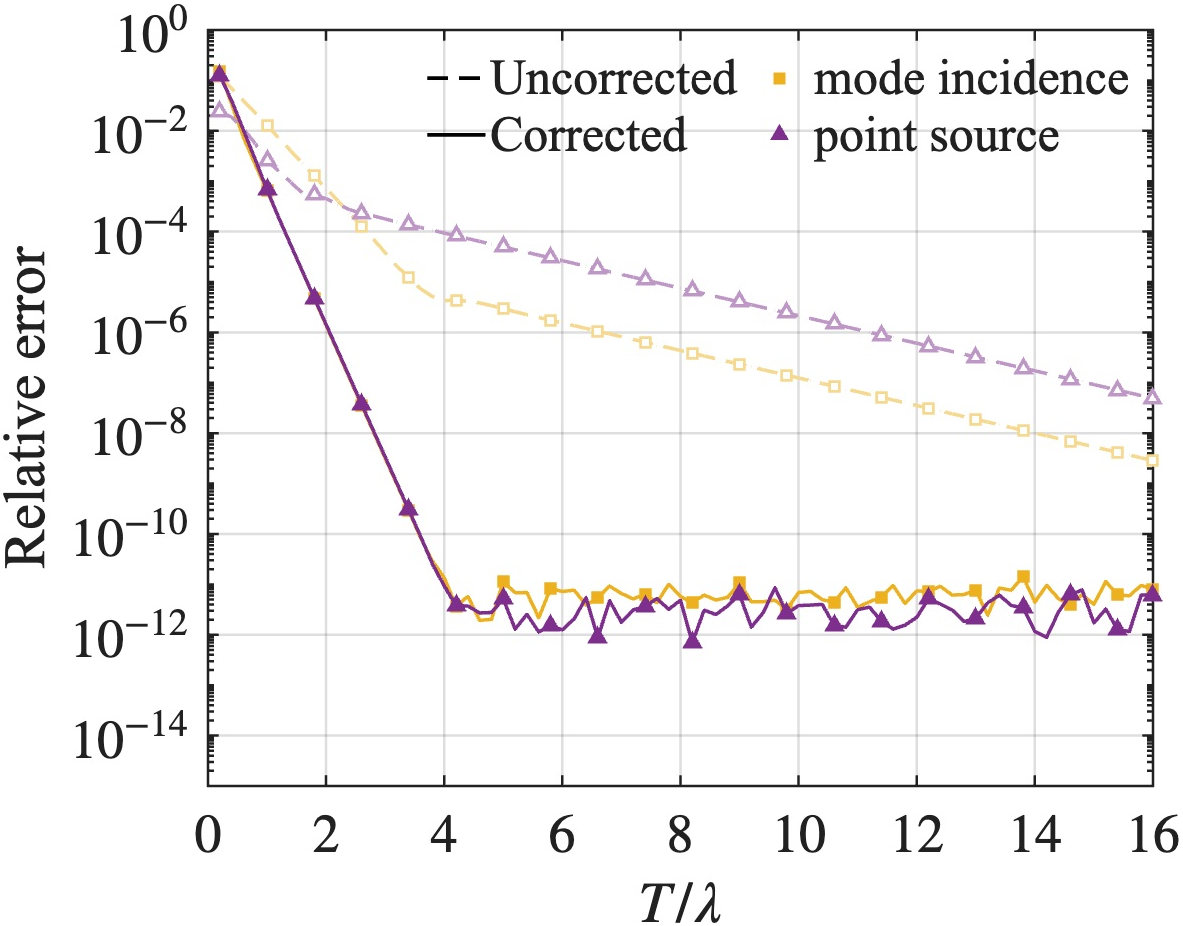}}
\end{tabular}
\caption{Example 2. Relative errors of the uncorrected and few-mode corrected RPML-BIE methods for problems in nonuniform waveguides.}
\label{fig:error-FD}
\end{figure}

{\bf Example~2} ({\it frequency-domain test for nonuniform waveguides}). We next consider the frequency-domain problems in nonuniform waveguides listed in Figure~\ref{fig:numerical-geometries} with a kite-shaped impedance obstacle inside, and choose $\kappa=7.99\pi$ and the same observation points as Example 1. The real parts of the total fields generating from the corrected RPML-BIE method are presented in Figure~\ref{fig:FD-utot}. The convergence of the relative errors \eqref{eq:FD-relative-error} with respect to $T$ shown in Figure~\ref{fig:error-FD} indicates that the corrected RPML-BIE method significantly improves the exponential convergence and can achieve a prescribed high accuracy with a substantially thinner layers near the cut-off frequencies.

\FloatBarrier

\begin{figure}[!ht]
\centering
\setlength{\tabcolsep}{1pt}
\begin{tabular}{cccc}
\subfloat[Sound-hard case]{\includegraphics[width=0.23\textwidth]{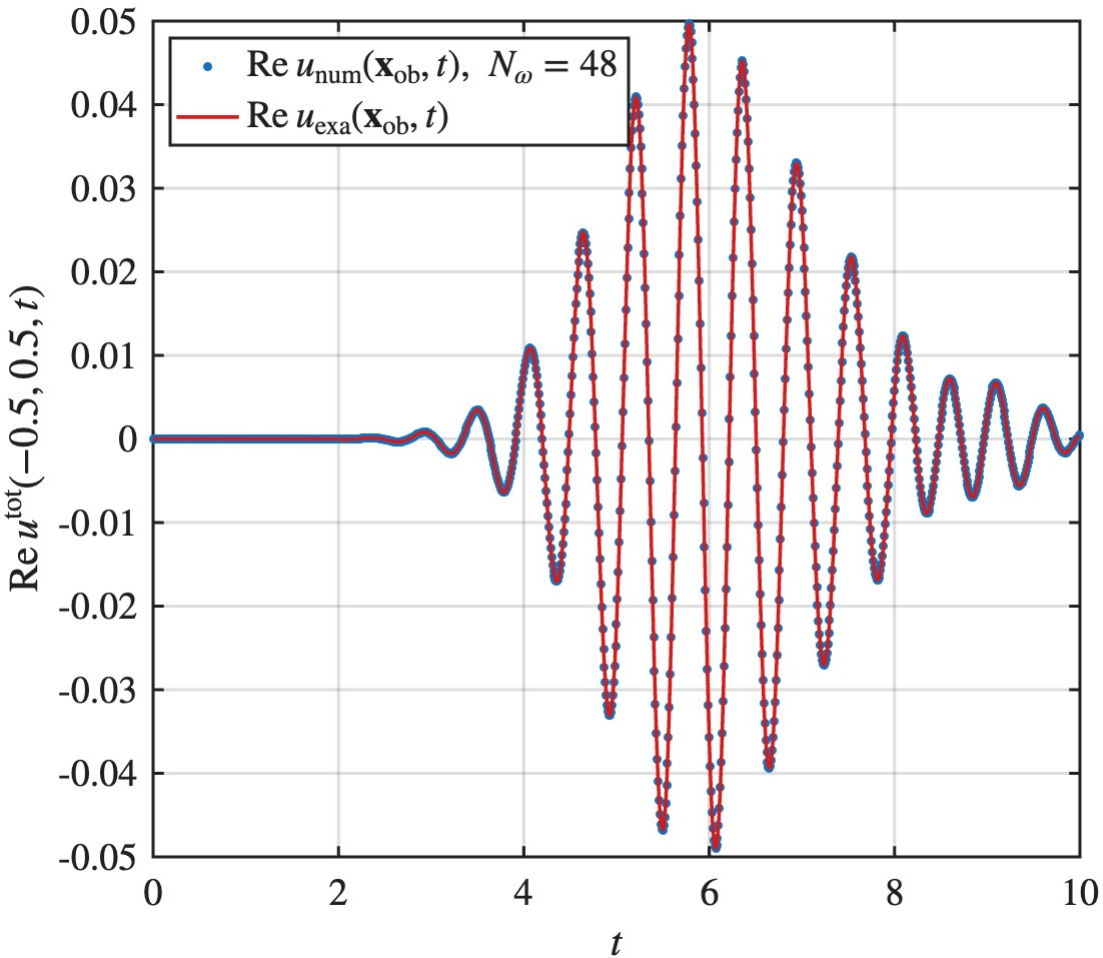}}&
\subfloat[Sound-hard case]{\includegraphics[width=0.23\textwidth]{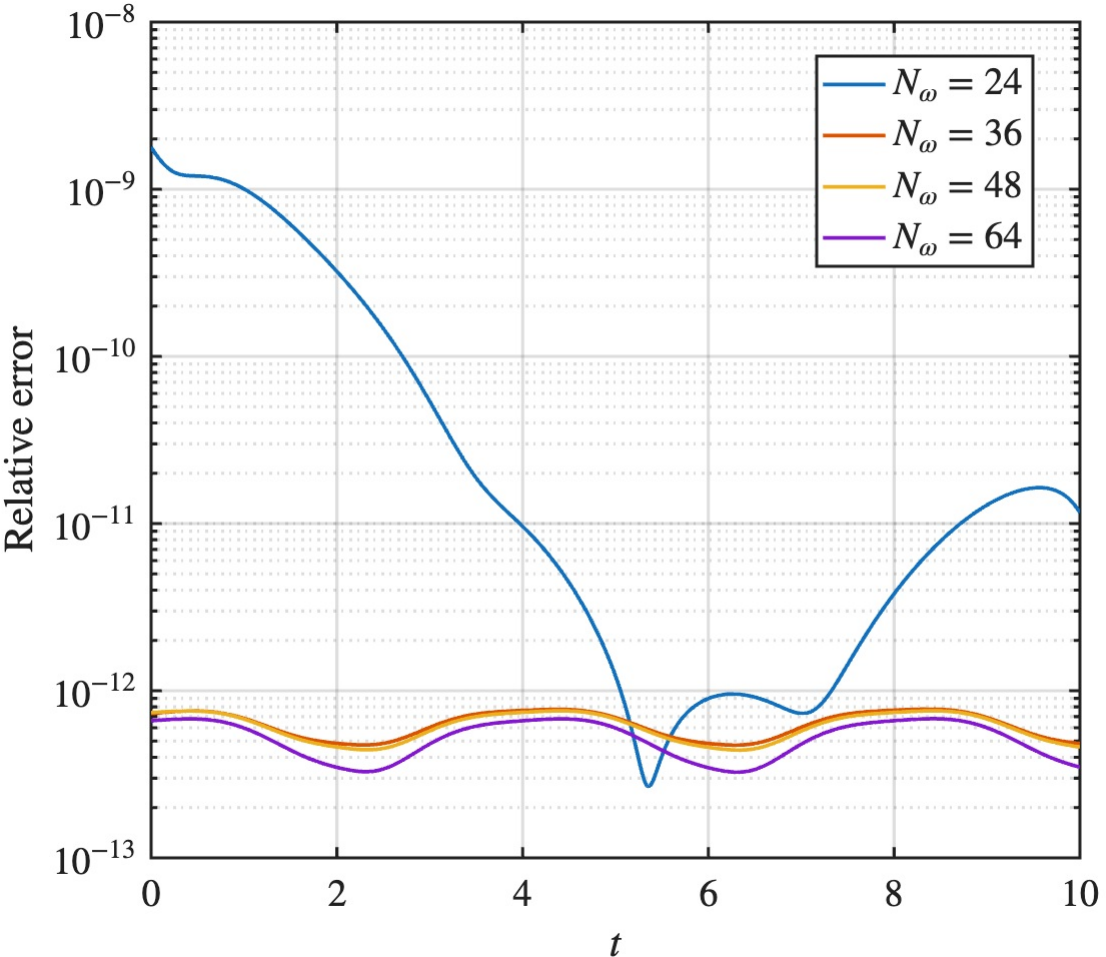}}&
\subfloat[Sound-soft case]{\includegraphics[width=0.23\textwidth]{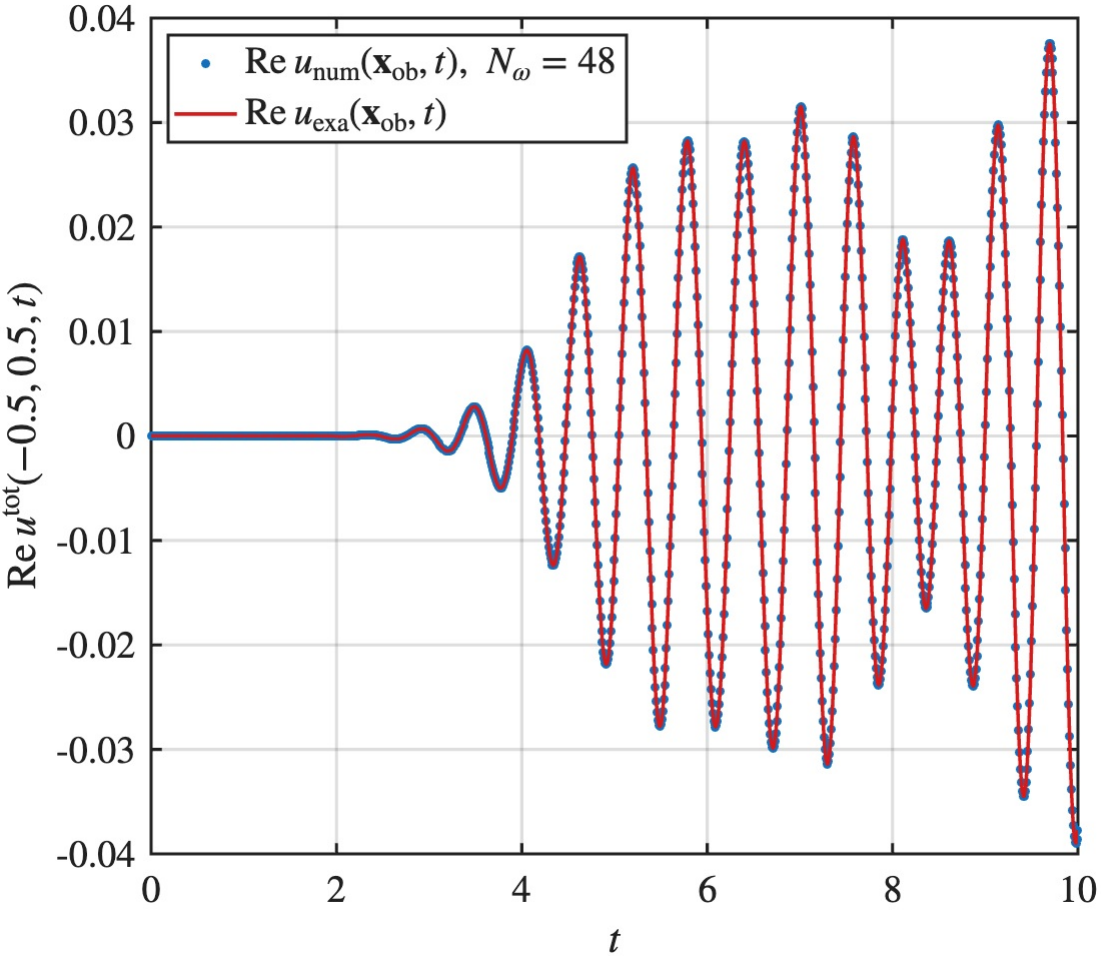}}&
\subfloat[Sound-soft case]{\includegraphics[width=0.23\textwidth]{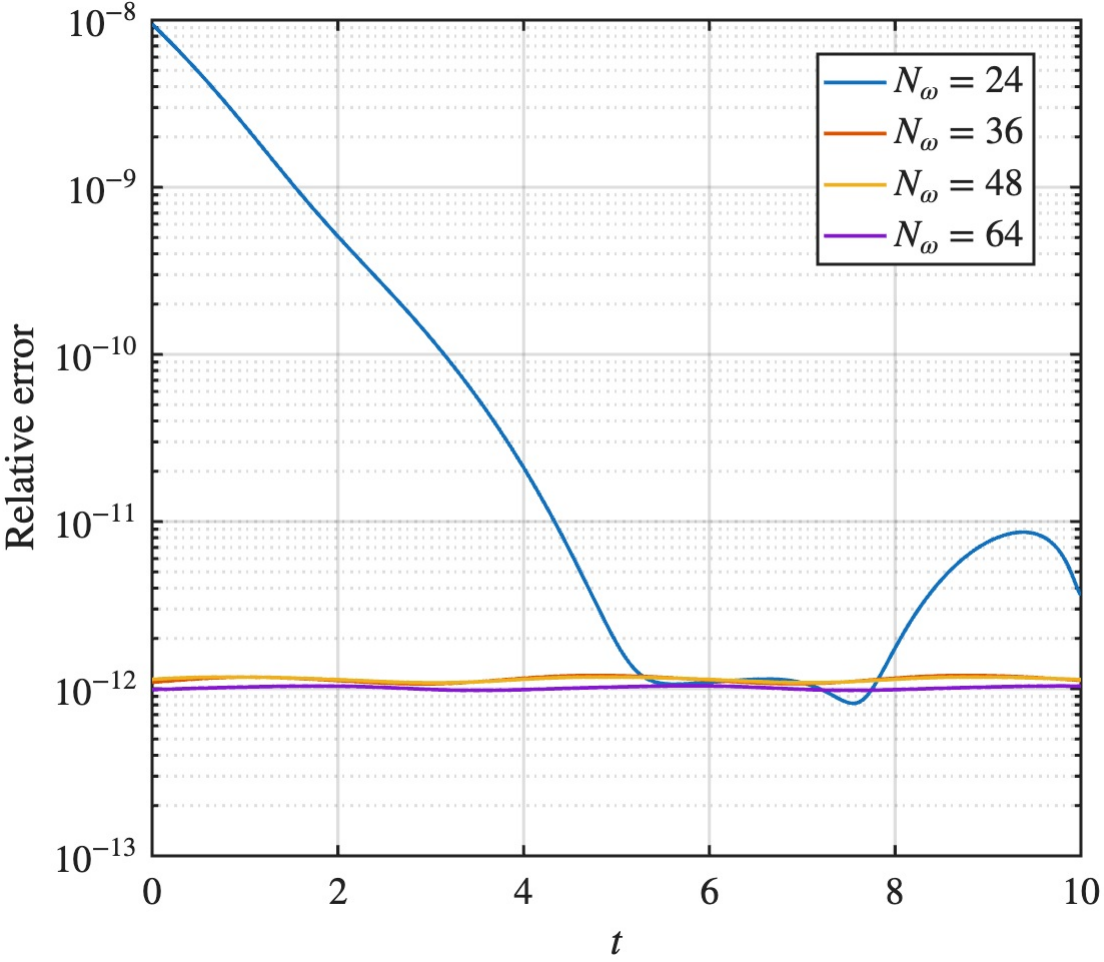}}
\end{tabular}
\caption{Example 3. Time traces of the exact and numerical solutions at the observation point $\bm x_{\rm ob}=(-0.5,0.5)$ (a,c) and numerical errors for different values of $N_\omega$ (b,d).}
\label{fig:Ex3-1-validation}
\end{figure}

{\bf Example 3} ({\it time-domain test for a planar waveguide}). We consider the propagation of an incident Gauss pulse of point source placed at $\bm z=(-2,1.3)$ in an obstacle-free planar waveguide with width $h=2$. Fix the RPML thickness $T=4\lambda_{\rm ref}=2$, where $\lambda_{\rm ref}=2\pi c/\omega_{\rm ref}$ with $\omega_{\rm ref}=4\pi$. The exact time-domain total field for the sound-hard planar waveguide is given by
\be\label{eq:straight-td-exact}
u_{\rm exa}(\bm x,t)
=\frac{A_0c}{2}\sum_{n\ge0}\psi_n(x_2)\psi_n(z_2)
\int_{-\infty}^{t-|x_1-z_1|/c}a(s)
J_0\!\left(\alpha_n\sqrt{c^2(t-s)^2-|x_1-z_1|^2}\right)\,\mathrm ds,
\en
where $J_0$ is the Bessel function of the first kind of order zero and $a$ is given by \eqref{eq:td_spectrum}. The exact solution for the sound-soft planar waveguide takes an analogous form. Figure~\ref{fig:Ex3-1-validation}(a,c) display the time-traces of the exact and numerical solutions at the observation point $\bm x_{\rm ob}=(-0.5,0.5)$ when $N_\omega=48$. The numerical errors $\varepsilon_{\rm TD}(\bm x_{\rm ob},t)$ over the time interval $0\le t\le 10$ shown in Figure~\ref{fig:Ex3-1-validation}(b,d) indicates the high accuracy of the proposed frequency-time hybrid RPML-BIE solver.

\FloatBarrier

\begin{figure}[!htb]
\centering
\setlength{\tabcolsep}{1pt}
\begin{tabular}{ccccc}
\includegraphics[width=0.19\textwidth]{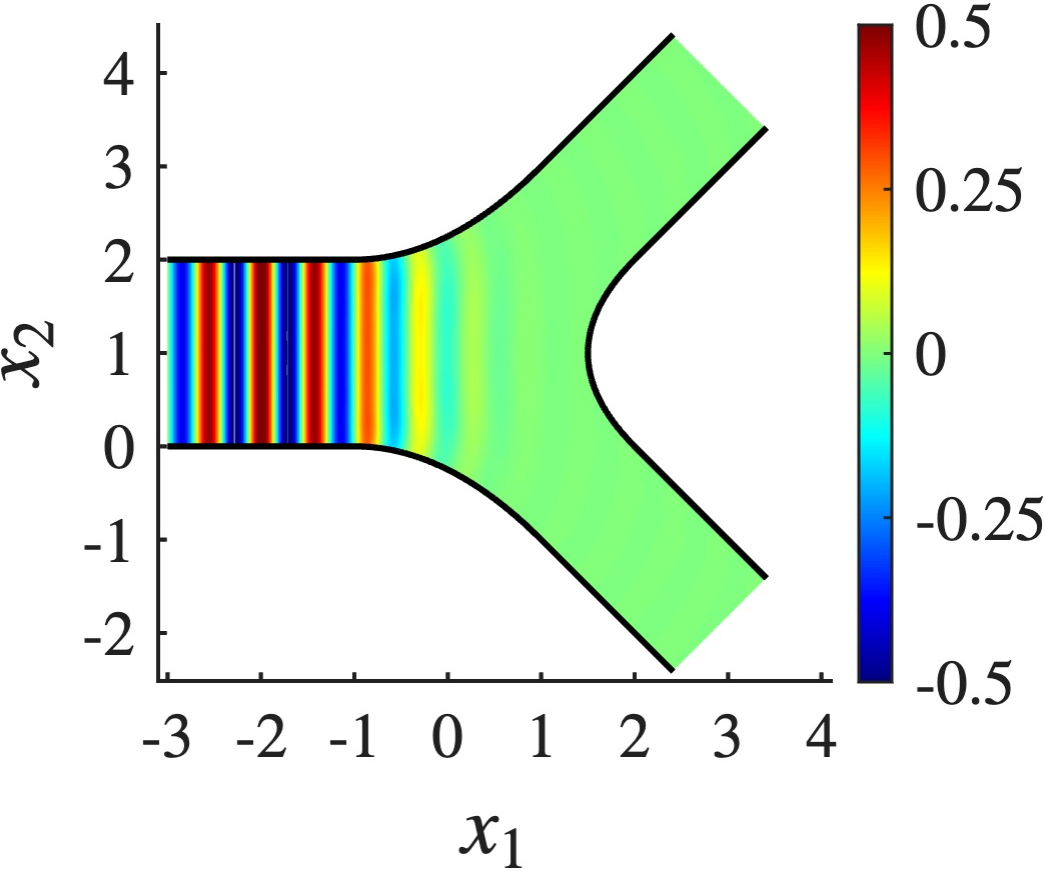}&
\includegraphics[width=0.19\textwidth]{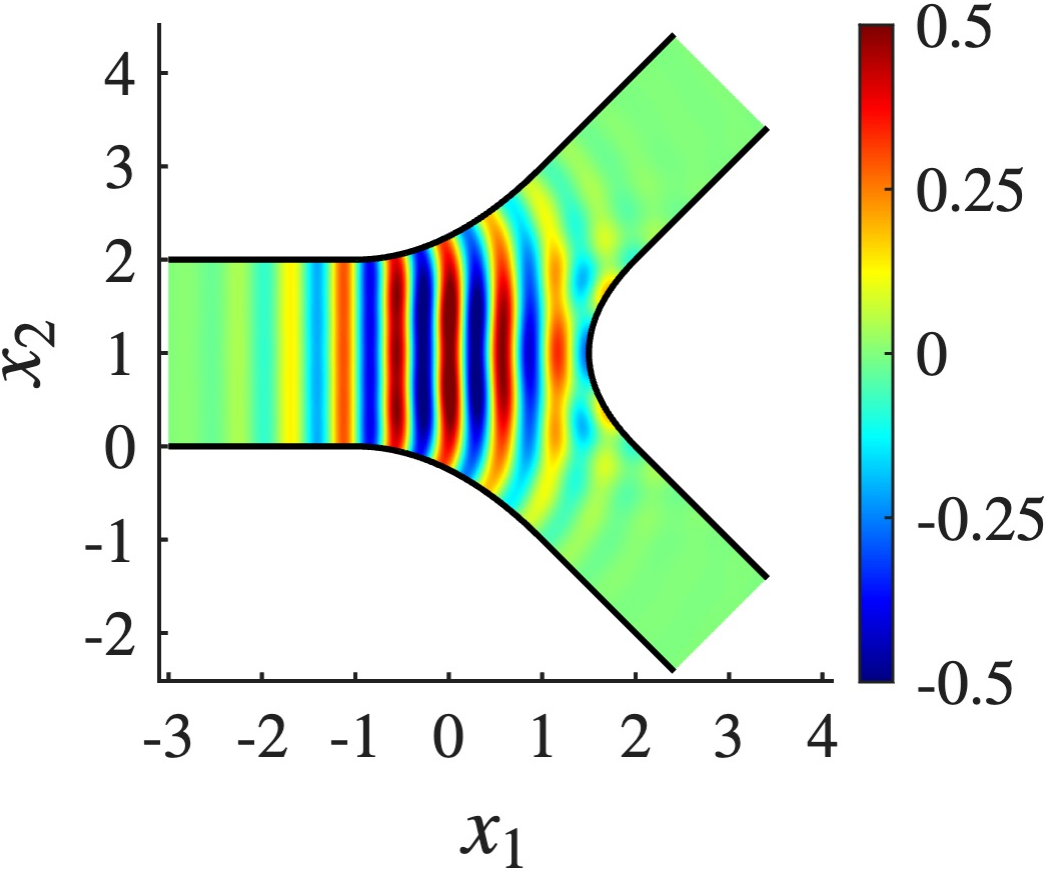}&
\includegraphics[width=0.19\textwidth]{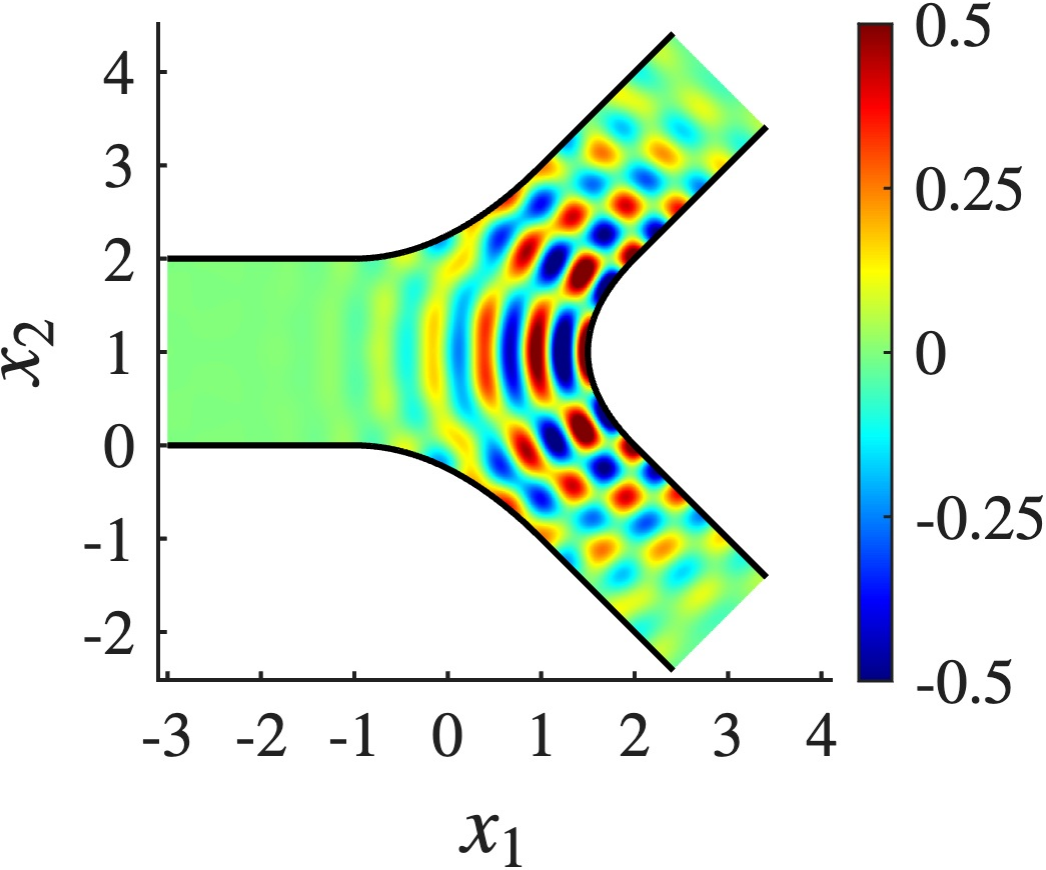}&
\includegraphics[width=0.19\textwidth]{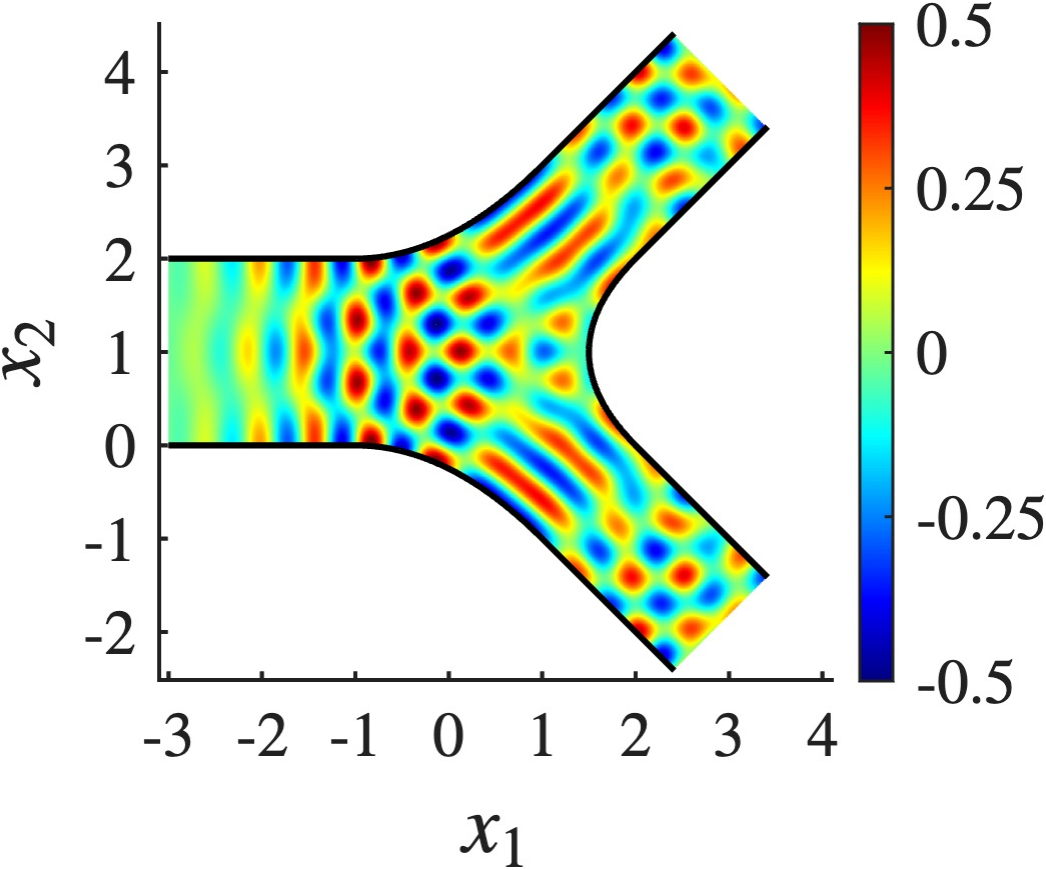}&
\includegraphics[width=0.19\textwidth]{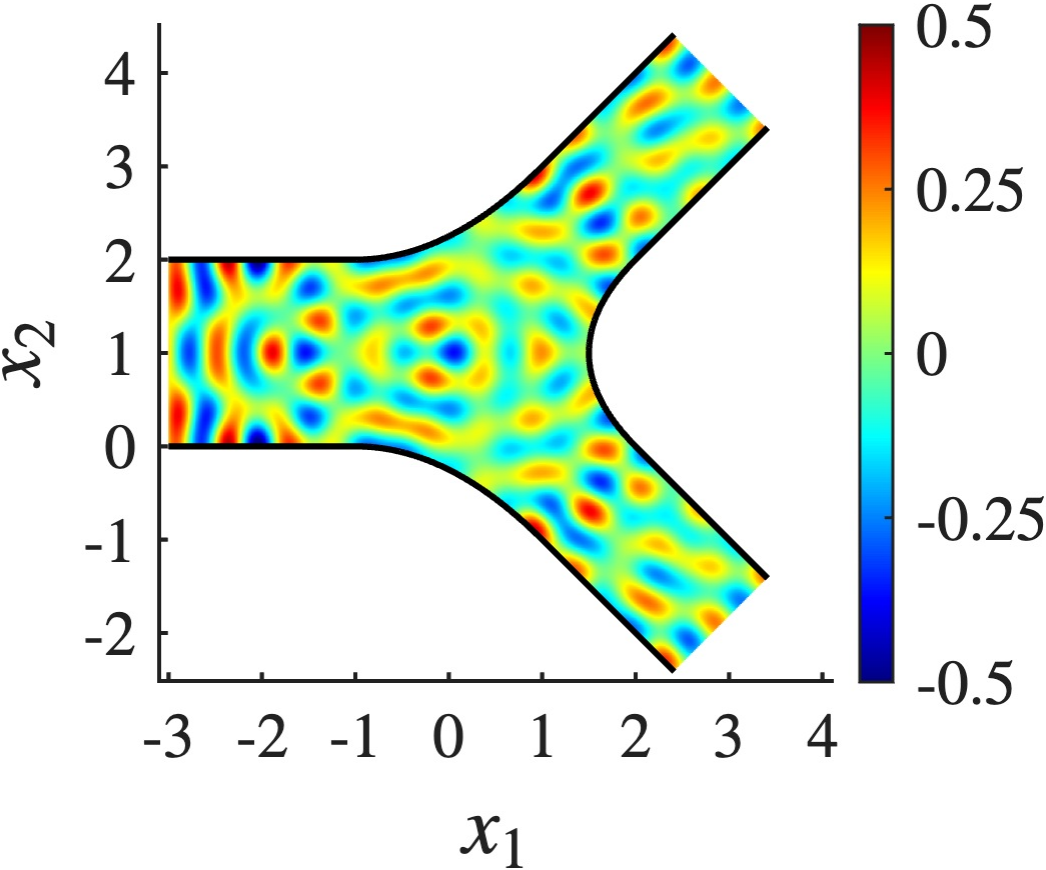}\\[-3mm]
\includegraphics[width=0.19\textwidth]{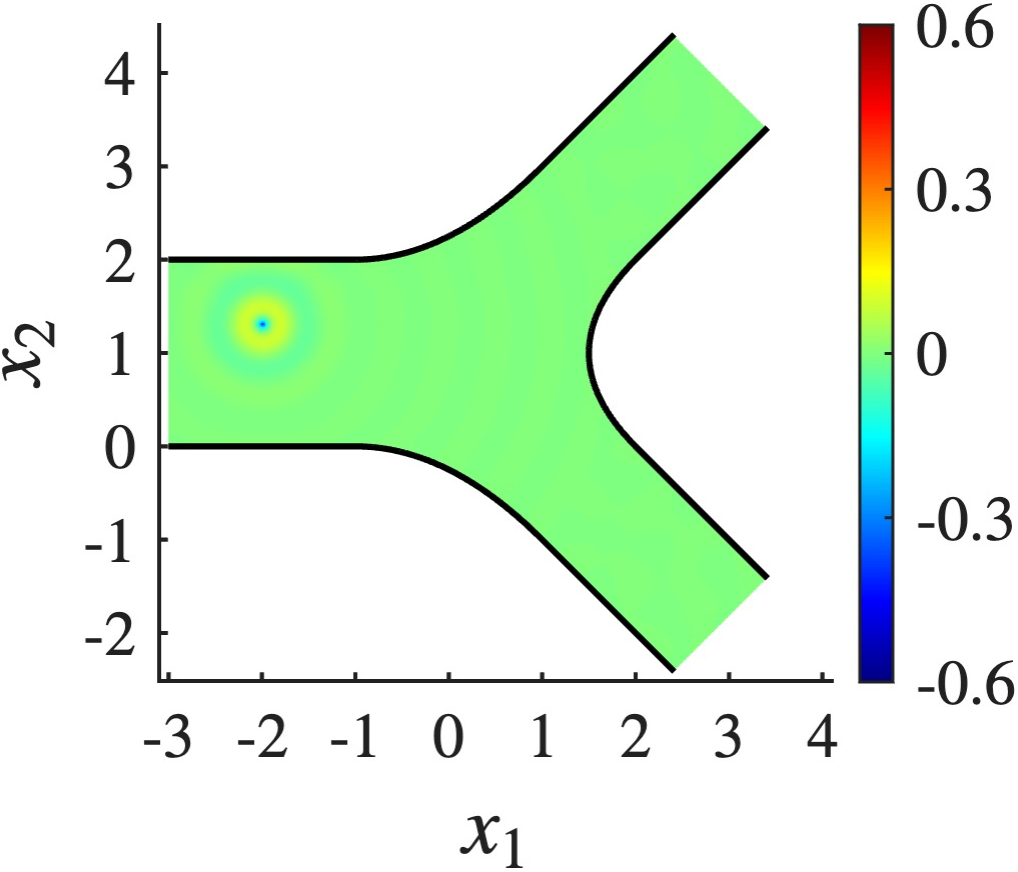}&
\includegraphics[width=0.19\textwidth]{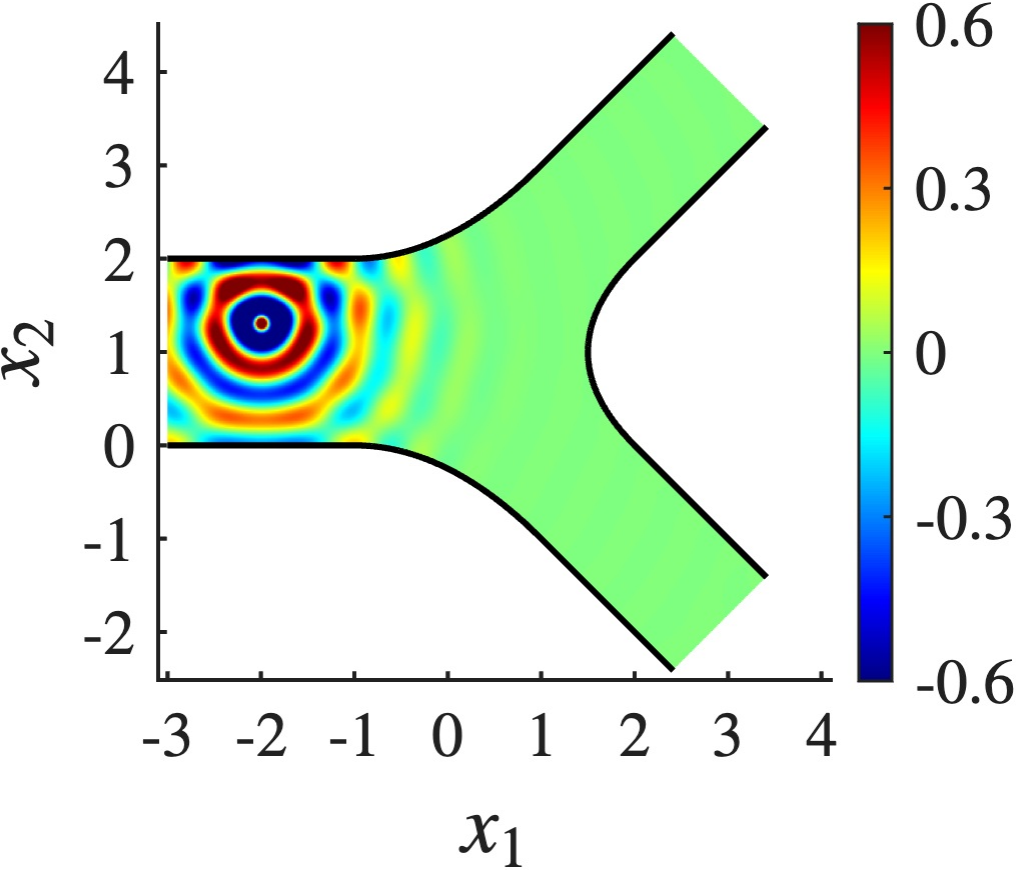}&
\includegraphics[width=0.19\textwidth]{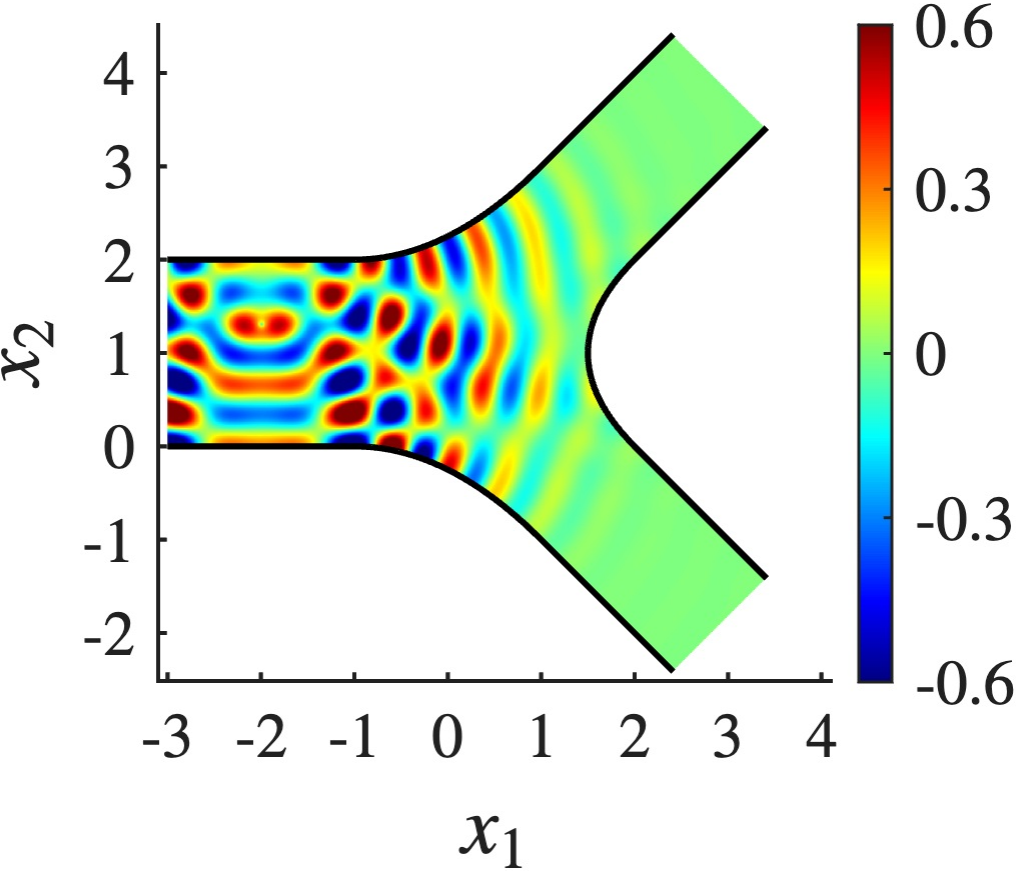}&
\includegraphics[width=0.19\textwidth]{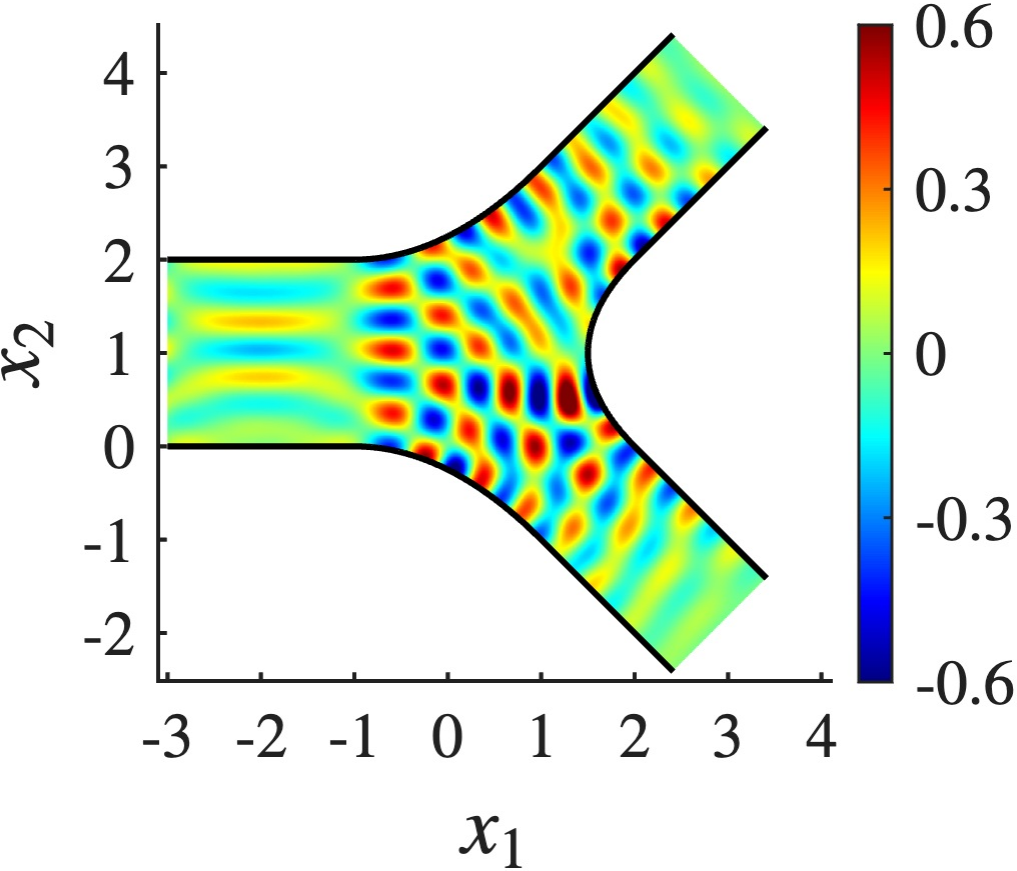}&
\includegraphics[width=0.19\textwidth]{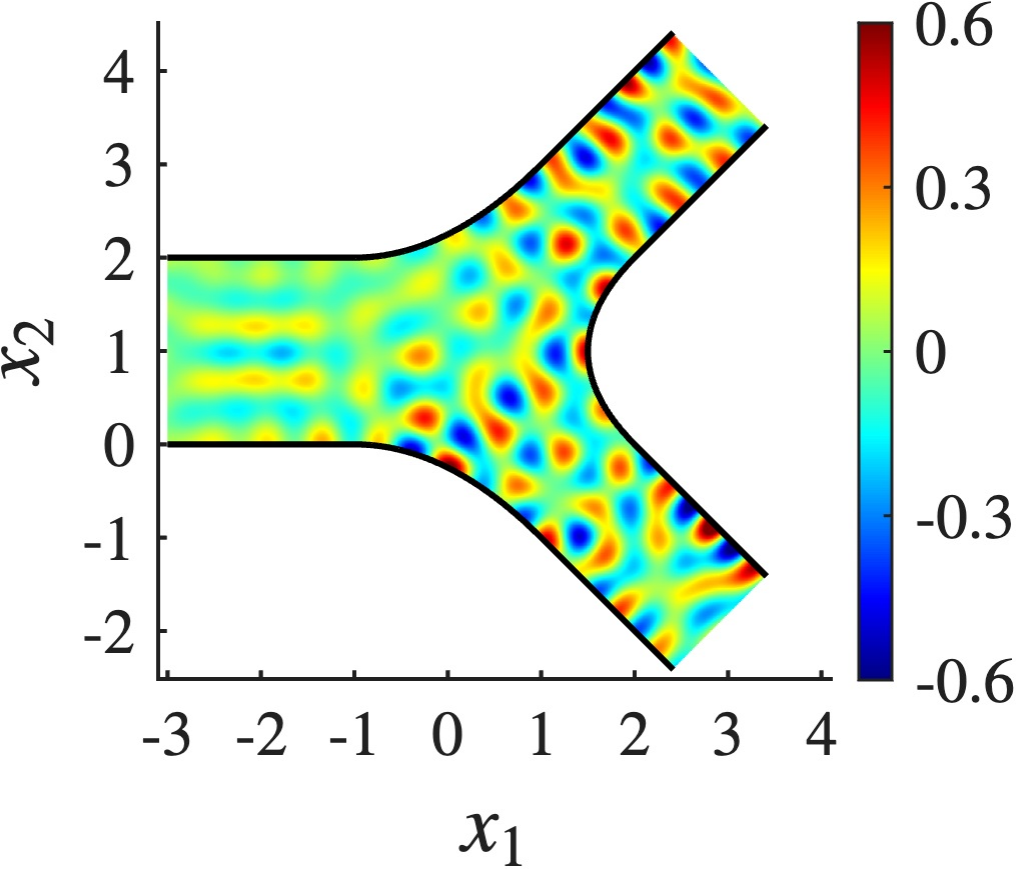}
\end{tabular}
\caption{Example 4. Real parts of the total fields at $t=2,4,6,8,10$ (left to right) for the problem of mode (top) and point-source (bottom) incidences in a Y-junction waveguide.}
\label{fig:Ex3-3-utot}
\end{figure}

\begin{figure}[!htbp]
\centering
\setlength{\tabcolsep}{2pt}
\def\FigureTwelveReady{0}
\IfFileExists{Ex3-3_short_time_mode_real_imag.pdf}{\IfFileExists{Ex3-3_short_time_mode_TD_abs_error.pdf}{\IfFileExists{Ex3-3_short_time_source_real_imag.pdf}{\IfFileExists{Ex3-3_short_time_source_TD_abs_error.pdf}{\def\FigureTwelveReady{1}}{}}{}}{}}{}
\ifnum\FigureTwelveReady=1
\begin{tabular}{cccc}
\subfloat[Mode incidence]{\includegraphics[width=0.24\textwidth]{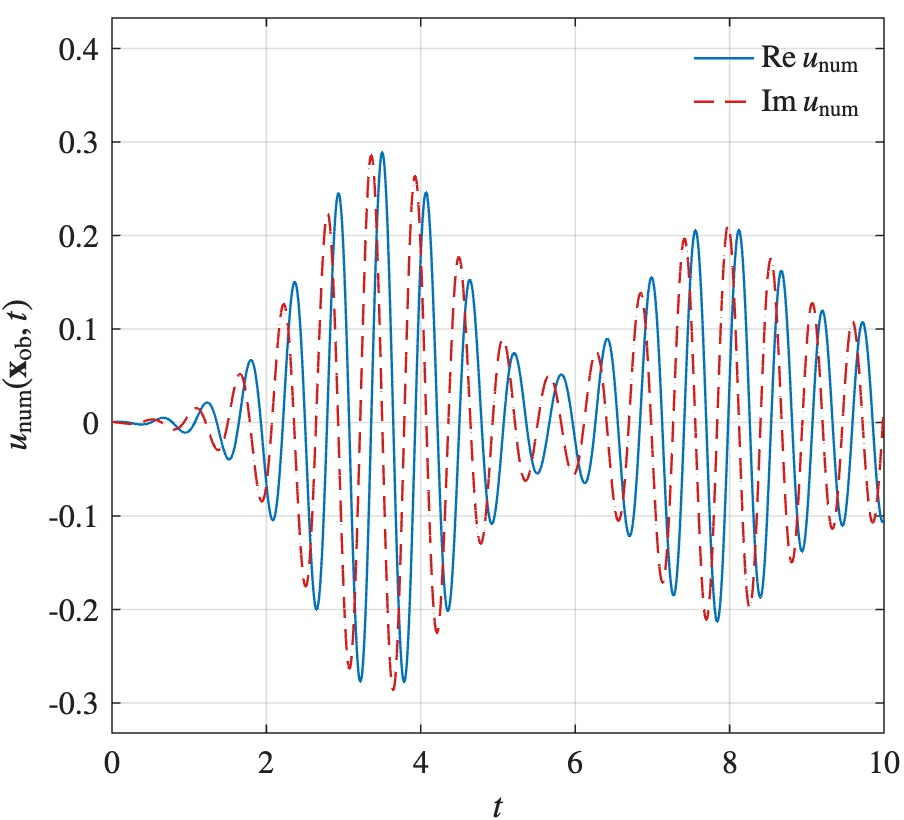}}&
\subfloat[Mode incidence]{\includegraphics[width=0.24\textwidth]{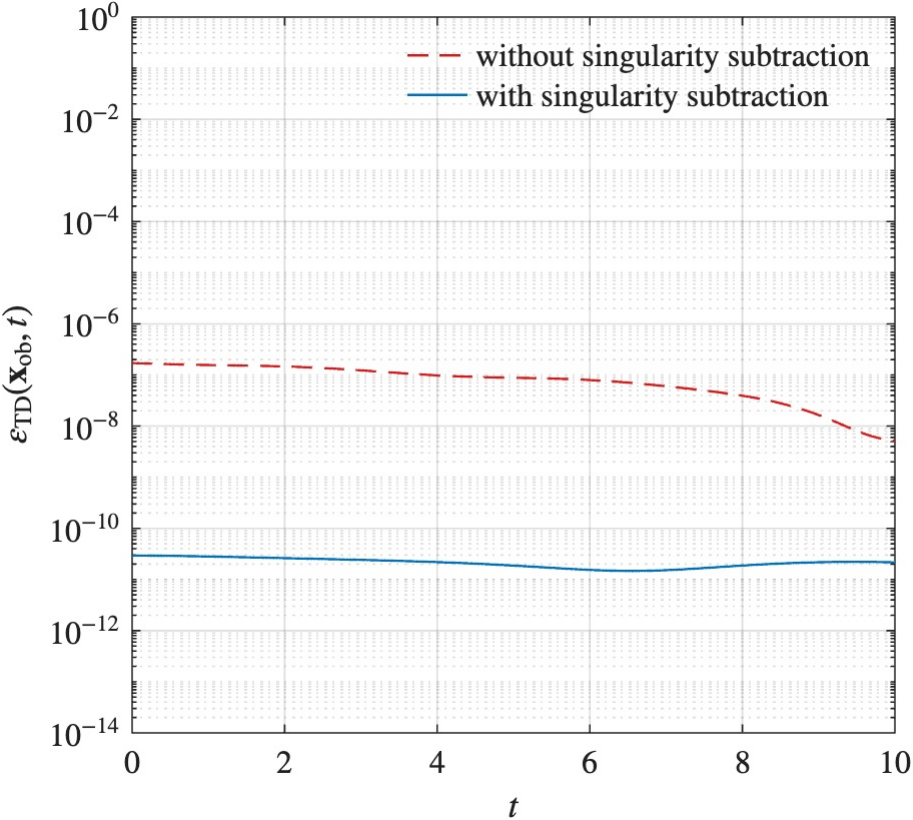}}&
\subfloat[Point source]{\includegraphics[width=0.24\textwidth]{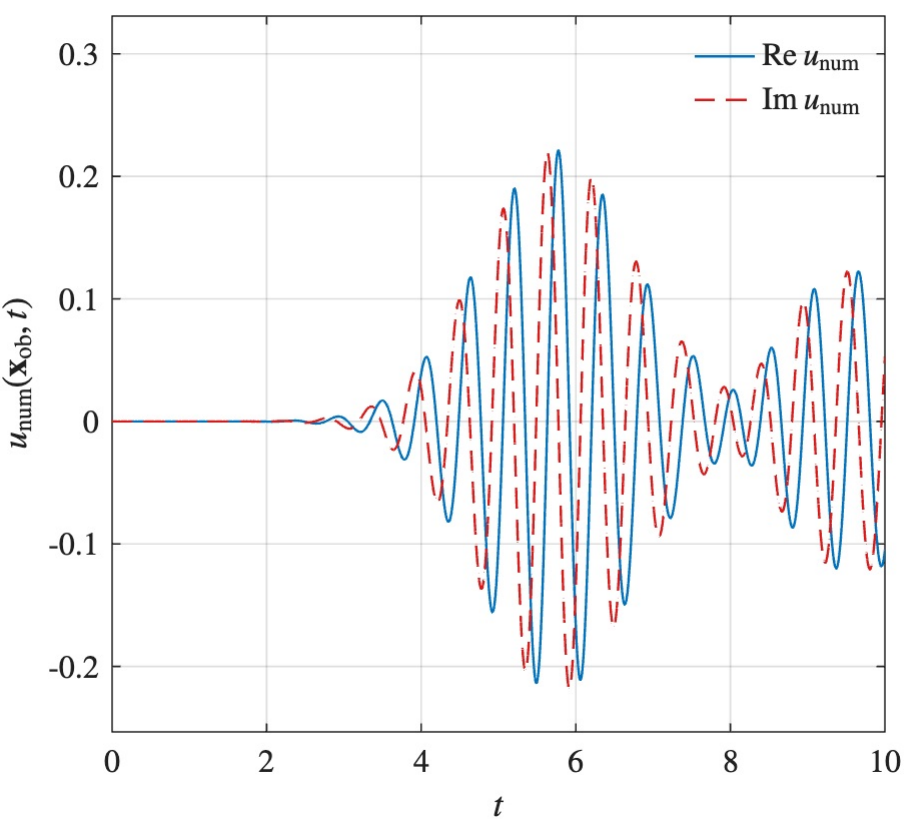}}&
\subfloat[Point source]{\includegraphics[width=0.24\textwidth]{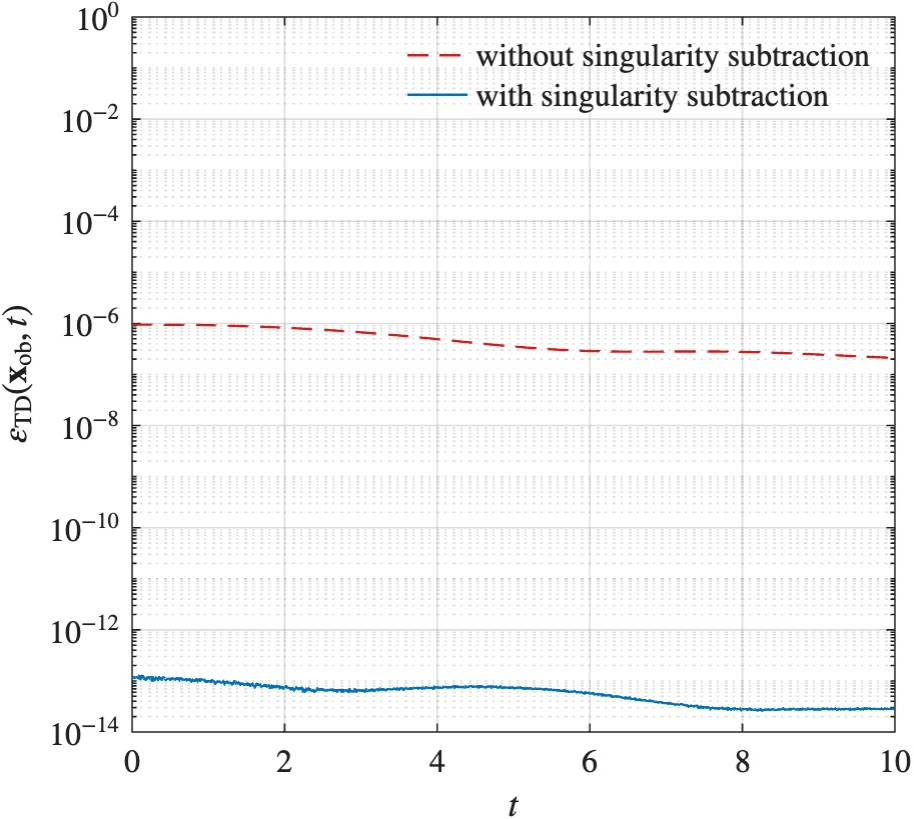}}
\end{tabular}
\caption{Example 4. Real and imaginary parts of $u_{\rm num}(\bm x_{\rm ob},t)$ (a,c) and numerical errors for the methods with or without singularity subtraction (b,d), for $N_\omega=48$.}
\fi
\label{fig:Ex3-3-single-pulse-AAA}
\end{figure}

\begin{figure}[!ht]
\centering
\setlength{\tabcolsep}{1pt}
\begin{tabular}{ccccc}
\includegraphics[width=0.19\textwidth]{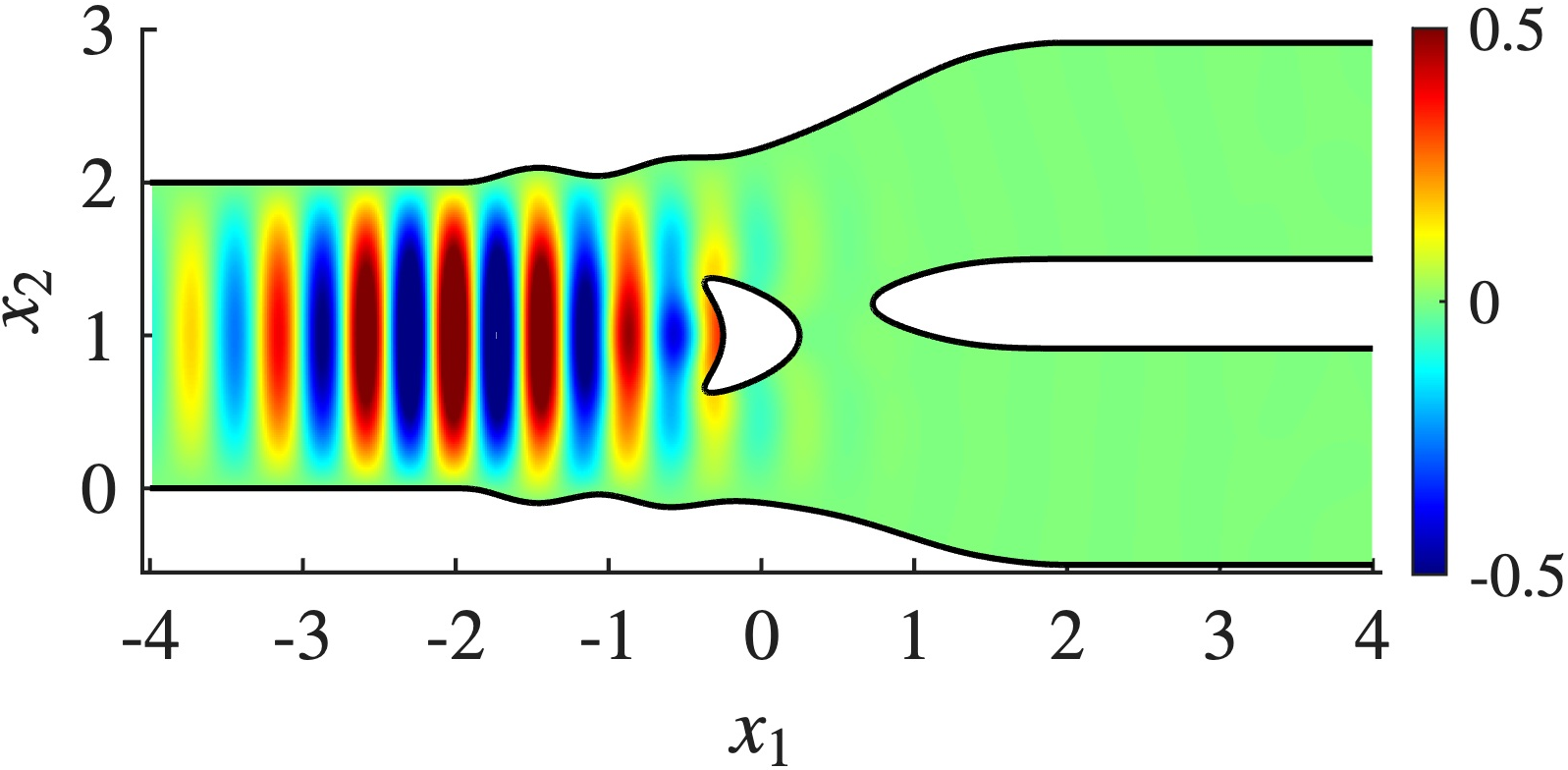}&
\includegraphics[width=0.19\textwidth]{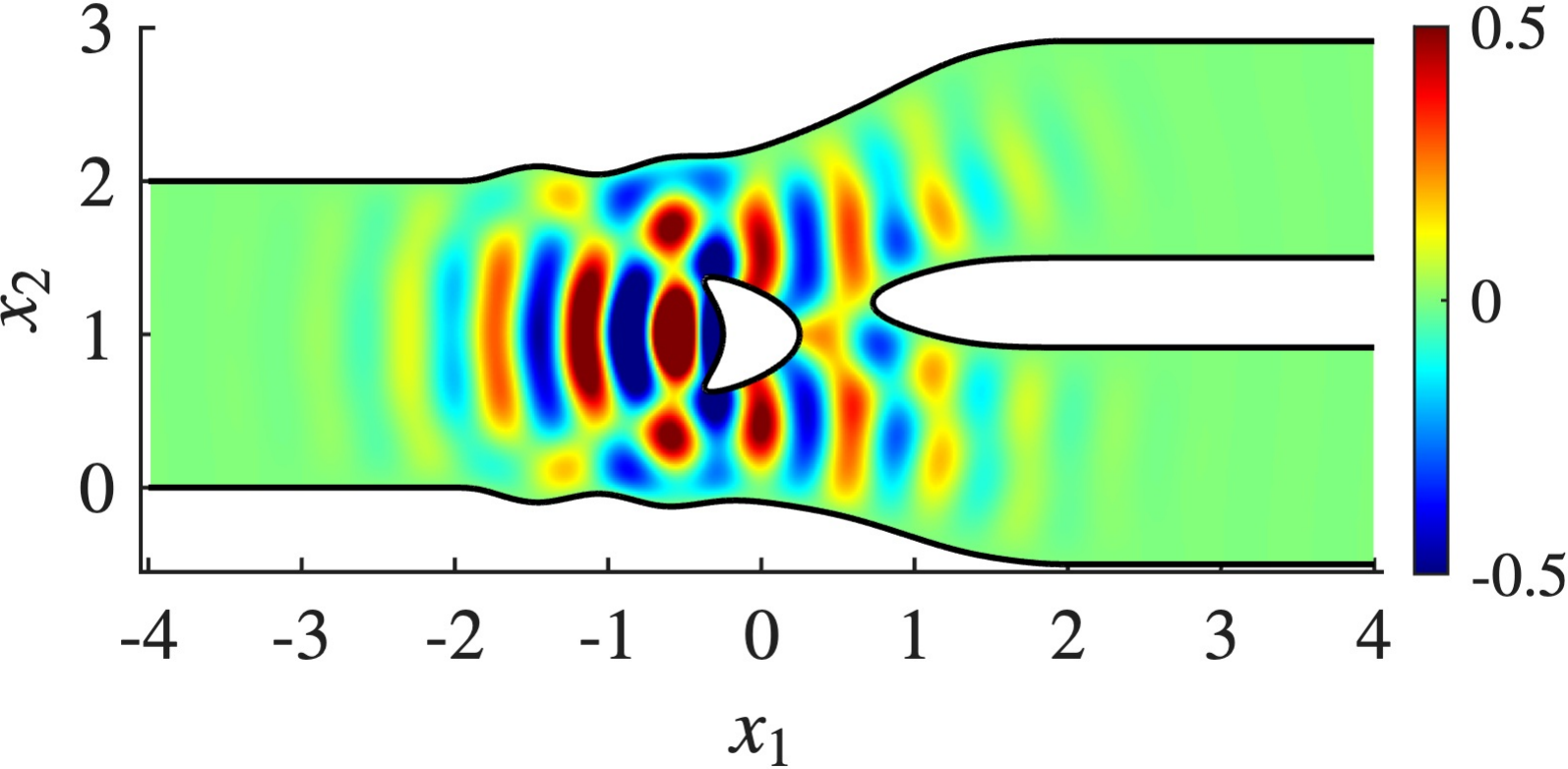}&
\includegraphics[width=0.19\textwidth]{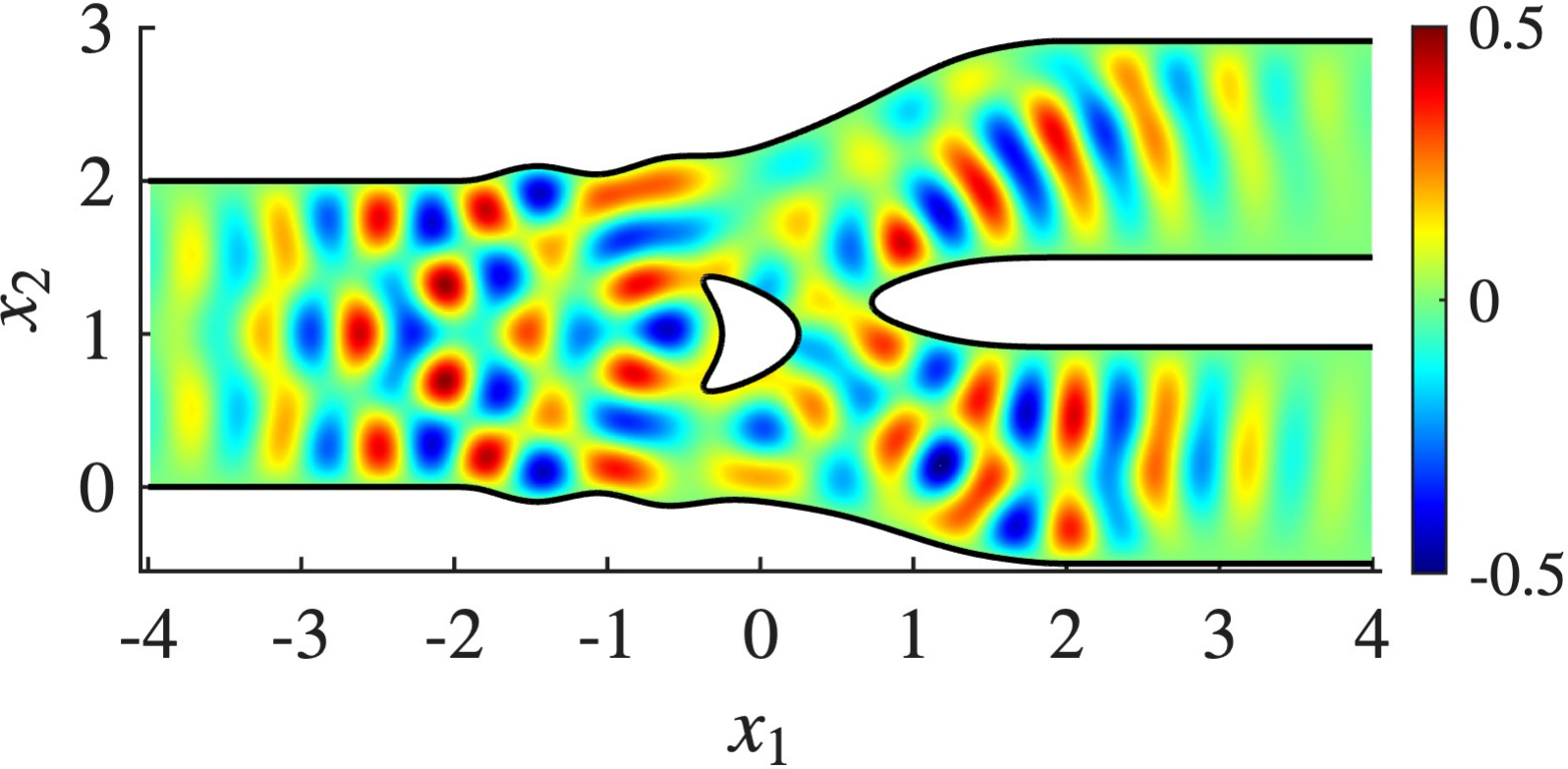}&
\includegraphics[width=0.19\textwidth]{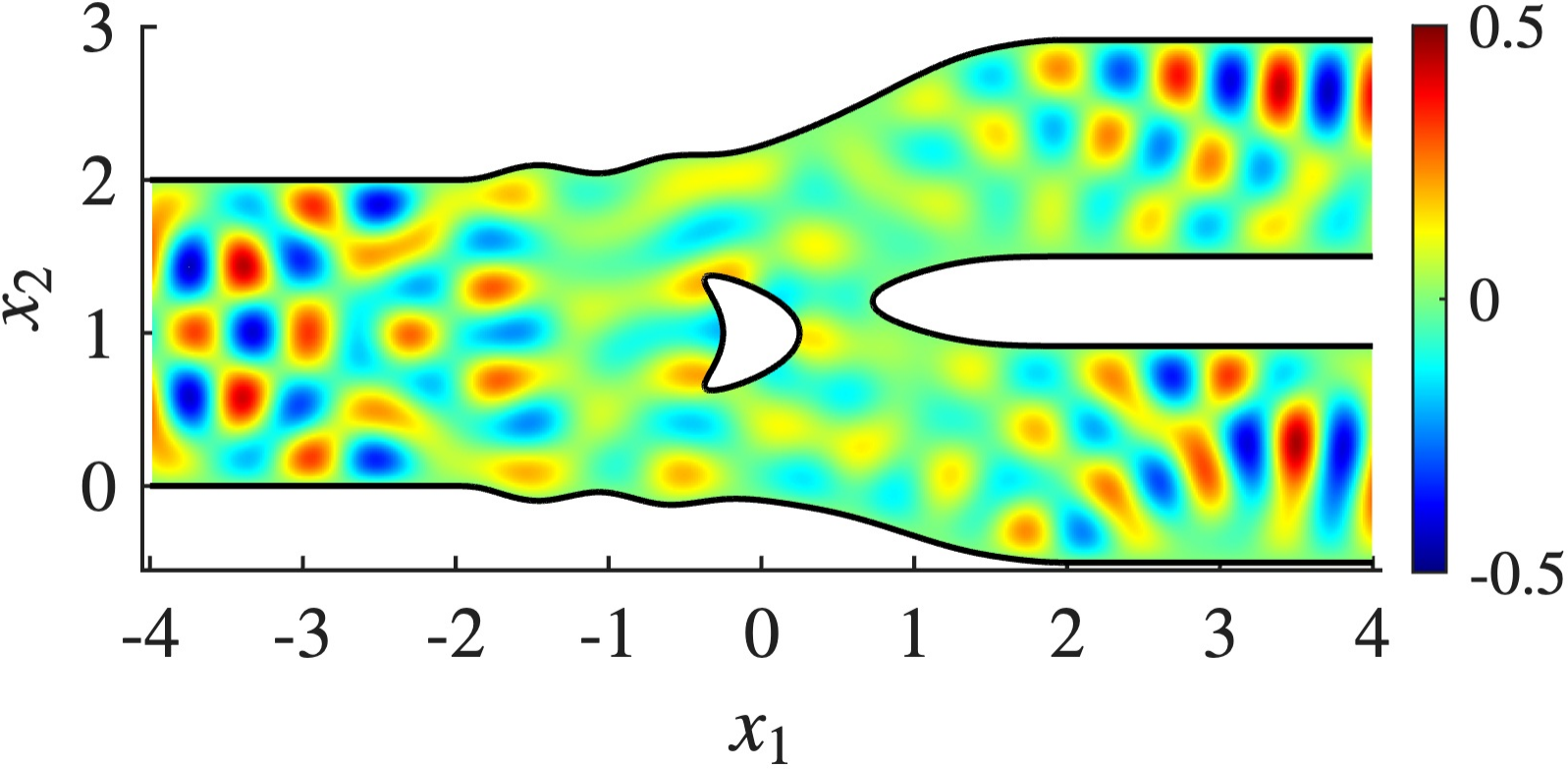}&
\includegraphics[width=0.19\textwidth]{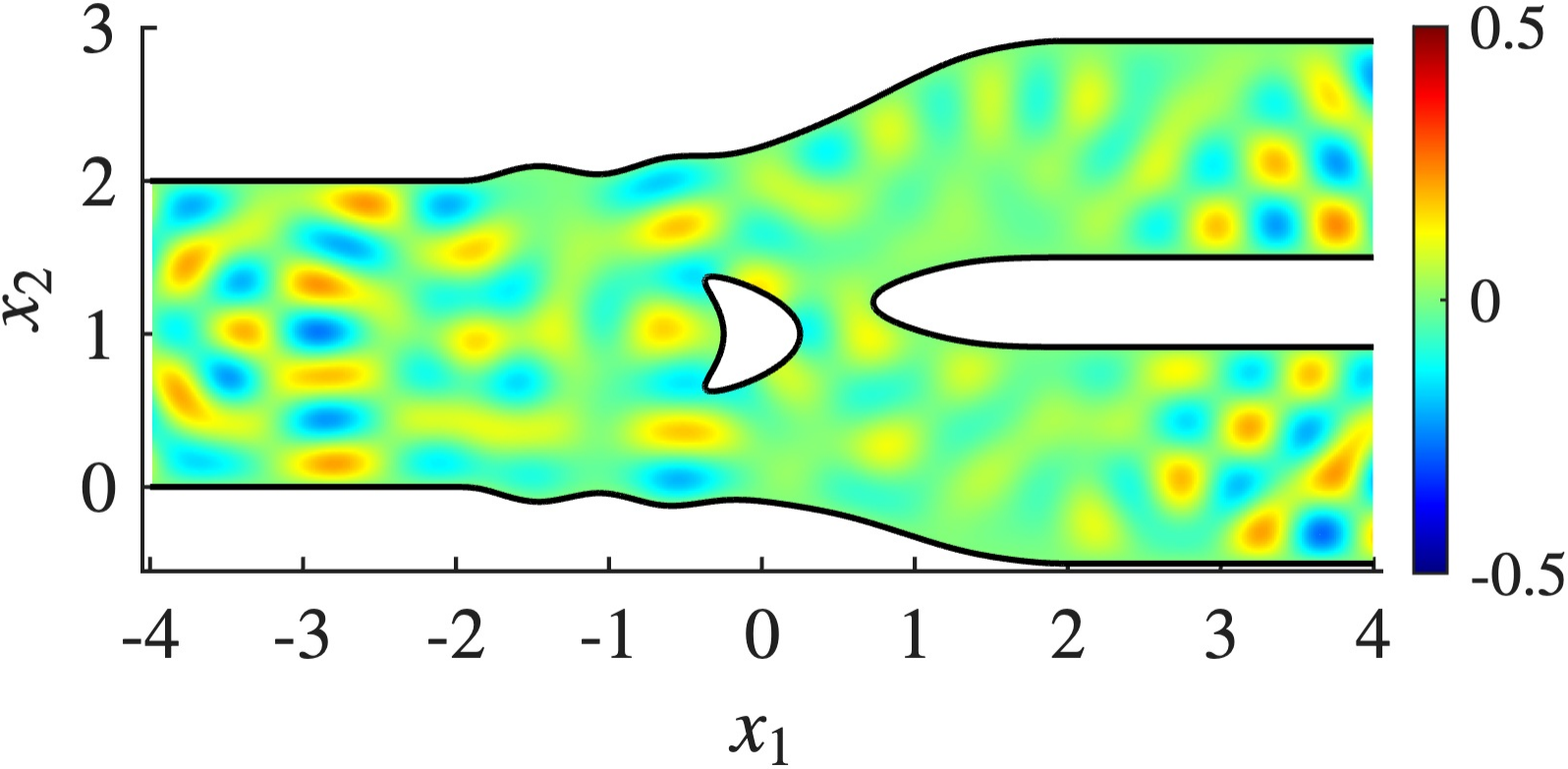}\\[-3mm]
\includegraphics[width=0.19\textwidth]{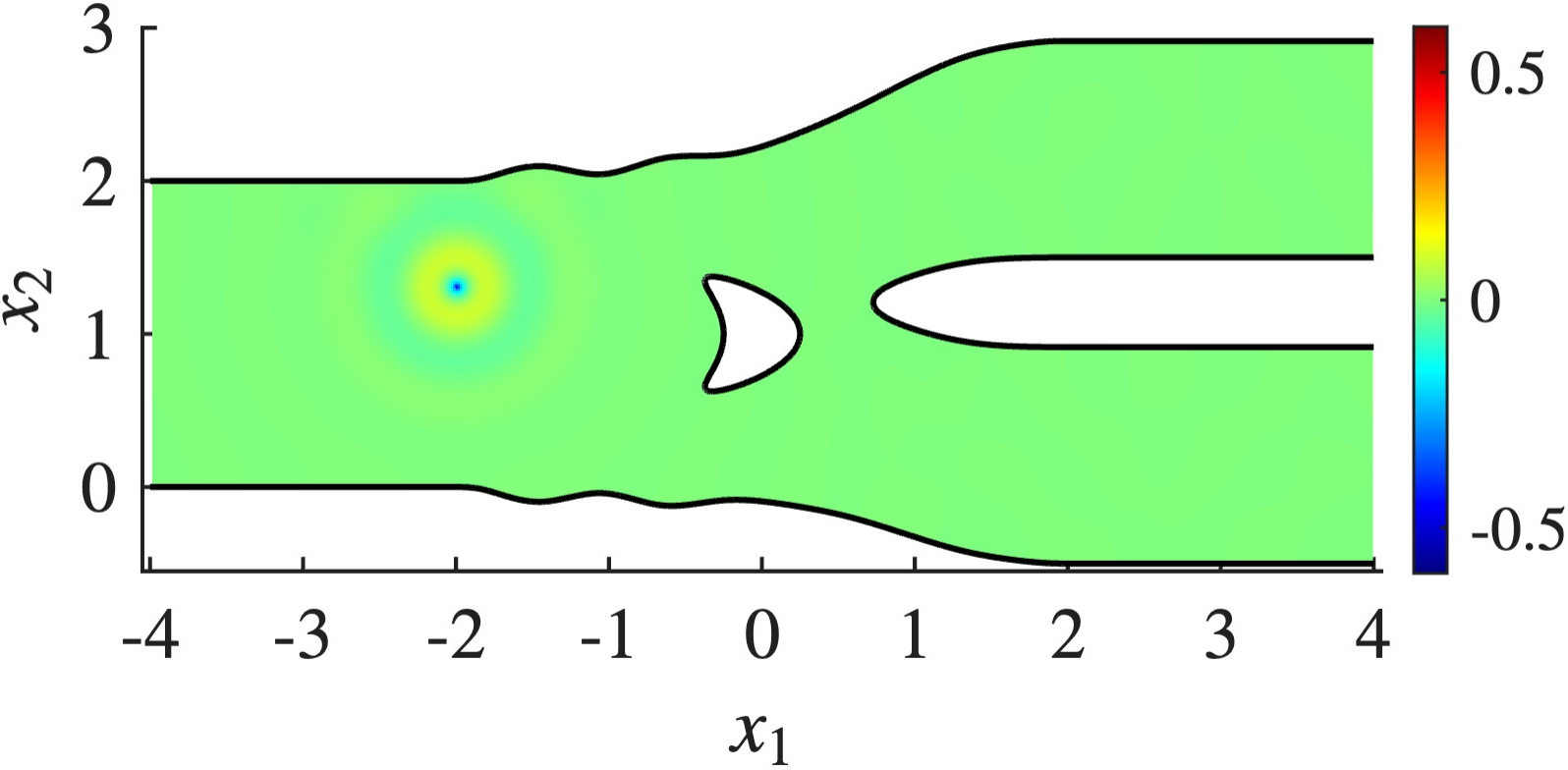}&
\includegraphics[width=0.19\textwidth]{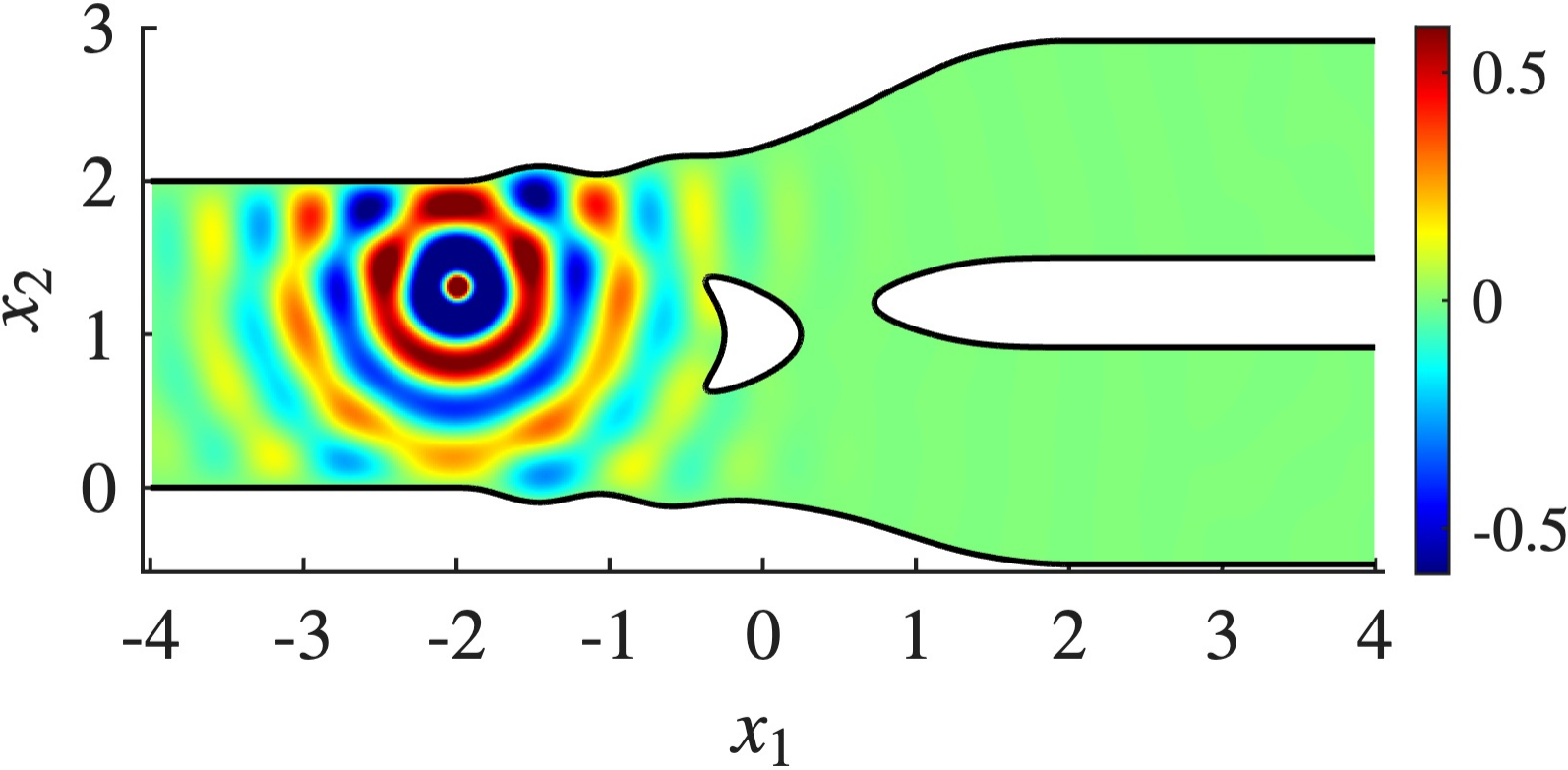}&
\includegraphics[width=0.19\textwidth]{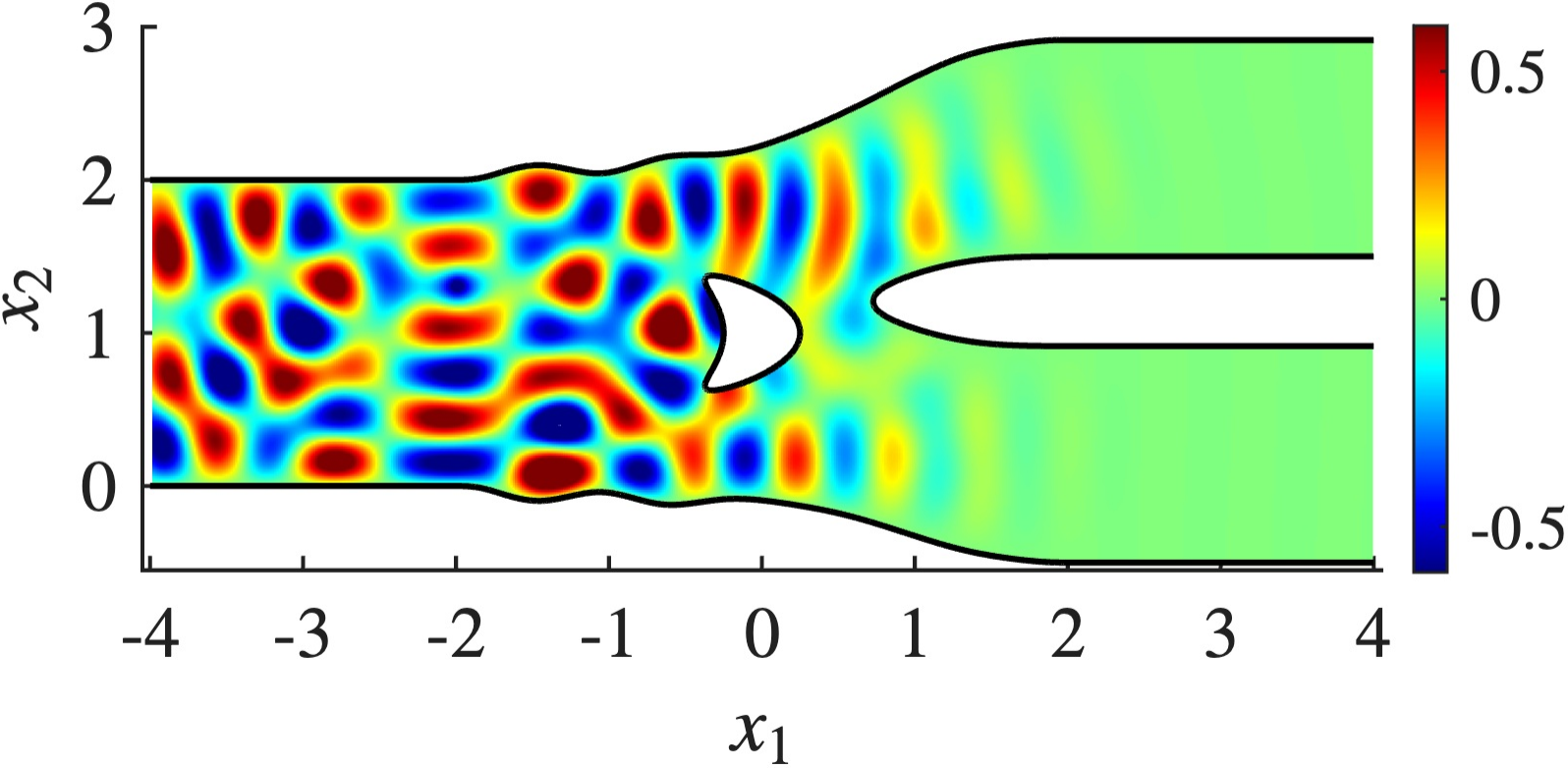}&
\includegraphics[width=0.19\textwidth]{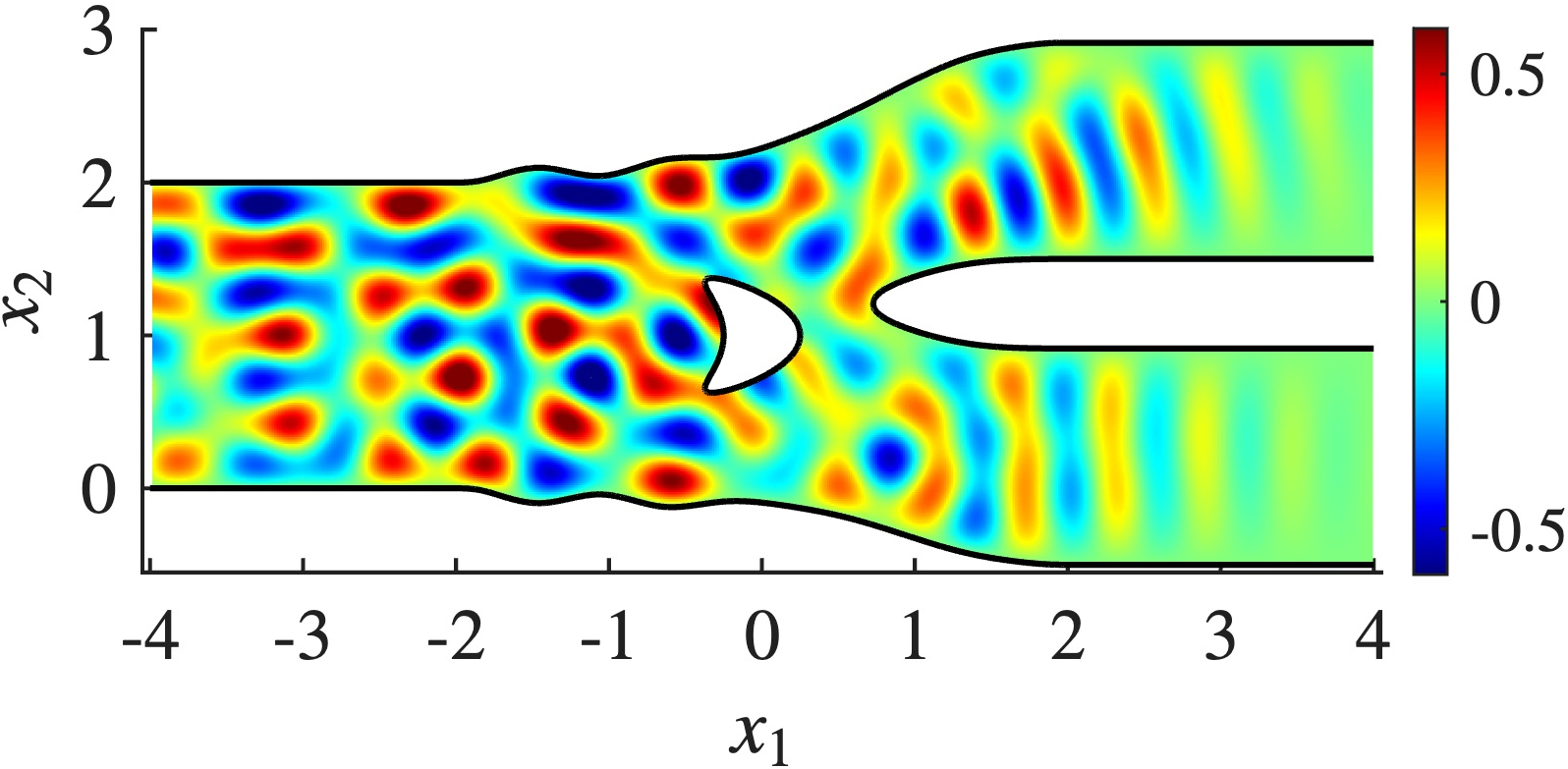}&
\includegraphics[width=0.19\textwidth]{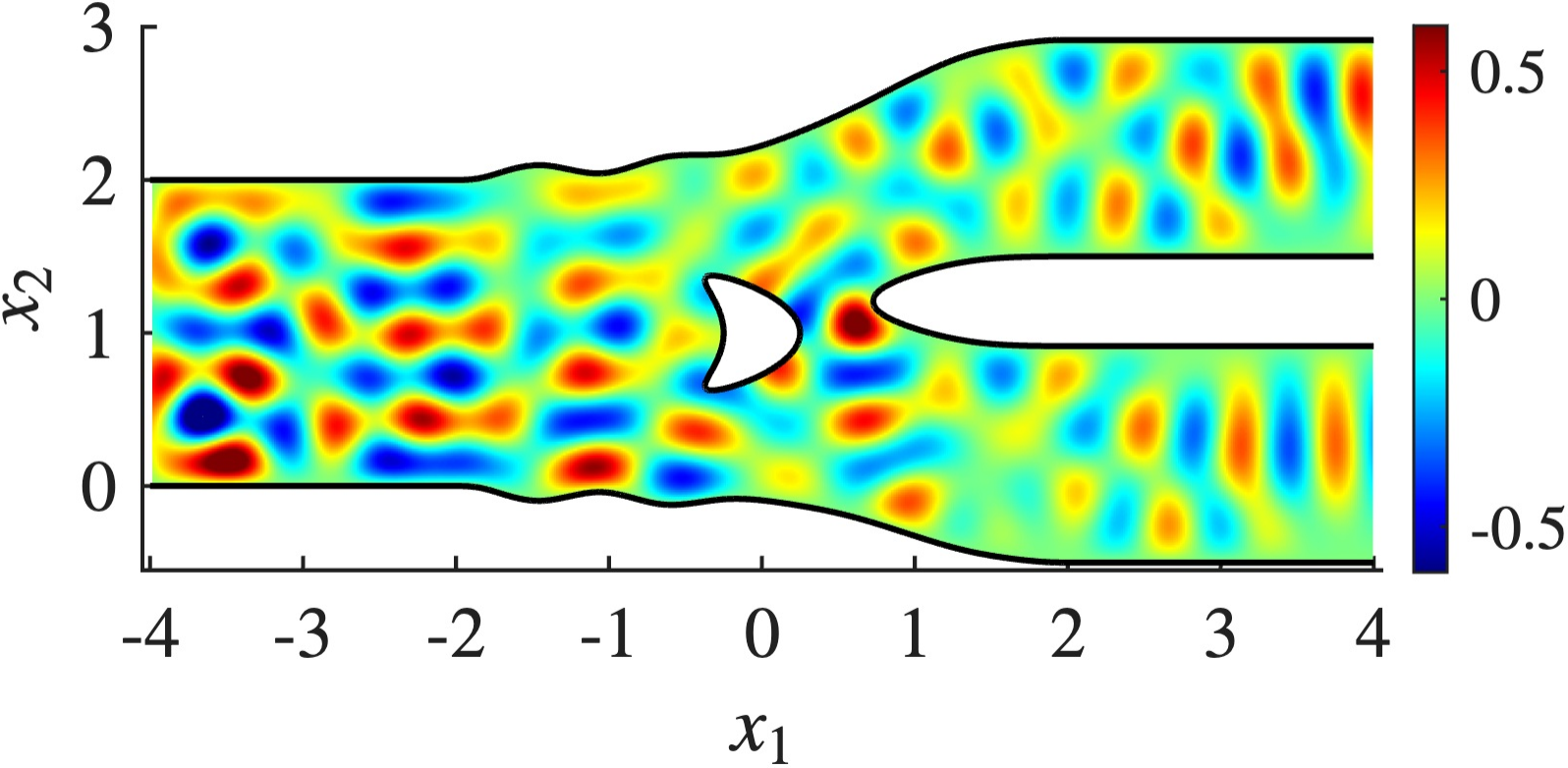}
\end{tabular}
\caption{Example 4. Real parts of the total fields at $t=2,4,6,8,10$ (left to right) for the problem of mode (top) and point-source (bottom) incidences in a parallel three-port sound-soft waveguide.}
\label{fig:Ex3-2-utot}
\end{figure}

{\bf Example 4} ({\it time-domain test for nonuniform waveguides}). In the last example, we study the performance of the method for the propagation of incident Gauss pulses of plane wave or point source located at $\bm z=(-2,1.3)$ in nonuniform waveguides. Considering the sound-hard Y-junction waveguide and choosing $n_0=0$, Figure~\ref{fig:Ex3-3-utot} presents the real parts of the numerical total fields at $t=2,4,6,8,10$. The time traces of the total fields at the observation point $\bm x_{\rm ob}=(-0.5,0.5)$ and the numerical errors $\varepsilon_{\rm TD}(\bm x_{\rm ob},t)$ over the time interval $0\le t\le 10$ are shown in Figure~\ref{fig:Ex3-3-single-pulse-AAA} which again illustrates the high accuracy of the proposed method. The comparison of the numerical errors in Figure~\ref{fig:Ex3-3-single-pulse-AAA}(b,d) for the methods with or without the singularity subtraction technique discussed in Section~\ref{sec:sine-SPFCC} clearly verify that the complex resonances can significantly effect the accuracy of Fourier transform algorithm and show the high accuracy and efficiency of our method. Finally, we set $n_0=1$ and perform the wave propagation in a parallel three-port sound-soft waveguide whose wall structure is optimized in \cite[Figure~1(a)]{SGB19} in Figure~\ref{fig:Ex3-2-utot}. 


\section{Conclusions}
\label{sec:6}
We have developed a frequency/time hybrid RPML-BIE solver for the problems of wave propagation in nonuniform waveguides. A few-mode DtR-correction technique was introduced to improve the exponential convergence of the RPML truncation for the frequency-domain problems which can ensure uniform high accuracy and allows us to set parameters of PML stretching to be independent of the distance between excited frequency and the cut-off frequencies. Then a highly accurate Fourier transform algorithm was designed to generate the time-domain solutions from the frequency-domain solutions by inverse Fourier transform. A sine transform was applied to regularize the square-root singularity near the cut-off frequencies and then an SP-FCC quadrature rule was utilized to treat the reduced moments which contains stationary points. Furthermore, a novel strategy of complex resonance extraction based on AAA rational approximation and singularity subtraction was proposed to avoid the influence of the complex resonances closing to the real axis. Several numerical examples were presented to illustrate the efficiency and high accuracy of the proposed methodology. Extensions to more complicated three-dimensional waveguides and considering additionally the low-frequency contribution and long-time simulation remain natural directions for further work.

\section*{Acknowledgments}

TY greatly acknowledge the support from the National Key R\&D Program of China (2024YFA1016000), Beijing Natural Science Foundation (JR26005), the Strategic Priority Research Program of the Chinese Academy of Sciences (XDB0640000). The research of LX is supported in part by the National Natural Science Foundation of China (Nos. 12431015, 62231016), the Foundation for Innovative Research Groups of the Natural Science Foundation of Sichuan Province (No. 2025NSFTD0004), and the Fundamental and Interdisciplinary Disciplines Breakthrough Plan of the Ministry of Education of China (No. JYB2025XDXM102).

\bibliographystyle{siamplain}
\bibliography{Manuscript}

\end{document}